\documentclass[review,17pt,reqno]{Elsevier}
\usepackage{amssymb,amsmath}
\usepackage{cases}
\usepackage{exscale}
\usepackage{relsize}
\usepackage{subfigure}
\usepackage{natbib}
\usepackage[colorlinks=true,linkcolor=blue,citecolor=blue,urlcolor=blue,]{hyperref}
\usepackage{caption}
\usepackage{booktabs}
\usepackage{makecell}
\usepackage{graphicx}
\usepackage{pythonhighlight}
\usepackage{array, threeparttable}
\usepackage{algorithmicx}
\usepackage[ruled]{algorithm2e}    
\usepackage{algpseudocode}
\usepackage{mathrsfs} 
\usepackage{arydshln}
\usepackage{multirow}
\usepackage{graphicx} 

\usepackage{epsfig}                          
\usepackage{epstopdf}                        
\usepackage{appendix}

\begin{document}
\newtheorem{mytheorem}{Theorem}[section]
\newtheorem{myproposition}{Proposition}[section]
\newtheorem{mylemma}{Lemma}[section]
\newtheorem{mycorollary}[mytheorem]{Corollary}
\newtheorem{myprop}[myproposition]{Proposition}
\newtheorem{myremark}{Remark}[section]
\renewcommand{\baselinestretch}{1.1}
\newcommand{\figref}[1]{Fig.~\ref{#1}}
\makeatletter
\newcommand\figcaption{\def\@captype{figure}\caption}
\newcommand\tabcaption{\def\@captype{table}\caption}
\makeatother
\begin{frontmatter}
\title{Primal-dual splitting methods for phase-field surfactant model with moving contact lines}

\author[label1]{Wei Wu}
\ead{wuwei837037@163.com}

\author[label2]{Zhen Zhang}
\ead{zhangz@sustech.edu.cn}

\author[label1]{Chaozhen Wei\corref{cor1}}
\ead{cwei4@uestc.edu.cn}\cortext[cor1]{Corresponding author.}

\address[label1]{School of Mathematical Sciences, University of Electronic Science and Technology of China, Chengdu, Sichuan 611731, China}


\address[label2]{Department of Mathematics, National Center for Applied Mathematics (Shenzhen), Southern University of Science and Technology (SUSTech), Shenzhen 518055, China}

\begin{abstract}
Surfactants have important effects on the dynamics of droplets on solid surfaces, which has inspired many industrial applications. Phase-field surfactant model with moving contact lines (PFS-MCL) has been employed to investigate the complex droplet dynamics with surfactants, while its numerical simulation remains challenging due to the coupling of gradient flows with respect to transport distances involving nonlinear and degenerate mobilities. We propose a novel structure-preserving variational scheme for PFS-MCL model with the dynamic boundary condition based on the minimizing movement scheme and optimal transport theory for Wasserstein gradient flows. The proposed scheme consists of a series of convex minimization problems and can be efficiently solved by our proposed primal-dual splitting method and its accelerated versions. By respecting the underlying PDE's variational structure with respect to the transport distance, the proposed scheme is proved to inherits the desirable properties including original energy dissipation, bound-preserving, and mass conservation. Through a suite of numerical simulations, we validate the performance of the proposed scheme and investigate the effects of surfactants on the droplet dynamics. 
\end{abstract}

\end{frontmatter}
\section{Introduction}
The wetting and dewetting dynamics of fluids on a solid substrate is ubiquitous in nature and has been employed in many industrial applications, ranging from spray coating, cleaning detergents, emulsifiers, dispersants, microfluids and ink-jet printing \cite{Probstein2005Physic,Eggleton2001Tip,Branger2002Exogenous,Baret2012Surfactants,Kommeren2018Combining}. The essence of the modulation of droplet dynamics is to control the wetting property of the solid surface. Emerging as the intersection of fluid-fluid interface and solid surface, the contact line plays an important role in such a physical process. And the dynamic process is further complicated when surfactants come into play.
Surfactants, as surface active agents, can alter the contact angle by reducing the interfacial tension between two fluids, and hence have attracted intense attention. 
Multiphase flows with moving contact lines, when coupled with the transport dynamics of surfactants, exhibits highly complex phenomena, which poses significant challenges for experimental investigation. Therefore, modeling and simulations serve as crucial alternatives for studying the complex droplet dynamics in the presence of surfactants.

The multiphase fluid-surfactant models have been intensely studied over the past few decades. Meanwhile, numerous numerical methods have been developed for their simulations. These include finite difference method \cite{GuBinary2014FluidSurfactant,Yang2017BinaryFluidSurfactant}, finite element method \cite{Ganesan2012FETwoPhaseFlows}, finite volume method \cite{Yon1998FVSurfactant,James2004FVSurfactant}, level set method \cite{Xu2006LevelSetSurfactant,Xu2012LevelSetSurfactant}, immersed boundary method \cite{Lai2008ImmersedSurfactant,Wu2024ImmersedMCL} and phase-field methods \cite{Liu2010PhaseField,Zhu2019NumericalPhaseField}. Given the capability of handling complex topological phase changes and the favored energetic variational structure, phase-field modeling has been successfully employed for simulating the multiphase flow with surfactants and recently generalized to account for  contact line dynamics \cite{Zhu2019ThermodynamicallyMCL}. 
The phase-field surfactant model with moving contact lines (PFS-MCL) employs two order parameters (or phase-field variables) to describe the volume fraction of two fluid phases and the concentration of surfactants, respectively. In addition to the classical Ginzburg-Landau free energy with double-well potential for the fluid-fluid interfacial energy and the Flory-Huggins free energy with logarithmic potential for the entropy of mixing surfactants, the model also incorporates a nonlinear coupling energy for the adsorption of surfactants on the fluid-fluid interface, as well as the wall free energy for the moving contact line. 

The PFS-MCL model can be viewed as a nonlinearly coupled system of multiple gradient flows with respect to different metrics for the dynamics of phase variable, surfactant concentration and moving contact line, rendering its numerical simulation a difficult task. Different from classical $\mathcal{L}^{2}$ or $\mathcal{H}^{-1}$ gradient flows with constant mobilities where the stiffness arises only from the nonlinear energy terms, the stability issue of PFS-MCL model is even worsen due to the presence of nonlinear, degenerate mobility for surfactant dynamics. Moreover, the degenerate mobility also present a global constraint on the bounds of solutions independent of the specific energy potential \cite{Elliott1996CahnHilliard}, bringing new challenges for the numerical simulations. Energy dissipative or energy stable schemes are favored to tackle the stability issue \cite{Tang2020EfficientPhaseField}, and a number of designing techniques have been developed. These numerical techniques include convex splitting method, stablized semi-implicit (SSI) method, invariant energy quadratization (IEQ) method and scalar auxiliary variable (SAV) method \cite{Eyre1998UnconditionallyGradient,Xu2006SSIEpitaxialGrowth,Yang2016IEQLinearFirstSecond,Yang2020IEQAllenCahn,Yang2021SIEQCahnHilliard,Shen2018SAVGradientFlows}. They have also been employed in the simulations of moving contact line problems \cite{Gao2012GradientMCL,Zhu2020PhaseFieldMCL,Kang2020SAVMultiple}. It is worth noting that these energy stable schemes may only possess the dissipation structure with respect to modified energies or dissipation functions rather than the original ones. Moreover, these energy stable schemes do not necessarily have the feature of bound-preserving. Without bound-preserving, the simulations may lead to unphysical solutions or even get stuck due to the singularity of the logarithmic Flory-Huggins potential. Several techniques involving cut-off method, flux limiting scheme and Lagrange multiplier \cite{Lu2013CutoffParabolic,Frank2020BoundPreservingCH,Cheng2022NewLagrangeMultiplier} were developed to establish bound-preserving schemes. Recently, a convex splitting based numerical method was proposed for the PFS-MCL model, leading to provable first-order unconditionally (original) energy stable and bound-preserving schemes \cite{Wang2022PFSMCLModel}. These schemes have demonstrated their success in the simulations of contact angle hysteresis and droplet impacting dynamics in the presence of surfactants \cite{Wang2024PhaseFieldNS}. 
However, due to the nonlinearity in the convex minimization step, the resulting schemes require solving a large coupled nonlinear system at each time step by Newton iteration with a damped step size, which could be time-consuming. Moreover, the original transport distance induced by the nonlinear mobility in the minimization was not solved directly but only approximated by an extrapolated metric from the solution in the previous time step, rendering a modified dissipation structure. 

In this paper, we develop a novel structure-preserving method based on the minimizing movement scheme that has been recently proposed for gradient flows with non-constant mobilities \cite{Carrillo2022PrimalDual,Carrillo2024StructurePD}, which are often referred as Wasserstein gradient flows. This approach is based on the so-called {\it JKO} scheme \cite{Jordan1998Variational,Lisini2012CahnHilliard} that is a fully implicit variational scheme for gradient flows with respect to transport distances (also referred as Wasserstein distances). By leveraging the underlying PDE's variational structure, this scheme is unconditionally energy stable and naturally inherits the favored properties of original energy dissipation, bound-preserving and mass conservation. The crucial improvements of this newly proposed approach over the one in \cite{Wang2022PFSMCLModel} lie in two essential aspects. First, the fully discrete scheme is a result of variation-then-discretization, thus is actually a consistent discretization of the variational scheme for the whole PFS-MCL system that accounts for the dynamics of fluids, surfactants and moving contact lines. Therefore, the energy dissipation property is naturally inherited from the variational structure itself, avoiding the specific and complicated design of convex splitting in the discretization of PDEs. In addition, the preservation of the original variational structure also bring benefits to the development of adaptive time stepping schemes based on the original dissipation functions. Secondly, instead of using an approximated metric in \cite{Wang2022PFSMCLModel}, we preserve the original transport distance and utilize its dynamic characterization \cite{Benamou2000Computational,Dolbeault2009NewTransport} to reformulate the non-convex minimization problem into a convex minimization with linear constraints, which can be solved by our proposed primal-dual splitting algorithms and its accelerated versions. By treating the nonlinear degenerate mobility implicitly in the dynamic formulation of transport distances and employing the proximal gradient splitting method, our proposed scheme guarantees the global bounds of solutions in both theoretical and practical perspectives, independent of the specific energy functional forms. This approach can be easily employed for gradient flows with degenerate mobilities and various energy functionals. 

The rest of the paper is organized as follows. In Sec.~\ref{sec:2}, we briefly review the phase-field surfactant model with moving contact lines and the governing equations of the dynamic system. In Sec.~\ref{sec:3}, we first introduce the JKO scheme for Wasserstein gradient flows and then propose the variational formulation of the PFS-MCL model with dynamic boundary conditions. We show that our fully discrete variational scheme has the desired structure-preserving properties. In Sec.~\ref{sec:4}, we mainly discuss the primal-dual splitting algorithm and its accelerated versions. In Sec.~\ref{sec:5}, we show the performance of our numerical method and employ it to investigate the effects of surfactants on the droplet dynamics. Finally, we conclude with a summary and outlook.

\section{Introduction of the PFS-MCL Model}\label{sec:2}
\begin{figure}[h]
	\centering			
	\includegraphics[width=0.5\textwidth]{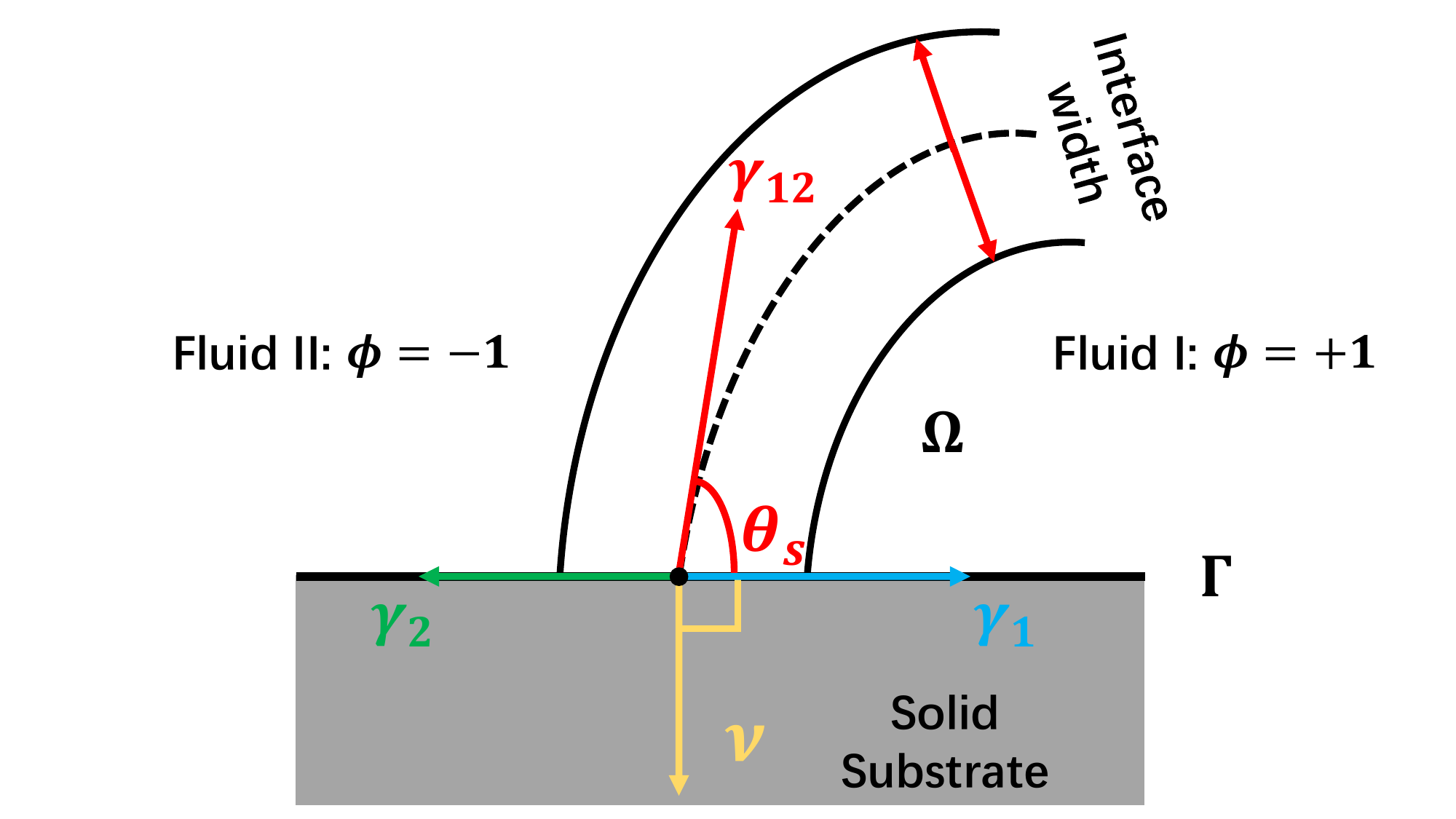}
	\caption{Phase-field description of the diffused interface between two fluids ($\phi=\pm 1$) in contact with a solid substrate $\Gamma$ with a contact angle $\theta_s$. $\gamma_{1}$ and $\gamma_{2}$ are the fluid-solid substrate interfacial tensions with Fluid I (with $\phi=+1$) and Fluid II (with $\phi=-1$), respectively, and $\gamma_{12}$ is the interfacial tension for the fluid-fluid interface. Young’s angle $\theta_{s}$ satisfies the Young–Dupré equation $\gamma_{12}\mathrm{cos}(\theta_{s})+\gamma_{1}=\gamma_{2}$.}
	\label{fig:TwoFluidsAngle}
\end{figure}

We firstly give a brief introduction of the dimensionless PFS-MCL model \cite{Zhu2019ThermodynamicallyMCL,Wang2022PFSMCLModel}. In the phase-field model, the order parameter $\phi$ takes $\pm 1$ in two fluids, between which the sharp interface is diffused and represented by a smooth transition between $\pm 1$ (see Fig.~\ref{fig:TwoFluidsAngle}); the concentration variable $\psi\in (0,1)$ describes the distribution of the surfactant.  

The total free energy of the PFS-MCL model consists of the energies associated with the interface, surfactant and their interactions, and the fluid-solid substrate,
\begin{eqnarray}\label{eq:totalEne}
	\begin{aligned}
    \mathcal{E}(\psi,\phi)=\mathcal{F}_{GL}(\phi)+\mathcal{F}_{sur}(\psi)+\mathcal{F}_{ad}(\phi,\psi)+\mathcal{F}_{\omega f}(\phi).
	\end{aligned}
\end{eqnarray}

In the above expression, $\mathcal{F}_{GL}$ is the Ginzburg-Landau free energy with the double-well potential $H(\phi)$ representing the energy of the diffused fluid-fluid interface
\begin{eqnarray}\label{eq:FGL}
	\begin{aligned}
		\mathcal{F}_{GL}(\phi)=\int_{\Omega}\Big(\dfrac{\mathrm{Cn}^{2}}{2}\left|\nabla \phi\right|^{2}+H(\phi)\Big)\mathrm{d}x:=\int_{\Omega}\Big(\dfrac{\mathrm{Cn}^{2}}{2}\left|\nabla \phi\right|^{2}+\dfrac{1}{4}(\phi^{2}-1)^{2}\Big)\mathrm{d}x,
	\end{aligned}
\end{eqnarray}
where $\mathrm{Cn}$ is a small parameter related to the dimensionless interface thickness. 

The term $\mathcal{F}_{sur}$ is the Flory-Huggins free energy of the mixing entropy in the binary surfactant-fluid system 
\begin{eqnarray}\label{eq:Fsur}
	\begin{aligned}
		\mathcal{F}_{sur}(\psi)=\int_{\Omega}\mathrm{Pi}G(\psi)\mathrm{d}x:=\int_{\Omega}\Big(\mathrm{Pi}(\psi\mathrm{log}\psi+(1-\psi)\mathrm{log}(1-\psi))\Big)\mathrm{d}x,
	\end{aligned}
\end{eqnarray}
where $\mathrm{Pi}$ is a constant temperature-dependent surfactant diffusion rate. 

The interaction energy between the surfactant and the fluids $\mathcal{F}_{ad}$ mainly models the adsorption of surfactant on the interface
\begin{eqnarray}\label{eq:Fad}
	\begin{aligned}
		\mathcal{F}_{ad}(\phi,\psi)=\int_{\Omega}\Big(\dfrac{1}{2\mathrm{Ex}}P_{1}(\phi,\psi)-\dfrac{1}{4}P_{2}(\phi,\psi)\Big)\mathrm{d}x:=\int_{\Omega}\Big(\dfrac{1}{2\mathrm{Ex}}\psi\phi^{2}-\dfrac{1}{4}\psi(\phi^{2}-1)^{2}\Big)\mathrm{d}x,
	\end{aligned}
\end{eqnarray}
where the term with $P_1=\psi\phi^{2}$ penalizes the free surfactant dissolved in the bulk phases (where $\phi=\pm 1$) and $\mathrm{Ex}$ represents the bulk solubility of surfactant (or equivalently, the inverse of adsorption rate on interface), while the negative term with $P_2=\psi(\phi^{2}-1)^{2}$ represents the inclination of the surfactant to adsorb on the interface. Other expressions for $\mathcal{F}_{ad}$ are also discussed in \cite{Engblom2013TwoPhase}.

The last term $\mathcal{F}_{\omega f}$ represents the wall free energy of the fluid wetting on the solid substrate \cite{Qian2003Hydrodynamics}
\begin{equation}\label{eq:Fwf}
    \mathcal{F}_{\omega f}(\phi)=\mathrm{Cn}\int_{\partial \Omega}\gamma_{\omega f}(\phi)\mathrm{d}s,
\end{equation}
where the wetting energy potential $\gamma_{\omega f}$ is a smooth interpolation between the fluid-solid interfacial tensions $\gamma_{1}$ and $\gamma_{2}$ defined along the solid substrate $\Gamma$ and here we extend it to the boundary of a finite region $\partial \Omega$
\begin{eqnarray}\label{eq:gamma_wf}
	\begin{aligned}
		\gamma_{\omega f}(\phi)
		=\left\{
		\begin{aligned}
			&-\dfrac{\sqrt{2}}{3}\mathrm{cos}(\theta_{s})\mathrm{sin}(\dfrac{\pi\phi}{2})+\dfrac{\gamma_{1}+\gamma_{2}}{2}\quad &&\text{at $\Gamma$},\\
			&0 \quad&&\text{at $\partial\Omega/\Gamma$},
		\end{aligned}
		\right.
	\end{aligned}
\end{eqnarray}
where $\theta_{s}$ is the prescribed static contact angle between the fluid-fluid interface and the substrate, determined by  the balance of local interfacial tensions as $\mathrm{cos}(\theta_{s})=3\sqrt{2}(\gamma_{2}-\gamma_{1})/4$. Other commonly used forms of $\gamma_{\omega f}$ are also discussed in \cite{Zhu2019ThermodynamicallyMCL, Xuxianmin2018Sharp, Huang2013Wetting}.

The coupled dynamics of the surfactant-fluid system can be written in the form of Wasserstein-like gradient flows with respect to specific metrics induced by the corresponding mobilities
\begin{align}\label{eq:govern_eqn}
&\dfrac{\partial}{\partial t}
\begin{bmatrix}
\phi \\ \psi
\end{bmatrix} 
= \nabla\cdot \Bigg(
\begin{bmatrix}
M_{\phi} & 0 \\
0 & M_{\psi}
\end{bmatrix}
\begin{bmatrix}
\nabla \dfrac{\delta\mathcal{E}}{\delta \phi} \\ \nabla \dfrac{\delta\mathcal{E}}{\delta \psi}
\end{bmatrix} \Bigg), \\ 
&\dfrac{\delta\mathcal{E}}{\delta\phi}=-\mathrm{Cn}^{2}\Delta\phi+\phi^{3}-\phi+\dfrac{1}{\mathrm{Ex}}\psi\phi-\psi\phi(\phi^{2}-1),\\ 
&\dfrac{\delta\mathcal{E}}{\delta\psi}=\mathrm{Pi}\ \mathrm{log}(\dfrac{\psi}{1-\psi})+\dfrac{1}{2\mathrm{Ex}}\phi^{2}-\dfrac{1}{4}(\phi^{2}-1)^{2},
\end{align}
where $M_{\psi}=\psi(1-\psi)/\mathrm{Pe}_{\psi}$ is a nonnegative, nonlinear, degenerate mobility leading to the Fickian diffusion of surfactants, and  $M_{\phi}=1/\mathrm{Pe}_{\phi}$ is a constant mobility for phase separation, where $\mathrm{Pe}_{\phi}$ and $\mathrm{Pe}_{\psi}$ are the P\'eclet numbers representing the magnitude of mobilities. 

The boundary conditions are a combination of no-flux condition that guarantees the mass conservation of $\phi$ and $\psi$
\begin{align}\label{eq:noflux_bc}
    \nabla \dfrac{\delta\mathcal{E}}{\delta \phi}\cdot\boldsymbol{\nu}=0,\quad \nabla \dfrac{\delta\mathcal{E}}{\delta \psi} \cdot\boldsymbol{\nu}=0,\quad\text{on $\partial\Omega$},
\end{align}
and contact angle conditions 
at the boundary
\begin{align}\label{eq:natural_bc} 
    &\nabla\phi\cdot\boldsymbol{\nu}=0\quad&&\text{on $\partial\Omega\setminus\Gamma$}, \\
    \dfrac{\partial\phi}{\partial t}=-\dfrac{1}{\mathrm{Pe}_{s}}L_{\phi},\quad 
    &L_{\phi}=\mathrm{Cn}^2\nabla\phi\cdot\boldsymbol{\nu}+\mathrm{Cn}\gamma^{'}_{\omega f}(\phi)\quad&&\text{on $\Gamma$},\label{eq:dynamic_bc}
\end{align}
where $\boldsymbol{\nu}$ is the unit outer normal vector, and $\mathrm{Pe}_{s}$ controls the mobility at the boundary. The natural boundary condition \eqref{eq:natural_bc} and the dynamic boundary condition \eqref{eq:dynamic_bc} are imposed to guarantee zero energy dissipation on the non-substrate boundaries $\partial \Omega \setminus\Gamma$ and positive energy dissipation on the solid substrate $\Gamma$ respectively.

The PFS-MCL model \eqref{eq:govern_eqn}-\eqref{eq:dynamic_bc} is actually a coupled system of $\mathcal{H}^{-1}$ gradient flow for $\phi$, the generalized Wasserstein ($\mathcal{W}^2_m$) gradient flow for $\psi$, and $\mathcal{L}^{2}$ gradient flow for the dynamic wetting boundary condition, which is equipped with the  following energy dissipation law 
\begin{eqnarray}\label{eq:dE_dt}
	\begin{aligned}
		\dfrac{d\mathcal{E}}{dt}=-M_{\phi}\Big\Vert\nabla\dfrac{\delta\mathcal{E}}{\delta \phi}\Big\Vert^{2}_2-\Big\Vert\sqrt{M_{\psi}}\nabla\dfrac{\delta\mathcal{E}}{\delta \psi}\Big\Vert^{2}_2-\dfrac{1}{\mathrm{Pe}_{s}}\Vert L_{\phi}\Vert^{2}_{2,\Gamma}\leq 0,
	\end{aligned}
\end{eqnarray}
where $\Vert \cdot \Vert_2$ denotes the $\mathcal{L}^{2}$ norm. It is challenging to design an efficient stable numerical method for the above coupled system of gradient flows due to the strong nonlinearity coming from both the nonlinear energy functional and nonlinear mobility, as well as the motion of the contact line.

\section{Variational formulation}\label{sec:3}
\begin{figure}[htbp]
	\centering
	\includegraphics[width=1.0\textwidth]{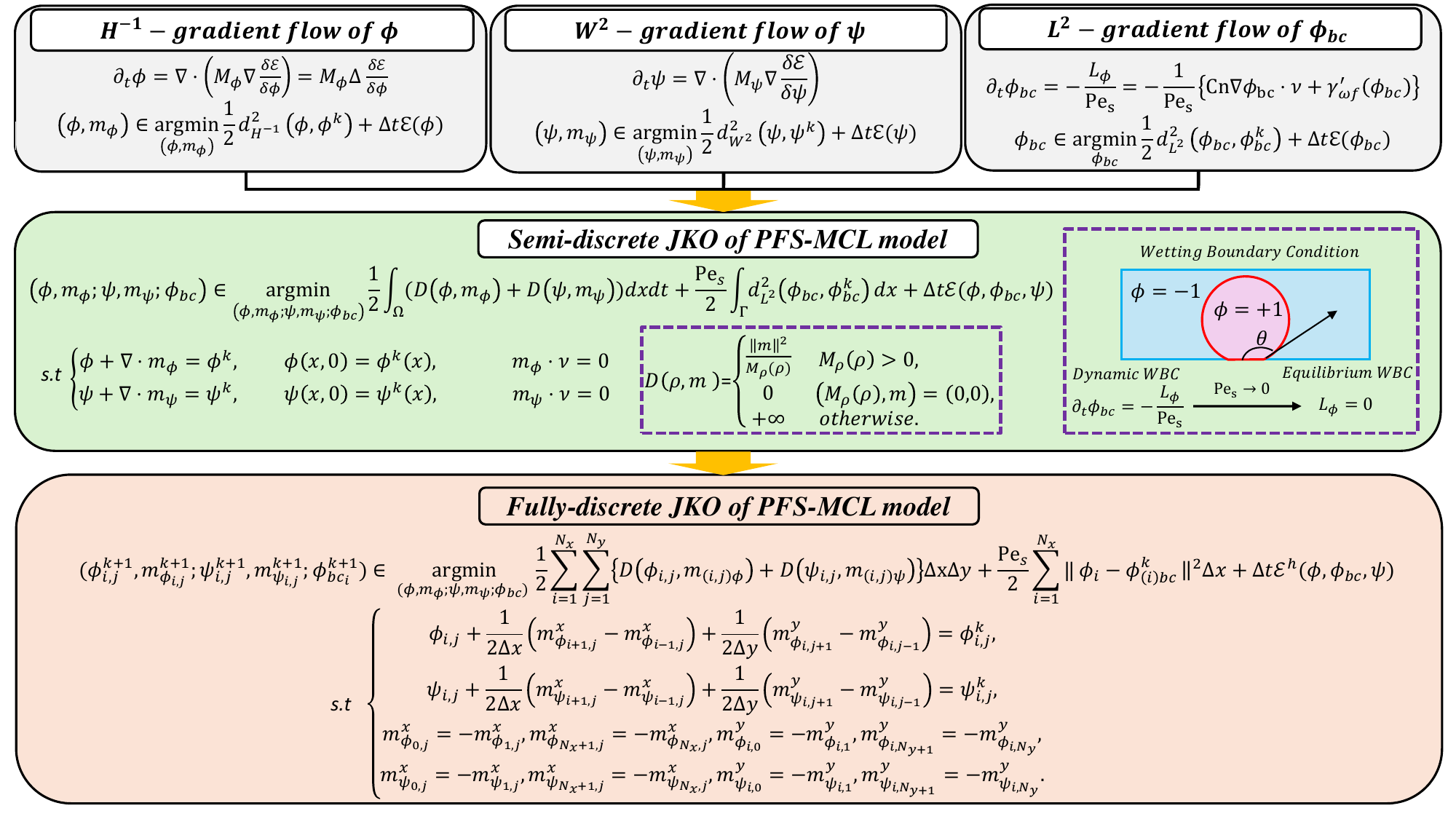}
	\caption{Illustration of the variational approach for the PFS-MCL model.}
	\label{fig:FlowChart}
\end{figure}

In this section, we introduce a novel approach based on the underlying variational structure and optimal transport theory to propose a stable numerical scheme that preserves the original energy dissipation structure, global bounds of solution and mass conservation law at the fully discrete level. The outline of this variational approach is summarized in Fig.~\ref{fig:FlowChart}.

\subsection{Semi-discrete JKO scheme for PFS-MCL model}\label{sec:semi_discrete}
Denoting the vector variable $\boldsymbol{\rho}=(\phi,\psi)^\text{T}$, the velocity $\textbf{v}_\rho=(v_{\phi},v_{\psi})^\text{T}\in \mathbb{R}^{2 \times d}$, and the mobility $\textbf{M}_{\rho}=\text{diag}(M_{\phi},M_{\psi})$, the coupled system \eqref{eq:govern_eqn} can be rewritten as the Wasserstein gradient flow of multiple variables with (diagonal) matrix-valued mobility as follows 
\begin{equation}\label{eq:govern_matrix}
\dfrac{\partial}{\partial t} \boldsymbol{\rho}
= -\nabla\cdot \big(\textbf{M}_{\rho} \textbf{v}_\rho\big), \quad \textbf{v}_\rho =- (\nabla \dfrac{\delta \mathcal{E}}{\delta \phi},  \nabla\dfrac{\delta \mathcal{E}}{\delta \psi})^\text{T}.
\end{equation}

Following the seminal work \cite{Jordan1998Variational} and its generalization \cite{Lisini2012CahnHilliard}, we consider the variational characterization of the Wasserstein gradient flows based on the minimizing movement scheme \cite{DeGiorgi1993} (i.e., the celebrated {\it JKO} scheme). The JKO scheme generates a temporal discrete sequence $\{\boldsymbol{\rho}^{k}(x)\}$ to approximate the solutions at $t_k$ by solving the following minimization problems
\begin{equation}\label{eq:JKO_vec}
    \boldsymbol{\rho}^{k+1}(x)\in\mathop{\text{arg min}}\limits_{\boldsymbol{\rho}}\dfrac{1}{2}d^{2}_{m}(\boldsymbol{\rho}(x),\boldsymbol{\rho}^{k}(x))+\Delta t\mathcal{E}(\boldsymbol{\rho}(x)),
\end{equation}
where the existence of the minimizer and its convergence to the solutions of the continues PDE, with first-order accuracy, as the time step $\Delta t\rightarrow0$ has been proved in \cite{Lisini2012CahnHilliard,Jordan1998Variational,Carrillo2010NonlinearMobility}. In particular, the generalized Wasserstein distance induced by the matrix mobility $\textbf{M}_{\rho}$ is defined in a dynamic formulation \cite{Benamou2000Computational,Dolbeault2009NewTransport} 
\begin{eqnarray}\label{eq:dw}
	\begin{aligned}
		d^{2}_{m}(\boldsymbol{\rho}_{0},\boldsymbol{\rho}_{1})
    	&=\mathop{\mathrm{inf}}\limits_{(\boldsymbol{\rho},\textbf{v}_{\rho})}\int_{0}^{1}\int_{\Omega}
        \textbf{v}_\rho:\big(\textbf{M}_{\rho} \textbf{v}_{\rho}\big)\mathrm{d}x\mathrm{d}s, \\
        &=\mathop{\mathrm{inf}}\limits_{(\phi,v_{\phi};\psi,v_{\psi})}\int_{0}^{1}\int_{\Omega}\Big(M_{\phi}|v_{\phi}|^{2}+M_{\psi}|v_{\psi}|^{2}\Big)\mathrm{d}x\mathrm{d}s, \\
        &=\mathop{\mathrm{inf}}\limits_{(\phi,m_{\phi};\psi,m_{\psi})}\int_{0}^{1}\int_\Omega \Big(\mathcal{D}(\phi, m_{\phi})+\mathcal{D}(\psi, m_{\psi})\Big)\mathrm{d}x\mathrm{d}s, 
	\end{aligned}
\end{eqnarray}
where $\mathbf{A}:\mathbf{B}=\mathrm{Tr}(\mathbf{A}\mathbf{B}^T)$ denotes the Frobenius inner product, 
in the third line we have introduced the momentum variables $\textbf{m}_{\rho}=(m_{\phi},m_{\psi})^\text{T}:=(M_{\phi}v_{\phi},M_{\psi}v_{\psi})^\text{T}=\mathbf{M}_\rho\mathbf{v}_\rho$, and the infimum is taken among $(\boldsymbol{\rho},\textbf{m}_{\rho})$ that satisfy a continuity equation with no-flux boundary conditions, connecting initial and terminal densities
\begin{eqnarray}\label{eq:constraint}
	\begin{aligned}
		\left\{\begin{aligned}
			&\partial_{t}\boldsymbol{\rho}+\nabla\cdot\textbf{m}_{\rho}=\textbf{0},\quad&&(x,s)\in\Omega\times\left[0,1\right],\\		
			&\textbf{m}_{\rho}\cdot\boldsymbol{\nu}=\textbf{0},\quad&&(x,s)\in\partial\Omega\times\left[0,1\right],\\
			&\boldsymbol{\rho}(x,0)=\boldsymbol{\rho}_{0},\quad \boldsymbol{\rho}(x,1)=\boldsymbol{\rho}_{1},\quad&&\text{$x\in \Omega$}.
		\end{aligned}
		\right.	
	\end{aligned}
\end{eqnarray}
By introducing the momentum variable in the dynamic formulation \cite{Benamou2000Computational}, the continuity equation becomes a linear PDE constraint with respect to $(\boldsymbol{\rho},\textbf{m}_{\rho})$. In addition, the distance function $\mathcal{D}$ becomes convex with respect to $(\rho, m_{\rho})$ (without ambiguity, we use the scalar $\rho$ to denote $\phi$ or $\psi$):
\begin{equation} \label{D_function}
\mathcal{D}(\rho,m_{\rho})=\left\{\begin{aligned}
&\dfrac{\Vert m_{\rho}\Vert^{2}}{M_{\rho}(\rho)}\quad&&M_{\rho}(\rho)>0,\\
&0 \quad\quad&&(M_{\rho}(\rho),m_{\rho})=(0,0),\\
&+\infty \quad&&\text{otherwise}.
\end{aligned}
\right. 
\end{equation}

In order to extend the JKO scheme to the system of PFS-MCL model, we need to account for the dynamic boundary condition on $\Gamma$ in \eqref{eq:dynamic_bc}.
For the purpose of illustration, we denote the boundary values of phase-field as a separate variable $\phi_{bc}$ (and this will be made clear in the full-discrete scheme in the next subsection). Then we can easily give its variational formulation based on the minimizing movement scheme for the $\mathcal{L}^{2}$ gradient flow of the energy $\mathcal{E}_{bc}$ associated with $\phi_{bc}$  
\begin{eqnarray}\label{eq:JKO_bc}
	\begin{aligned}
	\phi^{k+1}_{bc}\in\mathop{\text{arg min}}\limits_{\phi_{bc}}\dfrac{\mathrm{Pe}_{s}}{2}d^{2}_{\mathcal{L}^{2}}(\phi_{bc},\phi^{k}_{bc})+\Delta t\mathcal{E}_{bc}(\phi_{bc}),
	\end{aligned}
\end{eqnarray}
whose optimal condition is exactly the implicit Euler scheme of \eqref{eq:dynamic_bc}
\begin{align}
\mathrm{Pe}_{s}(\phi^{k+1}_{bc}-\phi^{k}_{bc})+\Delta t L_\phi=0.
\end{align}

Now we incorporate the coupled dynamics of $(\phi,\psi,\phi_{bc})$ in the JKO scheme, which is a combination of the $\mathcal{H}^{-1}$, $\mathcal{W}^{2}_{m}$ and $\mathcal{L}^{2}$ gradient flows (see Fig.~\ref{fig:FlowChart}). By combining the JKO schemes in \eqref{eq:JKO_vec} and \eqref{eq:JKO_bc} with the dynamics characterization of Wasserstein distance in \eqref{eq:dw} and \eqref{eq:constraint}, we propose the following semi-discrete generalized dynamic JKO scheme \cite{Carrillo2022PrimalDual, Carrillo2024StructurePD}.

\textbf{Problem 1 (Semi-discrete JKO scheme for PFS-MCL model).} Defining $u=(\phi,m_{\phi};\psi,m_{\psi};\phi_{bc})$, given $u^{k}$, solve for $u^{k+1}$ by the following scheme
\begin{eqnarray}\label{eq:semiJKO}
	\begin{aligned}
		&u^{k+1}\in\mathop{\text{arg min}}\limits_{u}\dfrac{1}{2}\int_{\Omega}\Big(\mathcal{D}(\phi,m_{\phi})+\mathcal{D}(\psi,m_{\psi})\Big)\mathrm{d}x+\dfrac{\mathrm{Pe}_{s}}{2}\int_{\Gamma}\Vert\phi_{bc}-\phi^{k}_{bc}\Vert^{2}\mathrm{d}x 
		+\Delta t\mathcal{E}(\phi,\psi,\phi_{bc}),\\
		&\text{s.t.}\left\{\begin{aligned}
		&\phi+\nabla\cdot m_{\phi}=\phi^{k},\quad \psi+\nabla\cdot m_{\psi}=\psi^{k}, \quad \text{in $\Omega$}, \\
		&m_{\phi}\cdot\boldsymbol{\nu}=0,\quad m_{\psi}\cdot\boldsymbol{\nu}=0, \quad \text{on $\partial \Omega$},
		\end{aligned}
		\right.
	\end{aligned}
\end{eqnarray}
where we apply the right endpoint rule for the integration of $\mathcal{D}$ in the inner time $s$ and implicit Euler method for the continuity equation with a one-step strategy ($\Delta s=1$), which does not violate the first-order accuracy of the JKO scheme in outer time \cite{LiWuchen2020Fisher}.

\begin{myremark} (Relation with other schemes)
    Our JKO scheme \eqref{eq:semiJKO} leverages the underlying variational structure of the coupled PFS-MCL model, which enables us to guarantee the original energy dissipation law \eqref{eq:dE_dt} at the discrete level, unlike many other energy-stable schemes only possessing modified dissipation structure (see Section \ref{subsec:comparison_schemes} for numerical comparison). Compared to those fully implicit schemes that treat the Wasserstein gradient flow \eqref{eq:govern_matrix} as $\mathcal{H}^{-1}$ gradient flow weighted by the mobility at the new time step, the JKO scheme inherits even more properties such as bound-preserving (aside from the energy dissipation) from Wasserstein metric (not from the logarithmic potential), which endows itself with extended applicability for general bound-preserving gradient flows with degenerate mobilities (see Remark \ref{rmk:primal_prox} for more discussion).
\end{myremark}

\subsection{Full-discrete variational scheme for PFS-MCL model}\label{sec:full_discrete}
We further discuss the spatial discretization for the semi-discrete variational formulation \eqref{eq:semiJKO} and provide a full-discrete variational formulation. 
Let us consider the two-dimensional case with a rectangular domain $\Omega=\left[a,b\right]\times\left[c,d\right]$, where the substrate boundary is $\Gamma=\left[a,b\right]\times\left\{y=c\right\}$ and the non-substrate boundary is $\partial\Omega\setminus\Gamma$. We divide $\Omega$ uniformly into $N=N_{x}\times N_{y}$ subdomains, denoted as $G_{i,j}=\left[x_{i-1/2},x_{i+1/2}\right]\times\left[y_{j-1/2},y_{j+1/2}\right]$ with the center located at $x_{i}=a+(i-1/2)\Delta x$ and $y_{j}=c+(j-1/2)\Delta y$ and grid spacing $\Delta x=\frac{b-a}{N_{x}}$ and $\Delta y=\frac{d-c}{N_{y}}$, for $i=1,\cdots,N_{x}$ and $j=1,\cdots,N_{y}$. 

We apply the central difference to $\nabla\cdot m_{\rho}$ in the continuity equation
\begin{align}
	\rho_{i,j}+\dfrac{1}{2\Delta x}\big((m^x_\rho)_{i+1,j}-(m^x_\rho)_{i-1,j}\big)+\dfrac{1}{2\Delta y}\big((m^y_\rho)_{i,j+1}-(m^y_\rho)_{i,j-1}\big)=\rho^k_{i,j}, \quad \text{for $\rho=\phi, \psi$}, 
\end{align}
where the involved values at ghost points $(m^{x}_{\rho})_{0,j}$, $(m^x_\rho)_{N_{x}+1,j}$, $(m^y_\rho)_{i,0}$ and $(m^y_\rho)_{i,N_{y}+1}$ can be obtained by applying central difference to the no-flux boundary conditions ($m_{\rho}\cdot\boldsymbol{\nu}=0$). For instance, $(m^x_\rho)_{0,j}=-(m^x_\rho)_{1,j}$. 

For the spatial discretization of the objective function, we perform the trapezoidal rule for the Dirichlet energy associated with $|\nabla\phi|^{2}$ while apply the midpoint rule with a piecewise constant approximation for the other energy terms and the distance functionals to achieve second-order accuracy \cite{Carrillo2024StructurePD}. Denoting $\mathcal{E}^{h}$ the discrete energy, the chemical potentials can be calculated by taking its derivatives with respect to $(\phi_{i,j},\psi_{i,j},(\phi_{bc})_i)$, where $(\phi_{bc})_i=\phi_{i,\frac{1}{2}}$ denotes the boundary values along the substrate
\begin{eqnarray}
	\begin{aligned}
		&\nabla_{\phi}\mathcal{E}^{h}=\Big(-\mathrm{Cn}^{2}L_{i,j}+\dfrac{\partial H}{\partial\phi}\Big|_{i,j}+\dfrac{1}{2\mathrm{Ex}}\dfrac{\partial P_{1}}{\partial\phi}\Big|_{i,j}-\dfrac{1}{4}\dfrac{\partial P_{2}}{\partial\phi}\Big|_{i,j}\Big)\Delta x\Delta y+S_{i,j}\Delta x,\\
		&\nabla_{\psi}\mathcal{E}^{h}=\Big(\mathrm{Pi}\dfrac{\partial G}{\partial\psi}\Big|_{i,j}+\dfrac{1}{2\mathrm{Ex}}\dfrac{\partial P_{1}}{\partial\psi}\Big|_{i,j}+\dfrac{1}{4}\dfrac{\partial P_{2}}{\partial\psi}\Big|_{i,j}\Big)\Delta x\Delta y,\\
        &\nabla_{\phi_{bc}}\mathcal{E}^{h}=\Big(-2\mathrm{Cn}^{2}\Big(\dfrac{\phi_{i,1}-(\phi_{bc})_i}{\Delta y}\Big)+\mathrm{Cn}\gamma^{'}_{\omega f}((\phi_{bc})_i)\Big)\Delta x.
	\end{aligned}
\end{eqnarray}
Here $L_{i,j}=L^x_{i,j}+L^y_{i,j}$ represents the numerical Laplacian of $\phi$
\begin{eqnarray}
	L^{x}_{i,j}=\left\{
	\begin{aligned}
		&\dfrac{\phi_{2,j}-\phi_{1,j}}{(\Delta x)^{2}}&&i=1,\\
		&\dfrac{\phi_{i+1,j}-2\phi_{i,j}+\phi_{i-1,j}}{(\Delta x)^{2}}&&\text{otherwise},\\
		&\dfrac{\phi_{i-1,j}-\phi_{i,j}}{(\Delta x)^{2}}&&i=N_{x}.
	\end{aligned}
	\right.\qquad
	L^{y}_{i,j}=\left\{
	\begin{aligned}
		&\dfrac{\phi_{i,2}-\phi_{i,1}}{(\Delta y)^{2}}&&j=1,\\
		&\dfrac{\phi_{i,j+1}-2\phi_{i,j}+\phi_{i,j-1}}{(\Delta y)^{2}}&&\text{otherwise},\\
		&\dfrac{\phi_{i,j-1}-\phi_{i,j}}{(\Delta y)^{2}}&&j=N_{y},
	\end{aligned}
	\right.
\end{eqnarray}
where we have applied the boundary condition \eqref{eq:natural_bc} on the non-substrate boundaries $\partial \Omega \setminus\Gamma$, and $S_{i,j}$ is associated with the gradient of the Dirichlet energy along the substrate
\begin{eqnarray}\label{Wij}
	\begin{aligned}
		S_{i,j}&=\left\{
		\begin{aligned}
			&2\mathrm{Cn}^{2}\Big(\dfrac{\phi_{i,1}-(\phi_{bc})_i}{\Delta y}\Big)\quad&&j=1,\\ 
			&0&&\text{otherwise}.
		\end{aligned}
		\right.
	\end{aligned}
\end{eqnarray}

Therefore, the full-discrete JKO scheme of the PFS-MCL model is as follows (Fig.~\ref{fig:FlowChart})

\textbf{Problem 2 (Full-discrete JKO scheme for PFS-MCL model).}\label{prob2}
Given $\left\{\phi^{k}_{i,j},\psi^{k}_{i,j},(\phi^k_{bc})_i\right\}$, solve $\left\{\phi^{k+1}_{i,j},\psi^{k+1}_{i,j},(\phi^{k+1}_{bc})_i\right\}$ by following scheme
\begin{eqnarray}\label{eq:fullJKO}
	\begin{aligned}
		&\{u^{k+1}_{i,j}\}\in\mathop{\text{arg min}}\limits_{u}\dfrac{1}{2}\sum_{\rho=\phi,\psi}\sum_{i,j=1}^{N_{x},N_y}\mathcal{D}(\rho_{i,j},{\textbf{m}_{\rho}}_{i,j})\Delta x\Delta y +\dfrac{\mathrm{Pe}_{s}}{2}\sum_{i=1}^{N_{x}}\big((\phi_{bc})_i-(\phi^k_{bc})_i\big)^{2}\Delta x+\Delta t\mathcal{E}^{h}(\boldsymbol{\phi},\boldsymbol{\psi},\boldsymbol{\phi_{bc}})\\
		&\text{s.t.}\left\{
            \begin{aligned}
			&\phi_{i,j}+\dfrac{1}{2\Delta x}\big((m^{x}_{\phi})_{i+1,j}-(m^{x}_{\phi})_{i-1,j}\big)+\dfrac{1}{2\Delta y}\big((m^{y}_{\phi})_{i,j+1}-(m^{y}_{\phi})_{i,j-1}\big)=\phi^{k}_{i,j},\\
			&\psi_{i,j}+\dfrac{1}{2\Delta x}\big((m^{x}_{\psi})_{i+1,j}-(m^{x}_{\psi})_{i-1,j}\big)+\dfrac{1}{2\Delta y}\big((m^{y}_{\psi})_{i,j+1}-(m^{y}_{\psi})_{i,j-1}\big)=\psi^{k}_{i,j},\\
			&(m^{x}_{\phi})_{0,j}=-(m^{x}_{\phi})_{1,j},(m^{x}_{\phi})_{N_{x}+1,j}=-(m^{x}_{\phi})_{N_{x},j},(m^{y}_{\phi})_{i,0}=-(m^{y}_{\phi})_{i,1},(m^{y}_{\phi})_{i,N_{y}+1}=-(m^{y}_{\phi})_{i,N_{y}},\\
			&(m^{x}_{\psi})_{0,j}=-(m^{x}_{\psi})_{1,j},(m^{x}_{\psi})_{N_{x}+1,j}=-(m^{x}_{\psi})_{N_{x},j},(m^{y}_{\psi})_{i,0}=-(m^{y}_{\psi})_{i,1},(m^{y}_{\psi})_{i,N_{y}+1}=-(m^{y}_{\psi})_{i,N_{y}}.\\
		\end{aligned}
		\right.
	\end{aligned}
\end{eqnarray}

The proposed full-discrete JKO scheme for PSF-MCL model preserves the desirable traits of the solution at the discrete level. The proof of the theorem is given in Appendix~\ref{sec:proof}.
\begin{mytheorem}\label{EMB} 
The full-discrete variational scheme has the following structure-preserving properties:

(i) Original energy dissipation structure;

(ii) Mass conservation of $\phi$ and $\psi$;

(iii) Bound-preserving of $\psi$, i.e., $0\leq \psi \leq 1$.
\end{mytheorem}

\begin{myremark} 
(Accuracy order) Since the JKO scheme is analogous to the implicit Euler method in the variational form \cite{Carrillo2022PrimalDual} and we have consistently applied second-order spatial discretization in the integration of the objective function and the PDE constraints, our full-discrete scheme should posesses first-order accuracy in time and second-order accuracy in space. Although we currently do not have rigorous error analysis of the full-discrete scheme (which is not the main focus of this paper), we verify the accuracy order with numerical tests in Section \ref{sec:5.1}.
\end{myremark}

\begin{myremark}\label{rmk:equil_bc}
	(Equilibrium boundary condition) The dynamic boundary condition reduces to the equilibrium boundary condition with infinite contact line mobility (i.e., $\mathrm{Pe}_{s}\rightarrow 0$):
	\begin{align}\label{eq:eq_bc}
	\mathrm{Cn}\nabla\phi\cdot\boldsymbol{\nu}=-\gamma^{'}_{\omega f}(\phi),\qquad \text{on $\Gamma$},
	\end{align}
	which is widely used for other moving contact line models (see \cite{Huang2013Wetting}). The numerical treatment of JKO scheme for the equilibrium boundary conditions is discussed in \cite{Carrillo2024StructurePD} and we will briefly demonstrate in Appendix~\ref{appd:eqbc}.
\end{myremark}

\section{Primal-Dual method}\label{sec:4}
The full-discrete JKO scheme \eqref{eq:fullJKO} can be rewritten as an optimization problem of three functions by introducing a penalty term for the constraints as follows
\begin{eqnarray}\label{eq:minproblem}
	\begin{aligned}
		\mathop{\mathrm{min}}_{u}\Phi(u)+E(u)+i_{\delta}(Au), \quad \text{where }
		i_{\delta}(y)=\left\{
		\begin{aligned}
			&0\quad\Vert y-b\Vert_{2}\leq\delta, \\
			&\infty\quad\text{otherwise},\\
		\end{aligned}
		\right.
	\end{aligned}
\end{eqnarray}
where we define the following variable and functions
\begin{eqnarray}\label{eq:disc_obj}
	\begin{aligned}
		\left\{\begin{aligned}
			&u=\left\{\phi_{i,j},(m^{x}_{\phi})_{i,j}, (m^{y}_{\phi})_{i,j},\psi_{i,j},(m^{x}_{\psi})_{i,j},(m^{y}_{\psi})_{i,j},(\phi_{bc})_i\right\}^{1\leq j\leq N_{y}}_{1\leq i\leq N_{x}},\\
			&\Phi(u)=\dfrac{1}{2}\sum_{\rho=\phi,\psi}\sum_{i,j=1}^{N_{x},N_y}\mathcal{D}(\rho_{i,j},{\textbf{m}_{\rho}}_{i,j})\Delta x\Delta y
            +\dfrac{\mathrm{Pe}_{s}}{2}\sum_{i=1}^{N_{x}}\big((\phi_{bc})_i-(\phi^k_{bc})_i\big)^{2}\Delta x,\\
			&E(u)=\Delta t\mathcal{E}^{h}(\boldsymbol{\phi},\boldsymbol{\psi},\boldsymbol{\phi_{bc}}).
		\end{aligned}
		\right.\\
	\end{aligned}
\end{eqnarray}
Here the constraints of the continuity equation and boundary conditions are written in the matrix-vector form of $Au=b$ (see Appendix~\ref{appd:disA} for details) and further relaxed by $\Vert Ax-b\Vert_{2}\leq \delta$, where the choice of the small relaxation parameter $\delta$ is justified upon the discretization of the continuity equation \cite{Carrillo2022PrimalDual} and the mass conservation can be preserved up to a reasonably small error in implementation (see Figs.~\ref{fig:DiffTwophipsiEne} (c) and \ref{fig:Threedroplet} (c)). In this section, we will mainly discuss the algorithms for solving the above three-operator optimization problem.

\subsection{Primal-Dual method and its acceleration}

A primal-dual method based on modern proximal splitting schemes \cite{Yan2018NewPDThreeFunctions} can be applied to solve the minimization problem \eqref{eq:minproblem} by rewriting it in the convex-concave formulation
\begin{equation}
		\mathop{\mathrm{min}}_{u}\mathop{\mathrm{max}}_{v}\Phi(u)+E(u)+\langle Au, v\rangle - i^{*}_{\delta}(v),
\end{equation}
where $i^*_{\delta}(v)$ is the convex conjugate of $i_{\delta}(Au)$. The proposed algorithm refers to the iteration 
\begin{eqnarray}\label{eq:PD}
	\begin{aligned}
		\begin{cases}
			v^{(k+1)}=\mathrm{Prox}_{\sigma i^{*}_{\delta}}(v^{(k)}+\sigma A\bar{u}^{(k)}),\\
			u^{(k+1)}=\mathrm{Prox}_{\lambda\Phi}(u^{(k)}-\lambda\nabla E(u^{(k)})-\lambda A^{\mathrm{T}}v^{(k+1)}),\\
			\bar{u}^{(k+1)}=2u^{(k+1)}-u^{(k)}+\lambda\nabla E(u^{(k)})-\lambda\nabla E(u^{(k+1)}),
		\end{cases}
	\end{aligned}
\end{eqnarray}
which converges to the saddle point provided that $\sigma\lambda<1/\lambda_{\mathrm{max}}(AA^{\mathrm{T}})$ and $\lambda<2/L_{\nabla E}$ (where $\lambda_{\max}(AA^T)$ is the maximum eigenvalue of $AA^T$ and $L_{\nabla E}$ represents the Lipschitz constant of $\nabla E$). Here $\mathrm{Prox}_{\lambda h}(u)=\mathop{\text{arg min}}\limits_{\tilde{u}}\dfrac{1}{2\lambda}\Vert \tilde{u}-u\Vert^{2}_{2}+h(\tilde{u})$ is the proximal operator. However, since the matrix $AA^{\mathrm{T}}$ is related to the discrete Laplacian whose norm increases with decreasing grid size, the iteration could converge very slowly with required small $\lambda$ and $\sigma$ for a two-dimensional problem with fine grids \cite{Carrillo2024StructurePD} (see Section~\ref{sec:5.1}). 

Inspired by the work \cite{Liu2021AccelerationPrimalDual}, the last author and his collaborators proposed a preconditioned primal-dual method for the dynamic JKO scheme \cite{Carrillo2024StructurePD} to relax the constraints on $\lambda$ and $\sigma$ and accelerate the convergence
\begin{eqnarray}\label{eq:PrePD}
	\begin{aligned}
		\begin{cases}
			v^{(k+1)}=\mathrm{Prox}^{C_{2}}_{i^{*}_{\delta}}(v^{(k)}+C^{-1}_{2}A\bar{u}^{(k)}),\\
			u^{(k+1)}=\mathrm{Prox}^{C_{1}}_{\Phi}(u^{(k)}-C^{-1}_{1}\nabla E(u^{(k)})-C^{-1}_{1} A^{\mathrm{T}}v^{(k+1)}),\\
			\bar{u}^{(k+1)}=2u^{(k+1)}-u^{(k)}+C^{-1}_{1}\nabla E(u^{(k)})-C^{-1}_{1}\nabla E(u^{(k+1)}).
		\end{cases}
	\end{aligned}
\end{eqnarray}
where the convergence condition is $C_2 \succeq A C_1^{-1} A^{\mathrm{T}}$, and the extended proximal
operators is defined as a minimization 
\begin{align}\label{eq:preprox}
\mathrm{Prox}^{C}_{h}(u)=\mathop{\text{arg min}}\limits_{\tilde{u}}\dfrac{1}{2}\Vert \tilde{u}-u\Vert^{2}_{C}+h(\tilde{u})\quad \text{with $\Vert z\Vert^{2}_{C}:=z^{\mathrm{T}}Cz$}.
\end{align}
Here we choose $C_{1}=\frac{1}{\lambda}I$ and $C_{2}=\lambda AA^{\mathrm{T}}$ (with $\lambda>0$ being a tuning parameter) in \eqref{eq:PrePD} from both perspectives of faster convergence and the ease of computing the extended proximal operators (see Section~\ref{subsec:prox_comp}). Since the preconditioned primal-dual method is quite robust, the parameter $\lambda$ is chosen the same for most of the numerical simulations in Section~\ref{sec:5}. We refer the reader to \cite{Carrillo2024StructurePD} for detailed discussion on the choice of $\lambda$.

\begin{algorithm}[htbp]
	\caption{Preconditioned primal-dual method}\label{Al1}
	\KwIn{$u^{(0)},v^{(0)},T_{\mathrm{end}},Iter_{\mathrm{\mathrm{max}}},\lambda,\Delta t>0$.}
	\KwOut{$u^{*},v^{*}$.}
	Let $\bar{u}^{(0)}=u^{(0)}$, $k=0$ and $N=\dfrac{T_{\mathrm{end}}}{\Delta t}$.\\
	\For{$j\leq N$}{
		\For{$k<Iter_{\mathrm{max}}$}{
			\textbf{repeat}\\
			$v^{(k+1)}=\mathrm{Prox}^{C_{2}}_{i^{*}_{\delta}}(v^{(k)}+C^{-1}_{2}A\bar{u}^{(k)})$\\
			$u^{(k+1)}=\mathrm{Prox}^{C_{1}}_{\Phi}(u^{(k)}-C^{-1}_{1}\nabla E(u^{(k)})-C^{-1}_{1}A^{\mathrm{T}}v^{(k+1)})$\\
			$\bar{u}^{(k+1)}=2u^{(k+1)}-u^{(k)}+C^{-1}_{1}\nabla E(u^{(k)})-C^{-1}_{1} \nabla E(u^{(k+1)})$\\
			where, $C_{1}=\frac{1}{\lambda}I$ and $C_{2}=\lambda AA^{\mathrm{T}}$.\\
			\textbf{until}\quad stopping criteria is achieved\\
			$u^{*}=u^{(k+1)}$\\
			$v^{*}=v^{(k+1)}$
		}
	}
\end{algorithm}

We will employ the preconditioned primal-dual method in Algorithm~\ref{Al1} to solve the optimization problem for one step of the full-discrete JKO scheme with the following initial guesses of $u$ and $v$
\begin{equation*}
		u^{0}=(\phi^{0},\mathbf{0}_{N},\mathbf{0}_{N},\psi^{0},\mathbf{0}_{N},\mathbf{0}_{N},\phi^0_{bc})^{\mathrm{T}}, \quad v^{0}=(\mathbf{0}_{N},\mathbf{0}_{N})^{\mathrm{T}}.
\end{equation*}
The stopping criteria for iterations consist of the satisfaction monitor of PDE constraints
\begin{eqnarray}\label{eq:stop_pde}
	\begin{aligned}
		\Vert Au^{(k+1)}-b\Vert_{2}&\leq\delta,
	\end{aligned}
\end{eqnarray}
and the convergence monitors of relative errors in variables $(u,v)$ and functionals $(E,\Phi)$
\begin{eqnarray}
	\begin{aligned}
		&\mathrm{max}\Biggl\{\dfrac{\Vert u^{(k+1)}-u^{(k)}\Vert}{\Vert u^{(k+1)}\Vert},\dfrac{\Vert v^{(k+1)}-v^{(k)}\Vert}{\Vert v^{(k+1)}\Vert}\Biggr\}\leq\epsilon_{1},\\
		&\mathrm{max}\Biggl\{\dfrac{\Vert E(u^{(k+1)})-E(u^{(k)})\Vert}{\Vert E(u^{(k+1)})\Vert},\dfrac{\Vert \Phi(u^{(k+1)})-\Phi(u^{(k)})\Vert}{\Vert \Phi(u^{(k+1)})\Vert}\Biggr\}\leq\epsilon_{2}.
	\end{aligned}
\end{eqnarray}

\subsection{Calculating the proximal operators} \label{subsec:prox_comp}
The success of our approach based on the JKO scheme depends on the ease of computing the extended proximal. We will demonstrate the computation of the extended proximal operators in the following sequence.

\subsubsection{Proximal operator for dual variables: $\mathrm{Prox}^{C_{2}}_{i^{*}_{\delta}}(v)$}

The proximal operator $\mathrm{Prox}^{C_{2}}_{i^{*}_{\delta}}(v)$ can be calculated based on the Moreau identity 
\begin{equation}\label{eq:moreau}
		\mathrm{Prox}^{C_{2}}_{i^{*}_{\delta}}(v)=v-C^{-1}_{2}\mathrm{Prox}_{i_{\delta}}^{C_{2}^{-1}}(C_{2}v),
\end{equation}
where $\mathrm{Prox}_{i_{\delta}}^{C_{2}^{-1}}$ is related to a problem of minimizing a quadratic form over a 2-sphere
\begin{align}
\mathrm{Prox}_{i_{\delta}}^{C_{2}^{-1}}(y)=\mathop{\text{arg min}}\limits_{\tilde{y}}\dfrac{1}{2}\Vert \tilde{y}-y\Vert^{2}_{C_2^{-1}} \quad \text{subject to $\Vert\tilde{y}-b \Vert_{2}\leq \delta$}.
\end{align}
This minimization is an important {\it trust-region subproblem} and has been intensively discussed in literature \cite{Sorensen1997minquad,Hager2001SSM,Hager2004SSM}. The minimizer $y^*=\mathrm{Prox}_{i_{\delta}}^{C_{2}^{-1}}(y)$ is a solution to the following linear system \cite[Lemma 2.1]{Hager2001SSM}
\begin{equation}\label{eq:linear_mu}
    (\mathcal{I}+\mu C_{2})(y^*-b)=y-b 
\end{equation}
where $\mathcal{I}$ is the identity matrix and the parameter $\mu$ satisfies the following conditions:
\begin{itemize}
    \item If $\Vert y-b \Vert_{2}< \delta$, $\mu=0$ and hence $y^{*}=y$;
    \item If $\Vert y-b \Vert_{2}\geq \delta$, $\mu$ is a suitable parameter such that $C_{2}^{-1}+\mu\mathcal{I}\succeq 0$ and $\Vert y^{*}-b\Vert_{2}=\delta$. 
\end{itemize}
The mainstream methods for this problem involve solving a series of eigenvalue problems and finding the parameter $\mu$ by successive iterations. The bottleneck of the algorithms is the costly matrix decomposition for solving the eigenvalues problems. Fortunately, the spectrum of the matrix $\mathcal{I}+\mu C_{2}$ can be easily obtained in our problem since the matrix $C_{2}=\lambda AA^{\mathrm{T}}$ corresponds to the rescaled discrete Laplacian operator with homogeneous Neumann boundary conditions with staggered grids (see details in Appendix~\ref{appd:inversion}). 

Suppose $\lambda_{1}\geq\lambda_{2}\geq\cdots\geq\lambda_{n}$ are the eigenvalues of $C_{2}$ and $\phi_{1}, \phi_{2}, \cdots, \phi_{n}$ are the corresponding orthonormal eigenvectors. Following \cite[Lemma 2.2]{Hager2001SSM}, we can obtain the solution as an expansion of eigenvectors.
\begin{align}
y^{*}=b+\sum c^{*}_{i}\phi_{i},
\end{align}
with the expansion coefficients $\{c^*_i\}$ chosen in the following way:
\begin{itemize}
    \item[(a)] Degenerate case: If $c_{i}=\langle y-b,\phi_{i}\rangle=0$ for $i\in\{i: \lambda_{i}=\lambda_{1}\}$ and $\sum_{\lambda_{i}<\lambda_{1}} \frac{c_{i}^{2}}{(1+\mu\lambda_{i})^{2}}\leq\delta^{2}$, then $\mu=-1/\lambda_{1}$, and $c^{*}_{i}=c_{i}/(1-\lambda_{i}/\lambda_{1})$ for ${i}\in\{{i}:\lambda_{i}<\lambda_{1}\}$; the $c^{*}_{i}$ for ${i}\in\{i: \lambda_{i}=\lambda_{1}\}$ are arbitrary scalars satisfying $\sum_{\lambda_{i}=\lambda_{1}}{c^{*}_{i}}^{2}=\delta^{2}-\sum_{\lambda_{i}<\lambda_{1}} {c^{*}_{i}}^{2}$.
    \item[(b)] Non-degenerate case: If (a) does not hold, then $c^{*}_{i}=c_{i}/(1+\mu\lambda_{i})$ for all $i$, where $\mu>-1/\lambda_{1}$ is determined by the relation $\sum_{i}{c^{*}_{i}}^{2}=\delta^{2}$.
\end{itemize}
Given the spectrum $\{\lambda_{i}\}$, the above algorithm is very efficient, where the parameter $\mu$ in the non-generate case can be easily calculated by Newton iterations. Furthermore, once $\mu$ is obtained, we can apply the fast Fourier transform (FFT) to solve the equation \eqref{eq:linear_mu}, which accelerates the computation of the expansion coefficients $\{c_{i}\}$ by the forward transform and the summation $y^{*}=b+\sum c^{*}_{i}\phi_{i}$ by inverse transform (see Appendix~\ref{appd:inversion}). 

\begin{myremark} \label{rmk:dual_prox1} (Reducing computational complexity of $\mathrm{Prox}^{C_{2}}_{i^{*}_{\delta}}$)
According to the update of dual variable in the iteration \eqref{eq:PrePD} and the formula for the proximal operator \eqref{eq:moreau}, we would need to compute the matrix inversion $C_{2}^{-1}$ possibly twice and the matrix-vector multiplication $C_{2}v$ multiple times, which is computational costly. To reduce the redundant computations, we introduce variables $y^{(k)}=C_{2}z^{(k)}$ and $y^{*}=\mathrm{Prox}_{i_{\delta}}^{C_{2}^{-1}}(y^{(k)})$, where $z^{(k)}=v^{(k)}+C_{2}^{-1}A\bar{u}^{(k)}$ is the input of $\mathrm{Prox}^{C_{2}}_{i^{*}_{\delta}}$ in \eqref{eq:PrePD}. By \eqref{eq:moreau} and \eqref{eq:linear_mu}, we have 
\begin{eqnarray}
    \begin{aligned}
        v^{(k+1)}&=C_2^{-1}(y^{(k)}-y^*)=\mu(y^*-b), \\
        y^{(k+1)}&=C_2v^{(k+1)}+A\bar{u}^{(k+1)} = y^{(k)}-y^*+A\bar{u}^{(k+1)}.
    \end{aligned}
\end{eqnarray}
By the above formulation, given $y^{(k)}$, we can apply FFT-based fast algorithm to obtain $y^*$ and update $v^{(k+1)}$ and $y^{(k+1)}$ by only one matrix-vector multiplication for $A\bar{u}^{(k+1)}$, without any matrix inversion or matrix-vector multiplication for $C_2$. 
\end{myremark}

\begin{myremark} \label{rmk:dual_prox2} (Inexact update of $\mathrm{Prox}^{C_{2}}_{i^{*}_{\delta}}$). 
Suggested by the work \cite{Liu2021AccelerationPrimalDual}, the dual-subproblem of computing $\mathrm{Prox}^{C_{2}}_{i^{*}_{\delta}}$ in Algorithm \ref{Al1} can be solved inexactly with certain accuracy without violating the convergence. Therefore, we can also approximates $\mathrm{Prox}^{C_{2}}_{i^{*}_{\delta}}$ by the projection onto the ball with radius $\delta$ centered at $b$
\begin{eqnarray}\label{eq:proxC2}
	\begin{aligned}
		\mathrm{Prox}^{C_{2}}_{i^{*}_{\delta}}(y)&\approx y-C^{-1}_{2}\left\{
		\begin{aligned}
			&C_{2}y\quad&&\text{$\Vert C_{2} y-b\Vert_{2}<\delta$},\\
			&\delta\dfrac{C_{2}y-b}{\Vert C_{2}y-b\Vert}_{2}+b\quad&&\text{otherwise}.\\
		\end{aligned}
		\right.
	\end{aligned}
\end{eqnarray}
When $\delta$ is very small, the above projection is good enough to guarantee the overall convergence of the preconditioned primal-dual method, which is has been employed and verified by various numerical examples in Section~\ref{sec:5}. 

We can still avoid redundant matrix inversions and matrix-vector multiplications by introducing axillary variables $z^{(k)}=\bar{v}^{(k)}+r^{(k)}$ where $\bar{v}^{(k)}=C_{2}v^{(k)}$ and $r^{(k)}=A\bar{u}^{(k)}-b$, which yields 
\begin{eqnarray}\label{eq:proxideltaybar}
    \begin{aligned}
        v^{(k+1)}=C_{2}^{-1}\bar{v}^{(k+1)}, \quad \text{where }
        \bar{v}^{(k+1)}&=\left\{
        \begin{aligned}
            &\textbf{0} \quad&&\text{$\Vert z^{(k)}\Vert_{2}<\delta$},\\
            &\Big(1-\dfrac{\delta}{\Vert z^{(k)}\Vert}_{2}\Big)z^{(k)} \quad&&\text{otherwise}.\\
        \end{aligned}
        \right.
    \end{aligned}
\end{eqnarray}
By the above formulation, given $\bar{v}^{(k)}$, we only need to compute one matrix-vector multiplication for $r^{(k)}$ to update $z^{(k)}$ and one matrix inversion to get $v^{(k+1)}$. Meantime, the matrix inversion can be efficiently computed in log-linear time by using FFT-based matrix inversion (see details in Appendix~\ref{appd:inversion}).
\end{myremark}

\subsubsection{Proximal operator for primal variables: $\mathrm{Prox}^{C_{1}}_{\Phi}(u)$}

With $C_1=\frac{1}{\lambda}I$, the extended proximal operator reduces to the standard case $\mathrm{Prox}^{C_{1}}_{\Phi}(u)=\mathrm{Prox}_{\lambda\Phi}(u)$. Since $\Phi(u)$ is separable with respect to the variables $\{(\phi_{i,j},{\textbf{m}_{\phi}}_{i,j});(\psi,{\textbf{m}_{\psi}}_{i,j});(\phi_{bc})_i\}$, its proximal operator is component-wise,
\begin{equation}
\mathrm{Prox}_{\lambda\Phi}(u)=
\left\{\mathrm{Prox}_{\lambda\mathcal{D}}(\phi_{i,j},{\textbf{m}_{\phi}}_{i,j});
\mathrm{Prox}_{\lambda\mathcal{D}}(\psi,{\textbf{m}_{\psi}}_{i,j});
\mathrm{Prox}_{\lambda\Vert\cdot\Vert}((\phi_{bc})_i)\right\}_{1\leq i \leq N_x}^{1\leq j \leq N_y},
\end{equation}
where the proximal operator for each component pair is defined as follows
\begin{eqnarray}
    \begin{aligned}
    &\mathrm{Prox}_{\lambda\mathcal{D}}(\phi,\textbf{m}_{\phi})=\mathop{\text{arg min}}\limits_{(\tilde{\phi},\tilde{\textbf{m}}_{\phi})}\left\{\dfrac{1}{2}(\tilde{\phi}-\phi)^{2}+\dfrac{1}{2}\Vert \tilde{\textbf{m}}_{\phi}-\textbf{m}_{\phi}\Vert^{2}+\lambda\dfrac{\Vert\tilde{\textbf{m}}_{\phi}\Vert^{2}}{2M_{\phi}}\right\}:=\mathop{\text{arg min}}\limits_{(\tilde{\phi},\tilde{\textbf{m}}_{\phi})}F_{1},\\
    &\mathrm{Prox}_{\lambda\mathcal{D}}(\psi,\textbf{m}_{\psi})=\mathop{\text{arg min}}\limits_{(\tilde{\psi},\tilde{\textbf{m}}_{\psi})}\left\{\dfrac{1}{2}(\tilde{\psi}-\psi)^{2}+\dfrac{1}{2}\Vert\tilde{\textbf{m}}_{\psi}-\textbf{m}_{\psi}\Vert^{2}+\lambda\dfrac{\Vert\tilde{\textbf{m}}_{\psi}\Vert^{2}}{2M_{\psi}(\tilde{\psi})}\right\}:=\mathop{\text{arg min}}\limits_{(\tilde{\psi},\tilde{\textbf{m}}_{\psi})}F_{2},\\
    &\mathrm{Prox}_{\lambda\Vert\cdot\Vert}(\phi_{bc})=\mathop{\text{arg min}}\limits_{\tilde{\phi}_{bc}}\left\{\dfrac{1}{2}(\tilde{\phi}_{bc}-\phi_{bc})^{2}+\lambda\dfrac{\mathrm{Pe}_{s}}{2}(\tilde{\phi}_{bc}-\phi_{bc})^{2}\right\}:=\mathop{\text{arg min}}\limits_{\tilde{\phi}_{bc}}F_{3}.
    \end{aligned}
\end{eqnarray}
	
The proximal operators for $(\tilde{\phi},\tilde{\textbf{m}}_{\phi})$ and $\tilde{\phi}_{bc}$ are trivial since $F_{1}$ and $F_{3}$ are quadratic. By taking derivatives, the optimality conditions gives the solutions of proximal operators
\begin{eqnarray}
	\begin{aligned}
		\left\{\begin{aligned}
			&\dfrac{\partial F_{1}}{\partial\tilde{\phi}}=\tilde{\phi}-\phi=0,
			\quad \dfrac{\partial F_{1}}{\partial\tilde{\textbf{m}}_{\phi}}=\tilde{\textbf{m}}_{\phi}-\textbf{m}_{\phi}+\lambda\dfrac{\tilde{\textbf{m}}_{\phi}}{M_{\phi}}=0,\\
			&\dfrac{\partial F_{3}}{\partial \tilde{\phi}_{bc}}=(1+\lambda\mathrm{Pe}_{s})(\tilde{\phi}_{bc}-\phi_{bc})=0.
		\end{aligned}
		\right.\\
	\end{aligned}
    \Rightarrow
	\begin{aligned}
		\left\{\begin{aligned}
			&\mathrm{Prox}_{\lambda\mathcal{D}}(\phi,\textbf{m}_{\phi})=(\phi,\dfrac{M_{\phi}\textbf{m}_{\phi}}{M_{\phi}+\lambda}),\\
			&\mathrm{Prox}_{\lambda\Vert\cdot\Vert}(\phi_{bc})=\phi_{bc}.
		\end{aligned}
		\right.\\
	\end{aligned}
\end{eqnarray}

The proximal operator $\mathrm{Prox}_{\lambda\mathcal{D}}(\psi,\textbf{m}_{\psi})$ can be computed by the following formula \cite{Carrillo2024StructurePD}
\begin{eqnarray}\label{eq:prox_psi}
	\begin{aligned}
		\mathrm{Prox}_{\lambda\mathcal{D}}(\psi,\textbf{m}_{\psi})=\left\{
		\begin{aligned}
			&(\psi^{*},\dfrac{M_{\psi}(\psi^{*})\textbf{m}_{\psi}}{M_{\psi}(\psi^{*})+\lambda})& \quad
			&\text{if $-\frac{\Vert \textbf{m}_{\psi}\Vert^{2}}{2\lambda}< \psi <1+\frac{\Vert \textbf{m}_{\psi}\Vert^{2}}{2\lambda}$}, \\
			&(0,\textbf{0})& \quad &\text{if $\psi\leq-\frac{\Vert \textbf{m}_{\psi}\Vert^{2}}{2\lambda}$},\\
                &(1,\textbf{0})& \quad &\text{if $\psi\geq1+\frac{\Vert \textbf{m}_{\psi}\Vert^{2}}{2\lambda}$}.\\
		\end{aligned}
		\right.
	\end{aligned}
\end{eqnarray}
where $\psi^{*}$ is only root of the following polynomial such that $\psi^{*}\in(0,1)$
\begin{eqnarray}
	\begin{aligned}
		f(\tilde{\psi})=(\tilde{\psi}-\psi)(\lambda+M_{\psi}(\tilde{\psi}))^{2}-\dfrac{\lambda }{2}M^{'}_{\psi}(\tilde{\psi})\Vert\textbf{m}_{\psi}\Vert^{2}=0,
	\end{aligned}
\end{eqnarray}
which can be efficiently computed by the Newton's method with a tailored plan for choosing initial guesses. The initial-value plan is designed based on the monotonicity and concavity of $f(\tilde{\psi})$ and is proved to guarantee the convergence of the Newton iterations to the desired root $\psi^{*}\in(0,1)$ \cite{Carrillo2024StructurePD}.

\begin{myremark}\label{rmk:primal_prox} (Bound-preserving and computational efficiency of $\mathrm{Prox}^{C_{1}}_{\Phi}$)
One main advantage of our approach is to guarantee the strict bound-preserving of $\psi\in [0,1]$ at each iteration of Algorithm \ref{Al1} by the computation of $\mathrm{Prox}_{\lambda\mathcal{D}}(\psi,\textbf{m}_{\psi})$, even with the inexact dual proximal solver in Remark \autoref{rmk:dual_prox2}. We have rigorously analyzed the proximal operator in our previous work \cite{Carrillo2024StructurePD}, where the solution of $\mathrm{Prox}_{\lambda\mathcal{D}}$ can be distinguished into cases whether a root $\psi^{*}\in(0,1)$ of $f(\tilde{\psi})$ exists or not, as given in \eqref{eq:prox_psi}.
If $\psi^{*}\in(0,1)$ exists, the solution of $\mathrm{Prox}_{\lambda\mathcal{D}}$ is attained at such $\psi^{*}$ and can be computed by a provable convergent Newton method; otherwise, $\mathrm{Prox}_{\lambda\mathcal{D}}$ is necessarily attained at the endpoints of $(0,1)$ and, by the definition of $\mathcal{D}$, we have $\textbf{m}_{\psi}^{*}=\textbf{0}$. Moreover, by leveraging the Wasserstein metric of the gradient flows, this bound-preserving approach does not depend on the energy potential and is applicable for the degenerate mobility paired with any energy functionals, for instance the double-well potential.

The Newton's method for finding the root does not decrease the overall computational efficiency given that it involves only element-wise computations for scalar functions and converges very fast with tailored initial values in implementation (different from the case of using Newton's method for large-scale coupled nonlinear systems). Furthermore, since the whole proximal operator $\mathrm{Prox}_{\lambda\Phi}(u)$ is also component-wise, it can be parallelized to enhance the computational efficiency. 
\end{myremark}

\subsection{Adaptive time stepping}
A main difficulty in the simulations of phase-field model is to accurately capture the instantaneous topological change of phase field and meantime achieve the correct long-time equilibrium state. Adaptive time stepping strategy is usually employed to increase the computational efficiency \cite{Qiao2011AdaptiveTime,Huang2023AdaptiveTime}, which however requires the stability property of the numerical scheme to allow for large time steps and a good monitor to detect when the solution change drastically or slightly. Our variational scheme possesses both unconditional stability (provided the well-posedness of the JKO scheme \cite{Lisini2012CahnHilliard}) and the original energy dissipation property, which makes it suitable for applying adaptive time stepping. 

Since our variational scheme respects the original energy dissipation structure \eqref{eq:dE_dt}, we use the relative rate of change of energy as the monitor for adaptive time stepping as follows: 
\begin{eqnarray}
	\begin{aligned}
		\Delta t =\mathrm{max}\Bigg(\Delta t_{\mathrm{min}},\dfrac{\Delta t_{\mathrm{max}}}{\sqrt{1+\beta\Big|\mathcal{R}|E^{'}(t_{n-1})|\Big|^{2}}}\Bigg),\quad \mathcal{R}|E^{'}(t_{n-1})|=\dfrac{E(t_{n})-E(t_{n-1})}{E(t_{n-1})(t_{n}-t_{n-1})}.
	\end{aligned}
\end{eqnarray}

\begin{itemize}
	\item $\Delta t_{\mathrm{min}}$ corresponds to the lower bound of step size for the fastest energy decay.
	\item $\Delta t_{\mathrm{max}}$ corresponds to the upper bound of step size for the slowest energy decay.
	\item $\beta$ is an adjustable parameter, typically taken within $\left[10^{1},10^{8} \right]$.
	\item $\mathcal{R}|E^{'}(t_{n-1})|$ is the relative rate of change of energy.
\end{itemize}
According to the adaptive strategy, the time step is set to be small when the energy rapidly changes (i.e., $\mathcal{R}|E^{'}(t)|$ is large) to capture the drastic phase change (e.g., topological change); while the time step is set to be large when the energy slowly changes (i.e., $\mathcal{R}|E^{'}(t)|$ is small) to speed up the simulation.

\begin{myremark}\label{rmk:adap_time} (Adaptive time stepping strategies)
There are several adaptive time stepping strategies for gradient flows with different monitors associated with the numerical energy change rate (i.e., LHS of \eqref{eq:dE_dt}) or the numerical energy dissipation rate (i.e., RHS of \eqref{eq:dE_dt}), which are not necessarily equivalent at discrete level. Since our scheme preserves the original energy dissipation structure, we can also propose another formulation of adaptive time stepping with respect to energy dissipation rate, which is equivalent to the present one. 
\end{myremark}

\section{Numerical experiments}\label{sec:5}

In this section, we will demonstrate the performance of our primal-dual splitting methods for the PFS-MCL model and investigate the effects of surfactants on the contact line dynamics through a series of numerical experiments. In Sec.~\ref{sec:5.1}, we first perform a number of benchmark experiments to demonstrate the performance of our methods based on the JKO scheme. In Sec.~\ref{sec:5.2}, we investigate the effects of the boundary mobility, surfactant concentration and its temperature-dependent diffusion rate on the moving conact line dynamics. Finally, in Sec.~\ref{sec:5.3} we will show some interesting examples illustrating the influence of surfactants on the wetting dynamics of droplets and the dewetting dynamics of liquid thin film on substrate.

Unless specifically stated, in the subsequent numerical experiments we take the following values for the important parameters of the model and our algorithm
\begin{eqnarray}
	\begin{aligned}
	   &\mathrm{Pe_{\phi}}=20, &\quad & \mathrm{Pe_{\psi}}=100,&\quad& \mathrm{Pe_{s}}=1/500,&\quad& \mathrm{Pi}=0.1481,&\quad&
        \mathrm{Ex}=1,\\
        &\lambda=100, &\quad&
        \delta=10^{-7},&\quad& \epsilon_1=10^{-5},&\quad&
        \epsilon_2=10^{-5}.
	\end{aligned}
\end{eqnarray}

\subsection{Benchmark experiments}\label{sec:5.1}
\subsubsection{Accuracy order}
\begin{figure}[htbp]
	\centering
	\subfigure[$\mathrm{Cn}=0.025$,\quad$\theta_{s}=120^{\circ}$]{
		\begin{minipage}{0.31\textwidth}
			\includegraphics[scale=0.35]{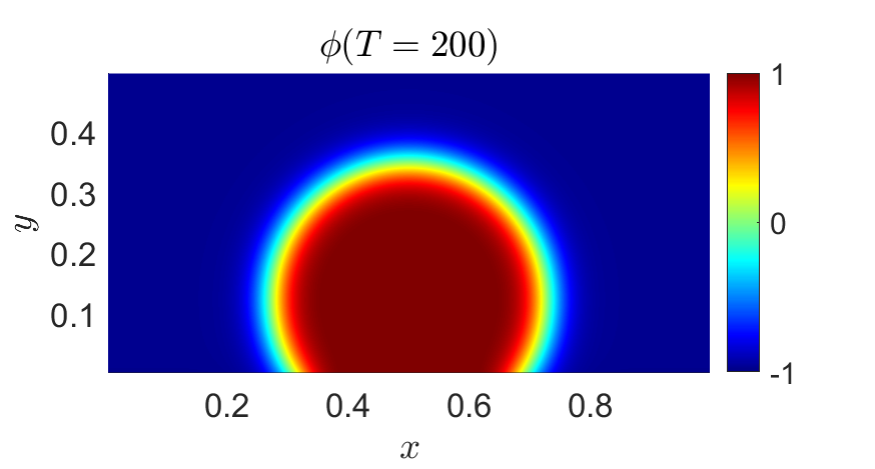}
			\includegraphics[scale=0.35]{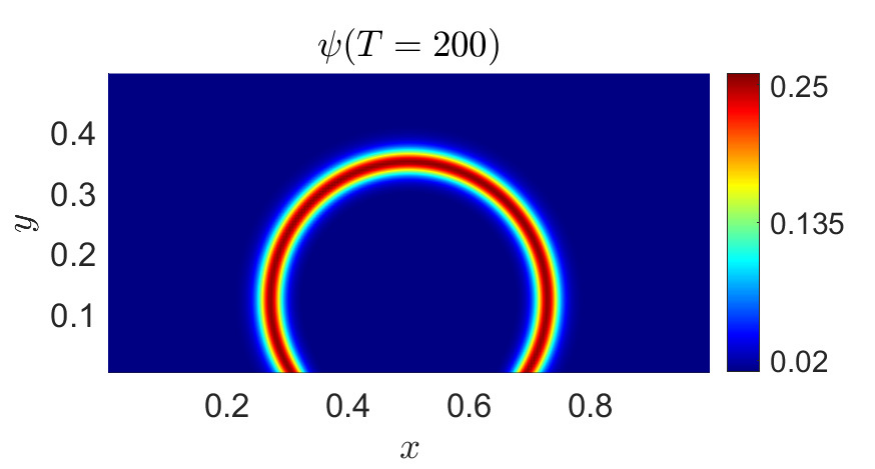}
		\end{minipage}
	}
	\subfigure[Energy dissipation]{
		\begin{minipage}{0.31\textwidth}			
            \includegraphics[scale=0.35]{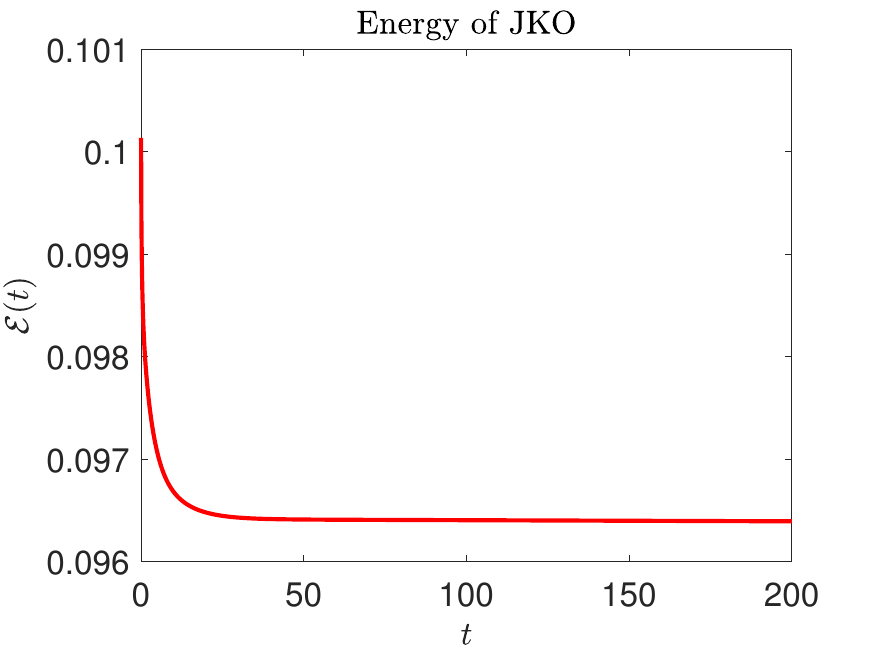}
		\end{minipage}
	}
    \subfigure[Mass conservation]{
		\begin{minipage}{0.31\textwidth}
			\includegraphics[scale=0.35]{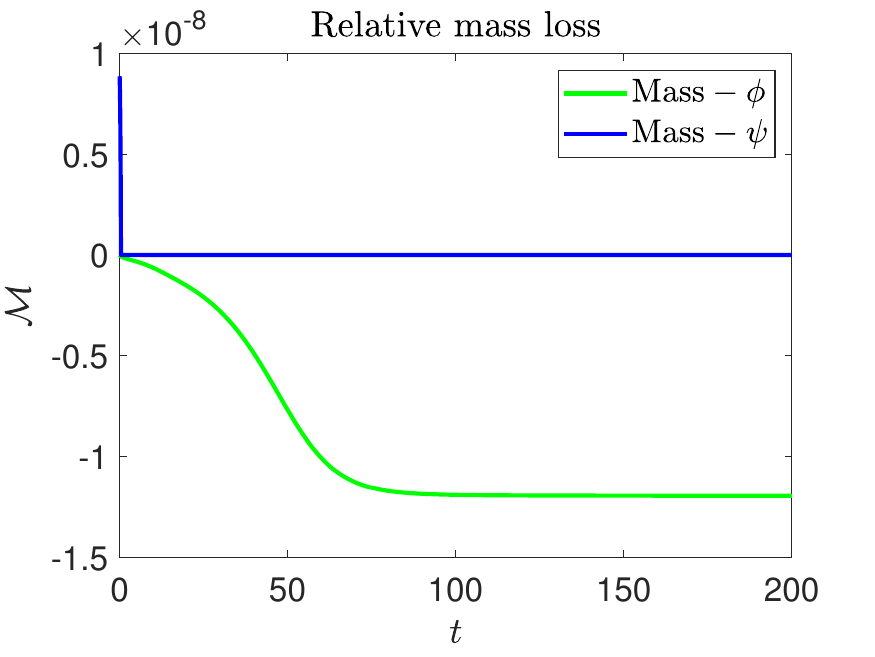}
		\end{minipage}
	} 			
	\caption{The equilibrium states of droplets and surfactants at $T = 200$ for $\mathrm{Cn}=0.025$ and $\theta_s=120^\circ$, with illustration of energy dissipation and mass conservation.}
	\label{fig:DiffTwophipsiEne}
\end{figure}

We first test the first-order accuracy in time for our JKO scheme. Consider the computational domain $\left[0,1\right]\times\left[0,0.5\right]$ and set the mesh size $\Delta x=\Delta y=0.005$. Choose the following initial values for $\phi$ and $\psi$
\begin{eqnarray}\label{eq:IC1}
	\begin{aligned}
		&\phi_{0}(x,y)=\mathrm{tanh}\Big(\dfrac{0.3-\sqrt{(x-0.5)^{2}+y^{2}}}{\sqrt{2}\mathrm{Cn}}\Big),\\
		&\psi_{0}(x,y)=0.02+0.001\xi,
	\end{aligned}
\end{eqnarray}
where $\xi$ is a random variable uniformly distributed in $\left[0,1\right]$. We take the interface thickness parameter $\mathrm{Cn}=0.025$ and the static contact angle $\theta_{s}=120^{\circ}$. 
The equilibrium states of $\phi$ and $\psi$ (at $T=200$) are plotted in Fig.~\ref{fig:DiffTwophipsiEne} (a), where we observe that surfactants are adsorbed and distributed along the diffuse interface. We also check the structure-preserving properties of our numerical schemes for energy dissipation in Fig.~\ref{fig:DiffTwophipsiEne} (b) and mass conservation, which is preserved up to the order of $\delta=10^{-8}$, in Fig.~\ref{fig:DiffTwophipsiEne} (c). 

\begin{figure}[htbp]
 	\centering
	\begin{minipage}[c][0.25\textheight][c]{0.53\textwidth}
		\centering
		\resizebox{1.0\columnwidth}{!}{
			\begin{tabular}{cccccc}
				\toprule[0.5mm]
				& $\Delta t$ & $\Vert\phi-\phi_{\mathrm{ref}}\Vert_{\infty}$ & Order & $\Vert\psi-\psi_{\mathrm{ref}}\Vert_{\infty}$ & Order \\
				\midrule
				\midrule
				\multirow{5}{*}{\rotatebox{90}{\text{First order}}}
				&$1/50$ & 1.32E-2  & - & 4.89E-3 & -  \\
				&$1/100$ & 6.61E-3  & 0.996 & 2.49E-3 & 0.974  \\
				&$1/200$ & 3.27E-3  & 1.015 & 1.23E-3 & 1.019  \\
				&$1/400$ & 1.62E-3  & 1.019 & 5.62E-4 & 1.129 \\
				&$1/800$ & 7.88E-4  & 1.035 & 2.56E-4 & 1.134 \\
				\midrule
				\bottomrule[0.5mm]
		\end{tabular}}
		\tabcaption{First-order accuracy in time of $\phi$ and $\psi$ at $t=0.1$ with fixed $\Delta x=\Delta y=0.005$ and different time steps $\Delta t$, where the reference solution is obtained with $\Delta t=10^{-5}$.}
		\label{tab:FirstTimeorder}
	\end{minipage}
	\begin{minipage}[c][0.25\textheight][t]{0.45\textwidth}
		\includegraphics[width=0.85\textwidth]{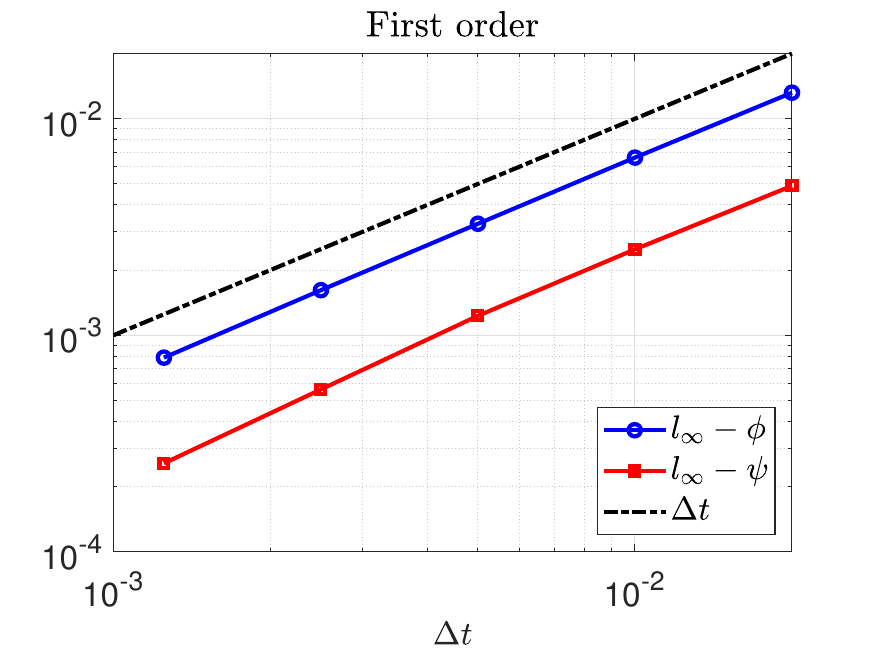}
		\centering	
		\figcaption{First-order accuracy for $\phi$ and $\psi$.}
		\label{fig:FirstOrderCurve}
	\end{minipage}
\end{figure}

We compare the numerical solutions (at $t = 0.1$) obtained using different time steps $\Delta t=1/50$, $1/100$, $1/200$, $1/400$, and $1/800$ with the reference solution obtained using $\Delta t=10^{-5}$. As expected, Table~\ref{tab:FirstTimeorder} and Fig.~\ref{fig:FirstOrderCurve} show the temporal first-order accuracy of the JKO scheme.

\subsubsection{Comparison with other energy-stable schemes}\label{subsec:comparison_schemes}
In this subsection, we demonstrate the importance of preserving the original energy dissipation at the discrete level for accurately capturing the interface dynamics for phase-field modeling. For illustrative purpose, we consider a benchmark test of the evolution of a clean quarter droplet initially sitting near the corner of a computational domain $[0,1]^2$, without wall energy effect (i.e. $\psi(x,t)=0$ and $\gamma_{\omega f}(\phi)=0$)
\begin{align}
    \phi(x,t)=\mathrm{tanh}\Big(\dfrac{\sqrt{x^{2}+y^{2}}-0.2}{\sqrt{2}\mathrm{Cn}}\Big), \quad \mathrm{Cn}=0.02.
\end{align}

\begin{figure}[htbp]
	\centering  
	\includegraphics[width=0.4\textwidth]{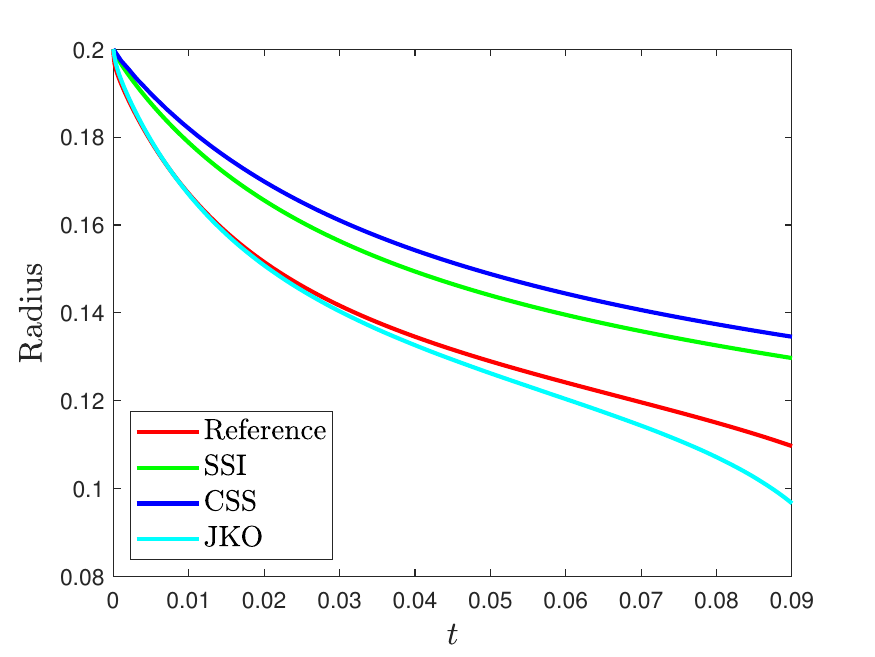}
	\caption{Comparson of SSI, CSS and JKO schemes. Here, $\Delta t = 10^{-4}$ and $\Delta x=\Delta y = 0.005$.}
	\label{fig:DiffCHR}
\end{figure}

We simulate the evolution of the droplet by applying our JKO scheme, and compare with the stabilized semi-implicit scheme (SSI) and convex splitting scheme (CSS) (see Appendix~\ref{appd:schemes}), which are among the most popular energy-stable (modified energy dissipative) schemes for phase-field models. The numerical results of the radius evolution by different schemes with same spatial resolution and time step are compared against the reference solution in Fig.~\ref{fig:DiffCHR}, where the reference solution is computed by a fully implicit scheme with finite element method on a highly refined spatial grid ($\Delta x=\Delta y=2.5\times10^{-3}$) with very small time steps ($\Delta t=10^{-7}$). 
We observe that the JKO scheme approximates the reference solution well while the SSI and CCS schemes exhibit a noticeable deviation from the beginning of evolution, indicating its superior performance in capturing the interface dynamics. This result is consistent with the theoretical analysis and numerical experiments in \cite{Xu2019CMAME} discussing the time-delay behavior of the SSI and CSS schemes for phase-field models when ratio of the time step v.s. the thickness parameter $\Delta t/\mathrm{Cn^{2}}$ is not small enough. Interestingly, the dynamics computed by the JKO scheme exhibits a slight time-advance behavior, similar to the result of the fully implicit scheme investigated in \cite{Xu2019CMAME}. This indicates the close connection of the JKO scheme to the fully implicit scheme, in which the $\mathcal{H}^{-1}$ distance weighted by the mobility at the new time step is used to approximate the Wasserstein distance in \eqref{eq:JKO_vec}. In the constant mobility case, these two distances are equivalent, leading to the same observation in Fig. \ref{fig:DiffCHR} as that in \cite{Xu2019CMAME}. However, in general cases when the variable mobility depends on the phase-field parameter itself, the approximation of the Wasserstein distance in the discrete energy dissipation by the weighted $\mathcal{H}^{-1}$ distance at the new time step may lead to additional errors in the fully implicity scheme.

\subsubsection{Adaptive time stepping and preconditioned primal-dual methods}
Now we test the performance of adaptive time stepping strategy and the effect of the parameter $\beta$. To check the capability of JKO scheme with adaptive time stepping to capture the topological phase change, we consider the wetting dynamics of two neighboring droplets in the computational domain $\left[0,1\right]\times\left[0,0.4\right]$, with the following initial values
\begin{eqnarray}
	\begin{aligned}
		&\phi_{0}(x,y)=1+\sum_{i=1}^{2}\mathrm{tanh}\Big(10-\dfrac{\sqrt{(x-x_{i})^{2}+y^{2}}}{\sqrt{2}\mathrm{Cn}}\Big),\quad\text{where $x_{1}=0.25$ and $x_{2}=0.75$}.\\
        &\psi_{0}(x,y)=0.07+0.001\xi,
	\end{aligned}
\end{eqnarray}
where we take $\mathrm{Cn}=0.01$ and $\theta_{s}=60^{\circ}$. 

\begin{figure}[htbp]
	\centering
	\subfigure[$\mathrm{Cn}=0.01$,\quad$\theta_{s}=60^{\circ}$]{
		\begin{minipage}{0.31\textwidth}
			\includegraphics[scale=0.35]{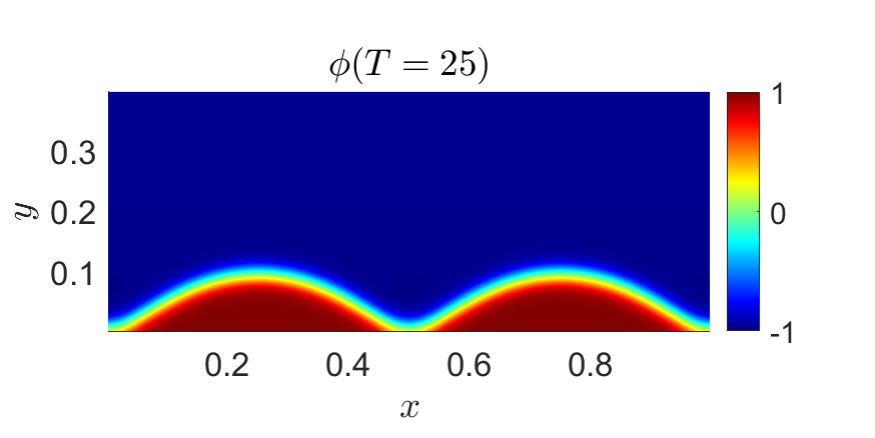}
			\includegraphics[scale=0.35]{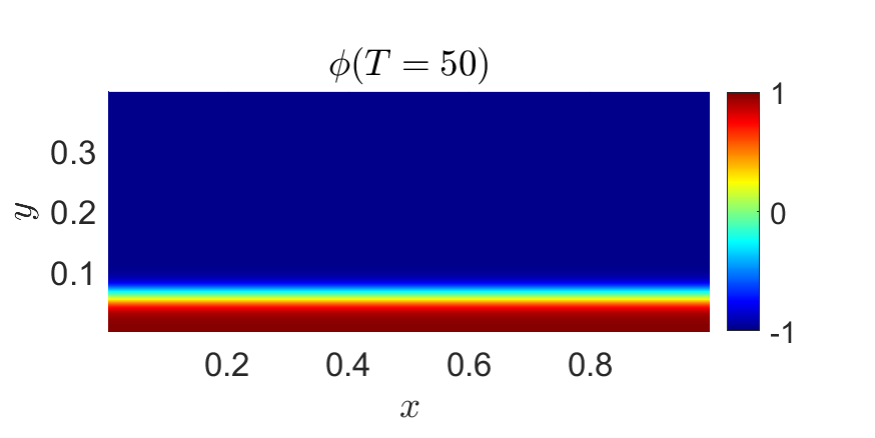}
		\end{minipage}
	}
	\subfigure[$\Delta t_{\mathrm{max}}=0.1$]{
		\begin{minipage}{0.31\textwidth}
			\includegraphics[scale=0.35]{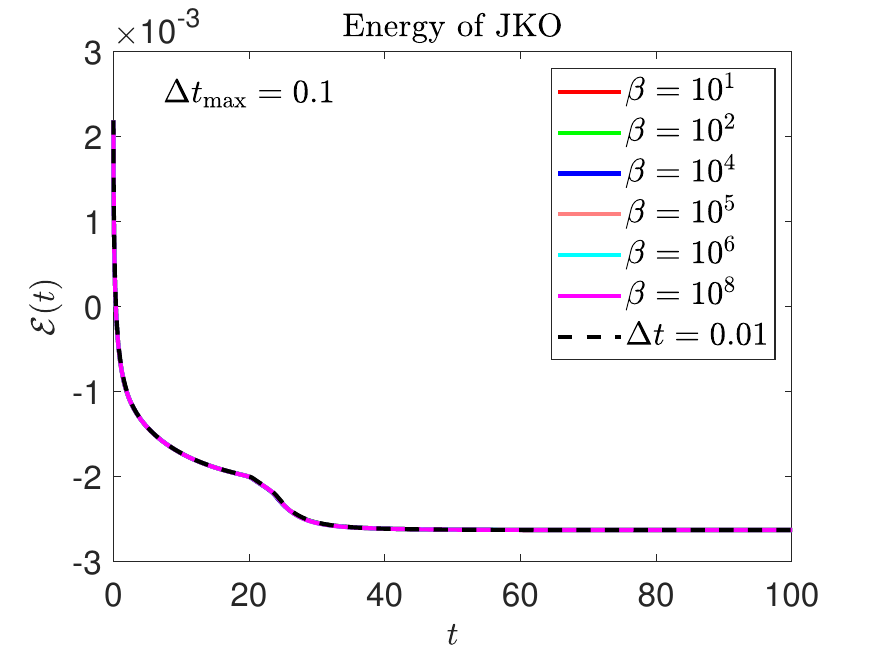}
		\end{minipage}
	}
	\subfigure[$\Delta t_{\mathrm{max}}=0.5$]{
		\begin{minipage}{0.31\textwidth}
			\includegraphics[scale=0.35]{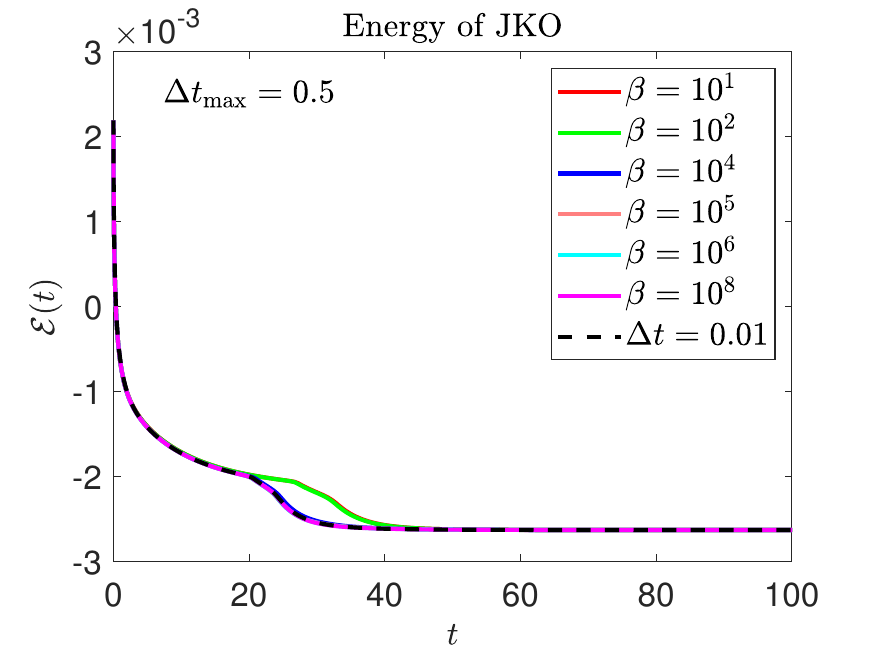}
		\end{minipage}
	}
	\caption{Time snapshots of $\phi$ at $T = 25$ and $T = 100$ illustrating the wetting dynamics of droplets and the energy evolution with different $\beta$ for $\Delta t_{\mathrm{max}}=0.1$ and $\Delta t_{\mathrm{max}}=0.5$. The solution computed with $\Delta t=0.01$ is also plotted as the reference.}
	\label{fig:TwodropletAdabeta}
\end{figure}

We simulate the wetting dynamics of the two droplets for $t\in[0,100]$ in the presence of surfactant using adaptive time stepping with $\Delta t\in[\Delta t_{\mathrm{min}}=0.01, \Delta t_{\mathrm{max}}]$. This dynamics concerning coalescence of inkjet droplets is widely studied and has important applications in industrial printing process \cite{Chai2024pressure}. As shown in Fig.~\ref{fig:TwodropletAdabeta} (a), two droplets that initially sat separately will merge and eventually completely wet the substrate. The corresponding energy evolution for $\Delta t_{\mathrm{max}}=0.1$ and $\Delta t_{\mathrm{max}}=0.5$ with different $\beta$ is shown in Fig.~\ref{fig:TwodropletAdabeta} (b) and (c). We find that the energy evolution contains two rapid decays and two plateaus. It is important to see that the second rapid energy decay around $t=25$ corresponding to the coalescence of two droplets is accurately captured by $\Delta t_{\mathrm{max}}=0.1$ with all $\beta\in[10,10^8]$, while the timing of the coalescence cannot be accurately predicted by $\Delta t_{\mathrm{max}}=0.5$ when $\beta<10^4$. 

\begin{figure}[htbp]
	\centering  
	\includegraphics[width=0.3\textwidth]{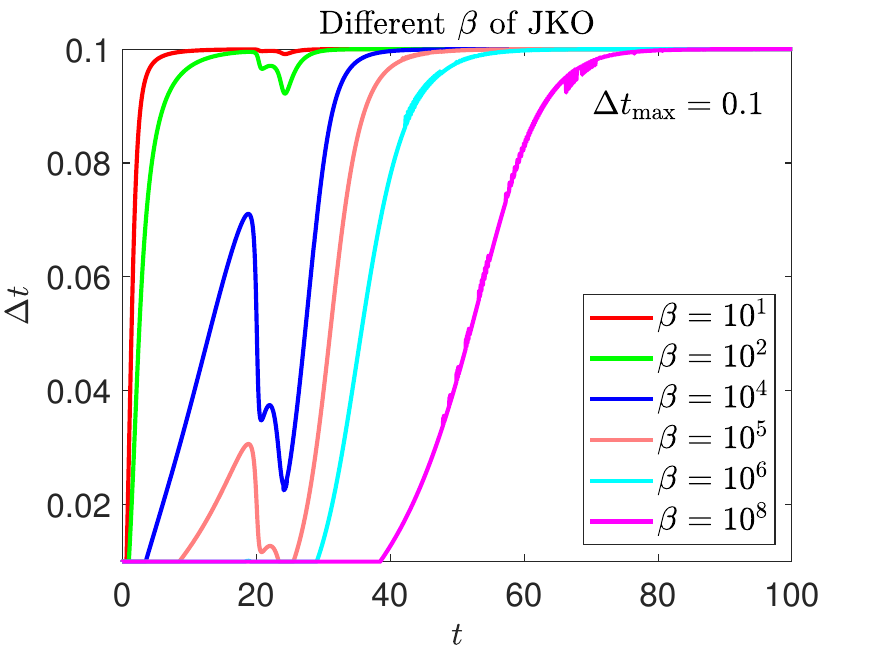}
	\includegraphics[width=0.3\textwidth]{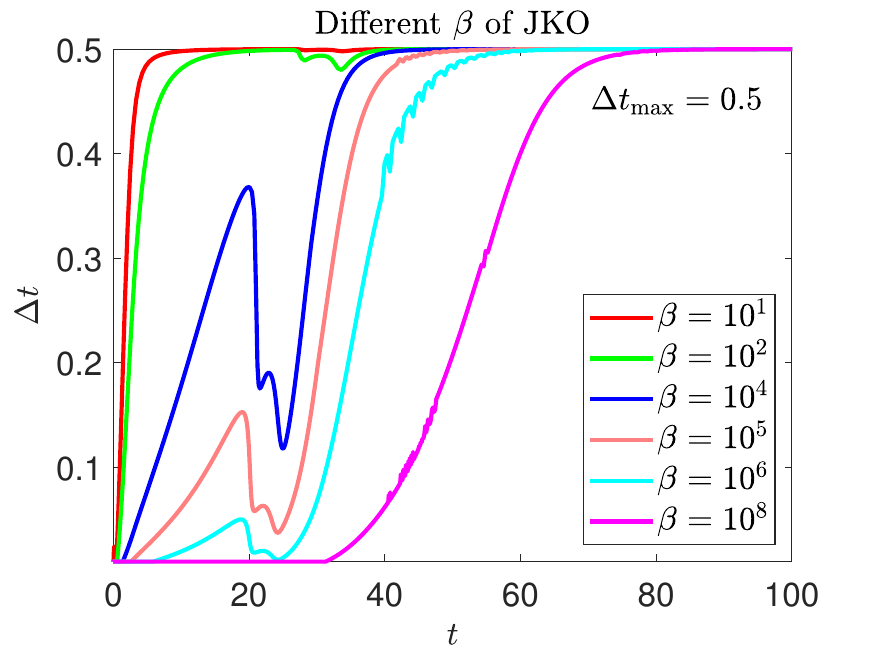}
	\includegraphics[width=0.3\linewidth]{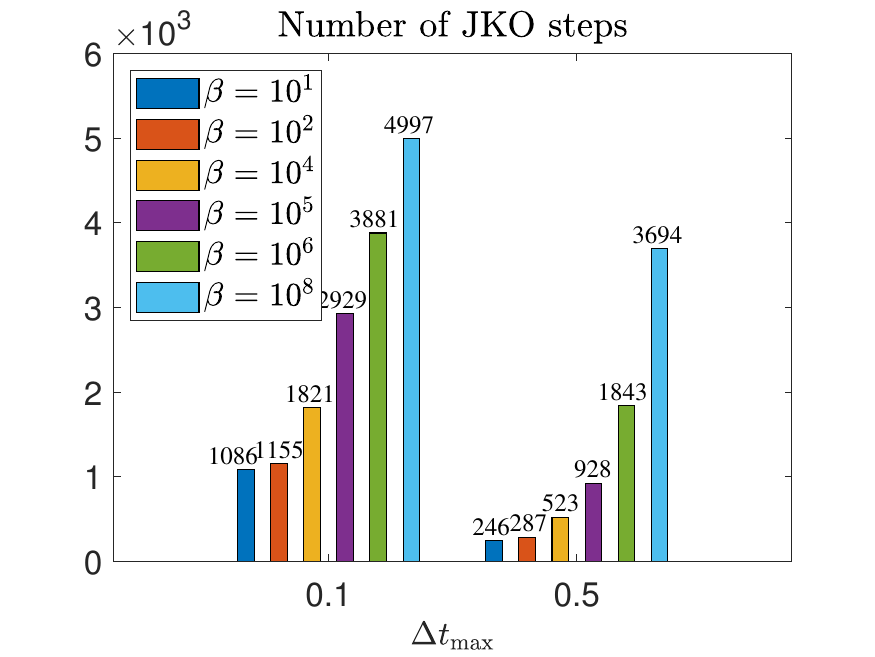}
	\caption{Evolution of adaptive time step and the total number of JKO steps for different $\beta$ and $\Delta t_{\mathrm{max}}$.}
	\label{fig:EneJKOstep}
\end{figure}

We show how the adaptive time step changes and the required total number of JKO steps for simulations in Fig.~\ref{fig:EneJKOstep}. We can clearly observe the adjustment of the time step (e.g., for $\beta=10^4$): $\Delta t$ first increases when the energy enters its first plateau and then adaptively decreases during the second rapid energy decay corresponding to the coalescence of droplets around $t=25$ (see Fig.~\ref{fig:TwodropletAdabeta} (b)) , and eventually increases to $\Delta t_{\mathrm{max}}$ as the dynamics approaches equilibrium. The histogram in Fig.~\ref{fig:EneJKOstep} shows that the total number of JKO steps is significantly reduced by using adaptive time stepping, especially when $\Delta t_{\mathrm{max}}$ is large and $\beta$ is small. However, $\Delta t_{\mathrm{max}}$/$\beta$ should not be too large/small to simulate the correct dynamics. The parameters should be chosen appropriately when using adaptive time stepping for the balance between accuracy and efficiency. In the PFS-MCL model under our consideration, we suggest choose $\beta=10^4$ as the best parameter for adaptive time stepping.

\begin{figure}[htbp]
	\vspace{-0.5cm}
	\centering
	\begin{minipage}[c][0.25\textheight][c]{1\textwidth}
		\centering
		\tabcaption{Comparisons of computational efficiency between PD3O, PrePD and PrePD-Ada, where $\Delta t_{\mathrm{min}}=0.01$, $\Delta t_{\mathrm{max}}=0.1$ or $t_{\mathrm{max}}=0.5$, $\beta=10^{4}$ and $T_{\mathrm{end}}=50$.}
		\label{tal:PDTimeComparison}
		\begin{tabular}{c|c|c|c|c|c|c}\hline
			& JKO Steps & Max Iter & Min Iter & Mean Iter & Total Iter & CPU Time \\ 
			\hline
			PD3O($0.01$) & 5000 & 9974 & 835 & 2260 & 11298422 & 445156 \\
			\hline
			PrePD($0.01$) & 5000 & 3509 & 12 & 268 & 1340300 & 37922 \\
			\hline
            PrePD-Ada($0.1$) & 1318 & 3509 & 19 & 127 & 167539 & 4965\\
            \hline
			PrePD-Ada($0.5$) & 416 & 3509 & 22 & 250 & 104149 & 3105\\
			\hline
		\end{tabular}
	\end{minipage}
	\vspace{-1cm}
\end{figure}

We also investigate the performance of the primal-dual method and its accelerated version by preconditioning used for solving the opitimization problem of JKO schemes. Table~\ref{tal:PDTimeComparison} shows the number of JKO steps, the number of iterations for convergence, and the corresponding CPU time for the primal-dual method for three operators (PD3O) \eqref{eq:PD}, the preconditioned primal-dual method (PrePD) (Algorithm \ref{Al1}) and the preconditioned primal-dual method with adaptive time stepping (PrePD-Ada). 
The table mainly demonstrates three important findings from the comparison between these algorithms. First, PrePD has a much faster convergence rate than PD3O and requires only $\sim 1/8$ iterations of PD3O for convergence. Second, although PrePD may involve matrix inversion for each iteration, unlike PD3O which does not need matrix inversion, the CPU time for PrePD versus PD3O still scales with their total iterations indicating the high efficiency and necessity of FFT-based fast algorithms. Third, the overall computational efficiency (in terms of CPU time) is improved by adaptive time-stepping strategy, despite the fact that it may require more iterations for one JKO step (revealed by Mean Iter) due to the larger time step. 

\subsection{Moving contact line dynamics}\label{sec:5.2}
\subsubsection{Effect of boundary mobility}
\begin{figure}[htbp]
	\centering  
	\includegraphics[width=0.35\textwidth]{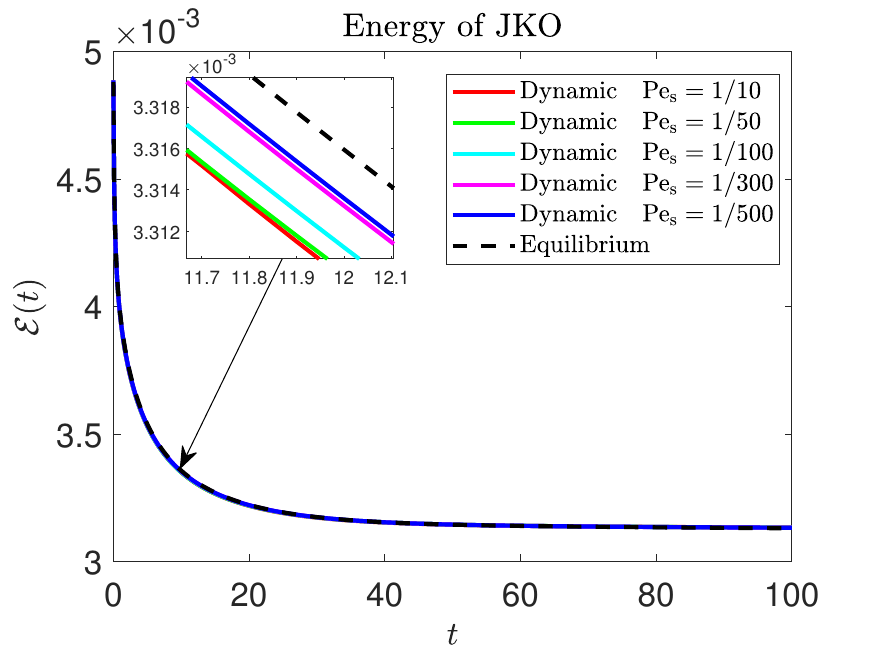}
	\includegraphics[width=0.35\textwidth]{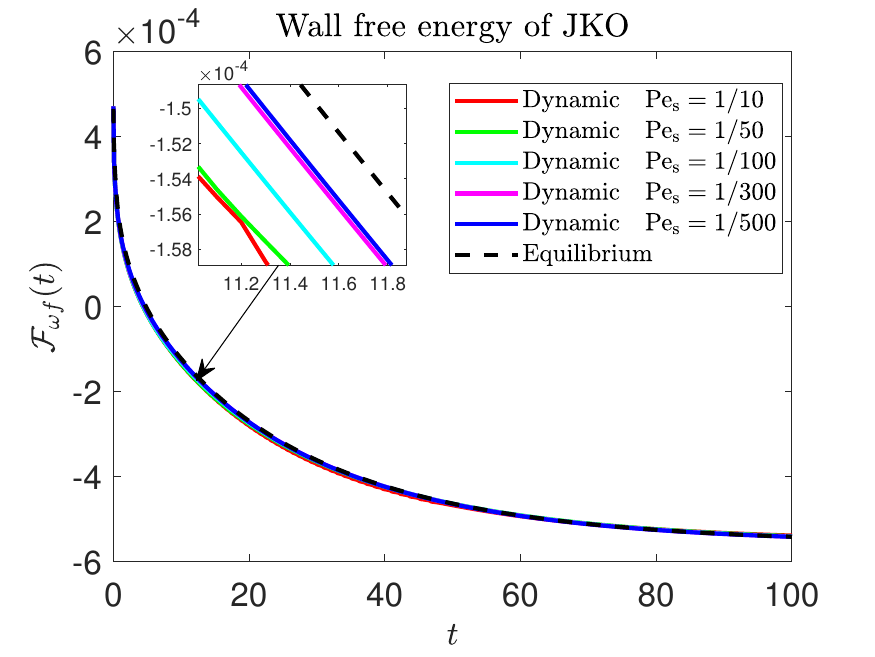}
	\caption{The evolution of total energy ($\mathcal{E}$) and wall free energy ($\mathcal{F}_{\omega f}$) for different $\mathrm{Pe}_{s}$.}
	\label{fig:DiffPesEneWall}
\end{figure}

We first investigate the effect of the mobility $\mathrm{Pe}_{s}$ on the dynamics of contact line using the initial conditions \eqref{eq:IC1} with $\mathrm{Cn}=0.01$ and $\theta_{s}=120^{\circ}$. Fig.~\ref{fig:DiffPesEneWall} shows the evolution of total energy ($\mathcal{E}$) and wall free energy ($\mathcal{F}_{\omega f}$) for dynamic boundary condition \eqref{eq:dynamic_bc} with different $\mathrm{Pe}_{s}$ and equilibrium boundary condition \eqref{eq:eq_bc} (corresponding to $\mathrm{Pe}_{s}\rightarrow 0$). We can only observe subtle differences between the results with different $\mathrm{Pe}_{s}$ indicating that the dynamics and the equilibrium state is not sensitive to the mobility of contact line. 
This is probably a result of the presence of surfactants, which introduces additional energy dissipation into the moving contact line system and has an effect of enhancing contact line dynamics \cite{Zhu2020PhaseFieldMCL,Zhao2021Surfactant}. A comprehensive investigation of the friction parameter $\mathrm{Pe}_{s}$ is beyond the scope of this work. For more details of the effects of $\mathrm{Pe}_{s}$ as well as $\mathrm{Pe}_{\phi}$, we refer the readers to the theoretical work \cite{Xu2018sharp} and the numerical study \cite{Kang2020multiple} for the phase-field model with MCL in the absence of surfactants.

\subsubsection{Effect of surfactants on droplet shapes and contact angles}
\begin{figure}[htbp]
	\centering
	\includegraphics[width=0.3\textwidth]{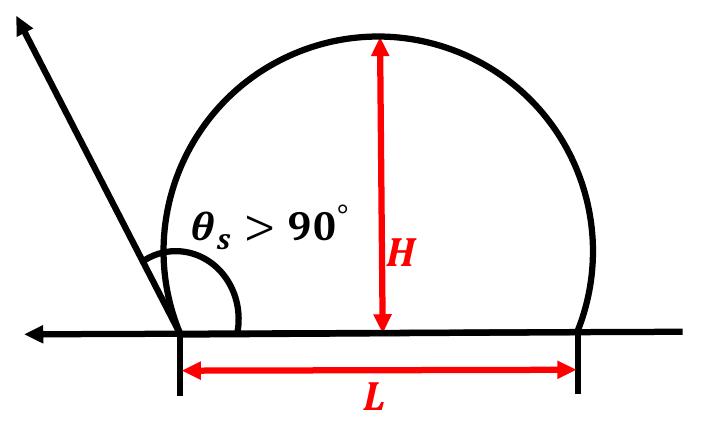}
    \includegraphics[width=0.3\textwidth]{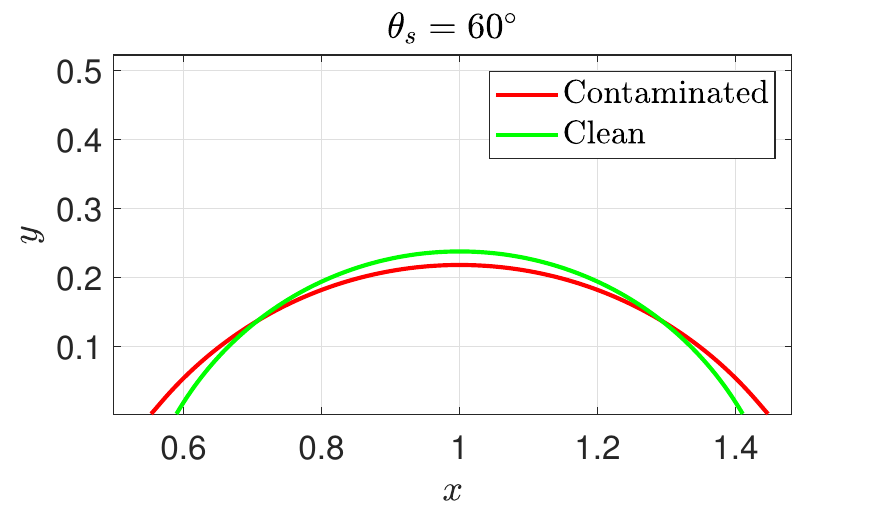}
    \includegraphics[width=0.3\textwidth]{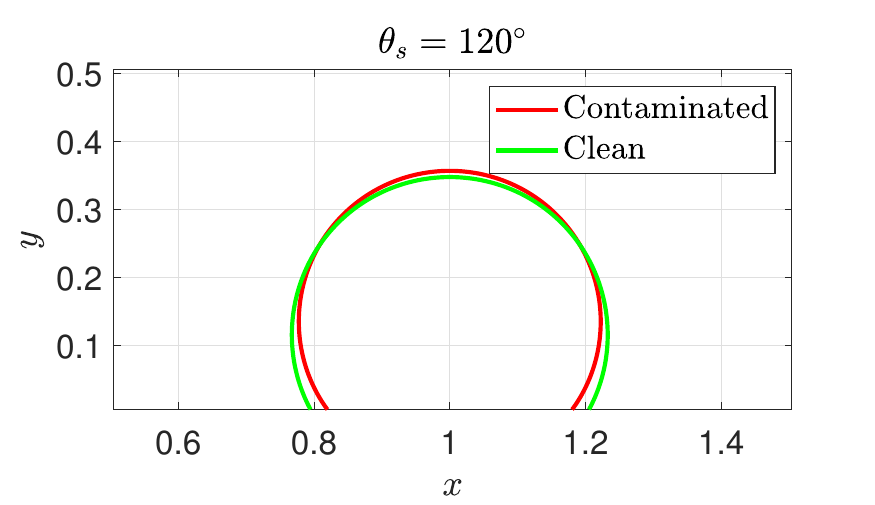}
	\caption{Left: Illustration of the relation between the droplet height $H$, spreading length $L$ and equilibrium contact angle $\theta_{s}$ of a clean droplet. Middle: Equilibrium shapes of clean and contaminated hydrophobic droplets for $\theta_s=60^{\circ}$. Right: Equilibrium shape of clean and contaminated hydrophobic droplets for $\theta_s=120^{\circ}$.}
	\label{fig:RHLangle}
\end{figure}

Next, we investigate the influence of surfactants on the equilibrium shape and the contact angle of a droplet. The dynamics of the moving contact line can be quantitatively translated into the relationship between the spreading length $L$ and droplet height $H$, as illustrated in Fig.~\ref{fig:RHLangle} (left). We start from an initial semicircular droplet with its radius being $R_{0}$ and its contact angle being $\pi/2$, sitting in the middle of the computational domain. For a clean droplet (no surfactant concentration), the relation between $L$, $H$ and the initial radius $R_0$, contact angle $\theta_s$ in equilibrium can be analytically derived using the principle of mass conservation \cite{Cai20143Dwetting} as follows
\begin{align}\label{eq:RHLangle}
    L=2R_{0}\sqrt{\dfrac{\pi}{2(\theta_{s}-\mathrm{sin}\theta_{s}\mathrm{cos}\theta_{s})}}\mathrm{sin}\theta_{s},\quad H=R_{0}\sqrt{\dfrac{\pi}{2(\theta_{s}-\mathrm{sin}\theta_{s}\mathrm{cos}\theta_{s})}}(1-\mathrm{cos}\theta_{s}).
\end{align}

\begin{figure}[htbp]
	\centering  
	\includegraphics[width=0.35\textwidth]{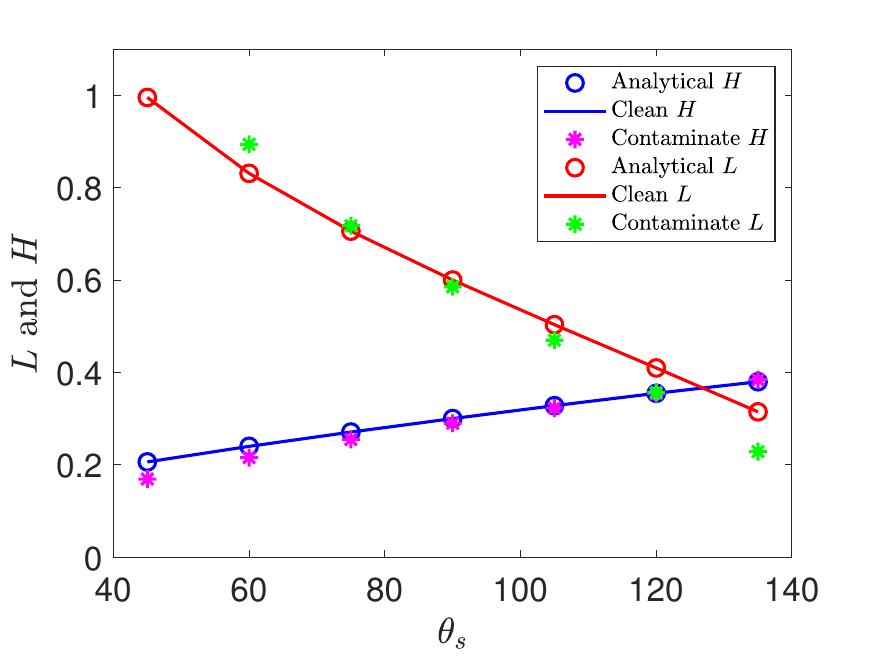}
    \includegraphics[width=0.35\textwidth]{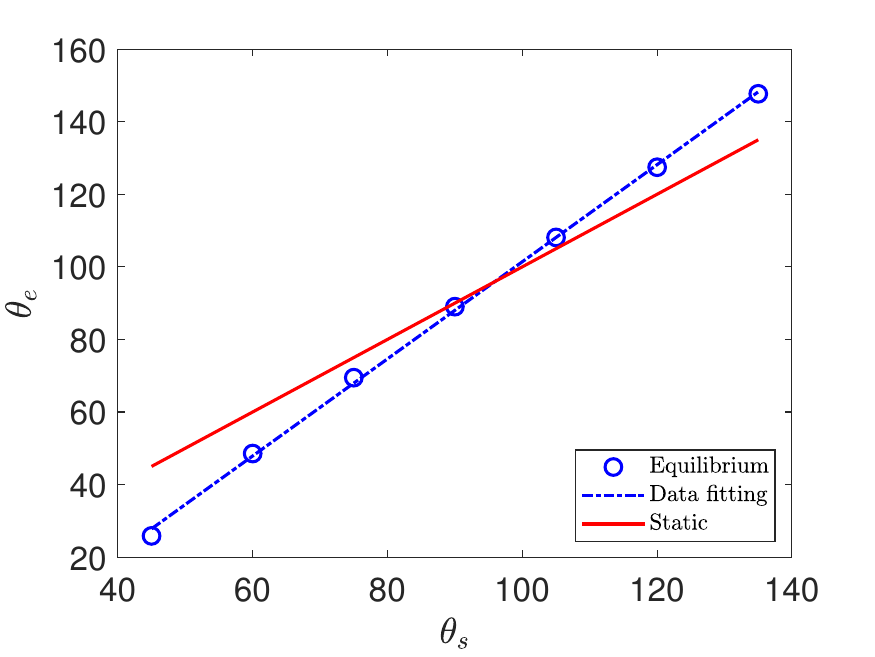}
	\caption{Left: The spreading length $L$ and height $H$ versus the effective equilibrium contact angle $\theta_s$ of clean droplets and contaminated droplets in equilibrium obtained by numerical simulations, compared against the analytical relations in \eqref{eq:RHLangle}. Right: Linear relationship between effective equilibrium contact angle $\theta_e$ and prescribed static contact angle $\theta_s$: $\theta_{e}=1.338\theta_{s}-32.396$.}
	\label{fig:RHL}
\end{figure}

We set $R_{0} = 0.3$ in the simulation, and consider the dynamics of a clean droplet ($\psi=0$) and a contaminated droplet ($\psi=0.02$) for both hydrophilic and hydrophobic conditions ($\theta_{s}$ ranging from $45^{\circ}$ to $135^{\circ}$), with $\mathrm{Cn}=0.01$ and mesh size $\Delta x=\Delta y=0.005$ for a long time $T=200$ with $\Delta t=0.1$ until the droplet reaches its equilibrium state. 
As shown in Fig.~\ref{fig:RHL} (left), for a clean droplet, the numerical results of $L$ and $H$ v.s. $\theta_s$ match well with the analytical relations in \eqref{eq:RHLangle}; while for a contaminated droplet, the presence of surfactants exhibits a significant influence on the equilibrium shapes \cite{Zhu2019ThermodynamicallyMCL,Zhu2020PhaseFieldMCL}. Moreover, from Fig.~\ref{fig:RHL} (right) we observe that the effective equilibrium contact angle $\theta_e$ of the contaminated droplets is deviated from the prescribed static contact angle $\theta_s$ due to the presence of surfactants. Specifically, the hydrophilic and hydrophobic properties of droplets can be enhanced in the presence of surfactants, i.e. $\theta_e<\theta_s<90^\circ$ for hydrophilic droplets (Fig.~\ref{fig:RHLangle} middle) and $\theta_e>\theta_s>90^\circ$ for hydrophobic droplets (Fig.~\ref{fig:RHLangle} right). 

\begin{figure}[htbp]
	\centering 
	\includegraphics[width=0.35\linewidth]{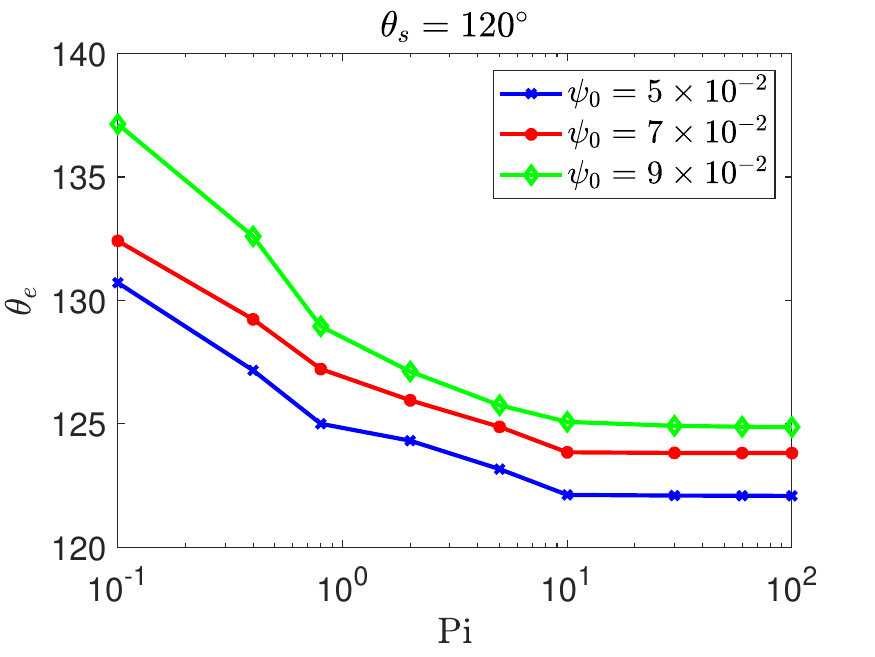}
	\includegraphics[width=0.35\linewidth]{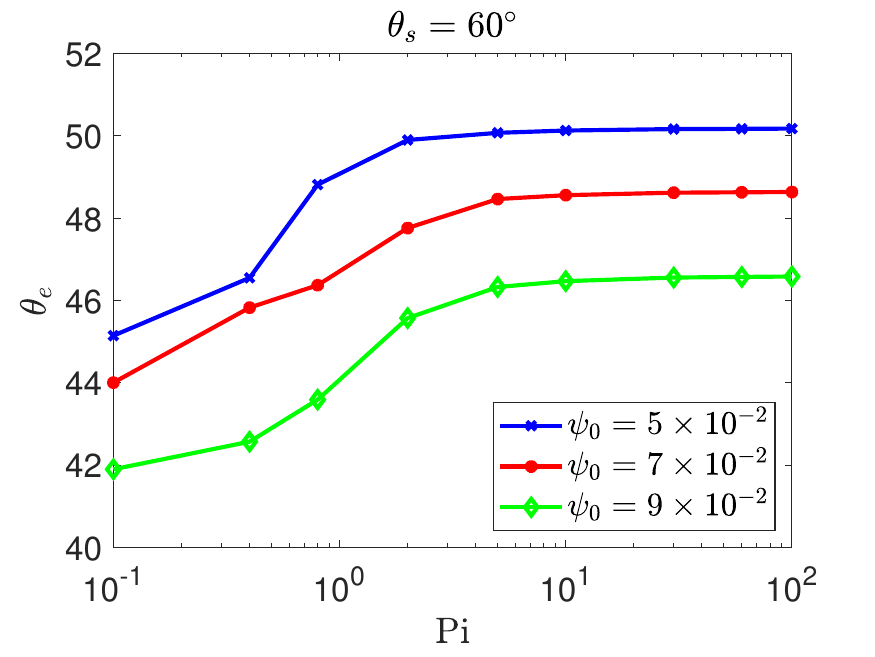}
	\caption{Plot of contact angle versus $\mathrm{Pi}$ with different initial surfactant concentrations $\psi_{0}$.}
	\label{fig:DiffPipsiAngle}
\end{figure}

The influence of different surfactant concentrations $\psi_0$ and the temperature-dependent surfactant diffusion rate $\mathrm{Pi}$ (with fixed inverse adsorption rate $\mathrm{Ex}=1$) on contact angle is further demonstrated in Fig.~\ref{fig:DiffPipsiAngle}. It can be seen that the change in the equilibrium contact angle is more pronounced as the surfactant concentration increases (with fixed $\mathrm{Pi}$). Meanwhile, the equilibrium contact angle for the hydrophobic/hydrophilic droplet increases/decreases with decreasing $\mathrm{Pi}$ with same amount of surfactants (i.e., fixed $\psi_0$). Given that the surfactant can adhere more tightly to the droplet interface with weaker diffusion (i.e., smaller $\mathrm{Pi}$), this result indicates that it is the distribution of surfactants concentrated along the droplet interface, rather than the overall concentration, controls the degree in which the contact angle changes. Since the surfactant distribution is determined by the competition of diffusion and adsorption, we expect that the contact angle also depends on the inverse adsorption rate $\mathrm{Ex}$ \cite{Wang2022PFSMCLModel}.

\subsection{Diffusion experiments concerning surfactants}\label{sec:5.3}
While we have studied the influence of surfactants on the equilibrium shape and contact angle of droplets in the previous section, we will investigate the controllability of the wetting and dewetting dynamics of droplets by tuning the concentration of surfactants.
\subsubsection{Wetting dynamics of droplets with surfactants}
\begin{figure}[htbp]
	\centering
	\subfigure[Initial profile]{
		\begin{minipage}{0.31\textwidth}
			\includegraphics[scale=0.35]{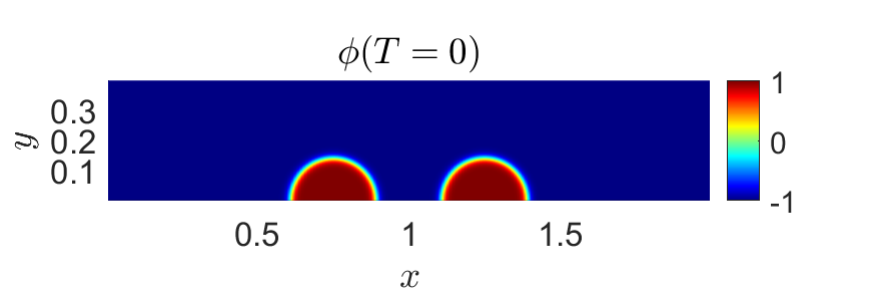}
			\includegraphics[scale=0.35]{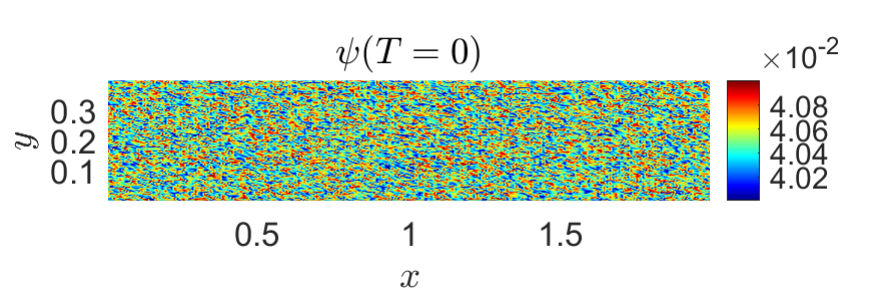}
		\end{minipage}
	}
	\subfigure[Evolution snapshots for $\langle\psi_{0}\rangle=0.04$]{
		\begin{minipage}{0.31\textwidth}
			\includegraphics[scale=0.35]{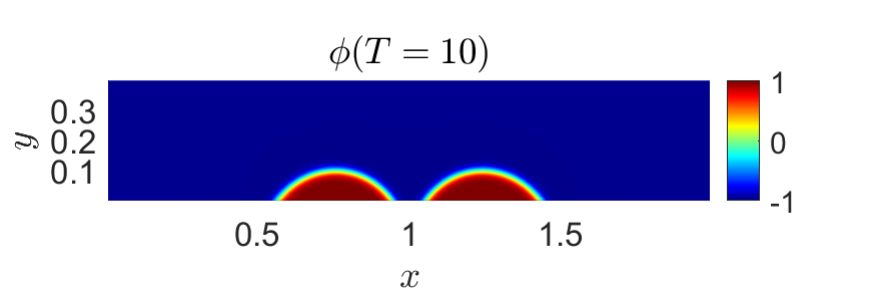}
			\includegraphics[scale=0.35]{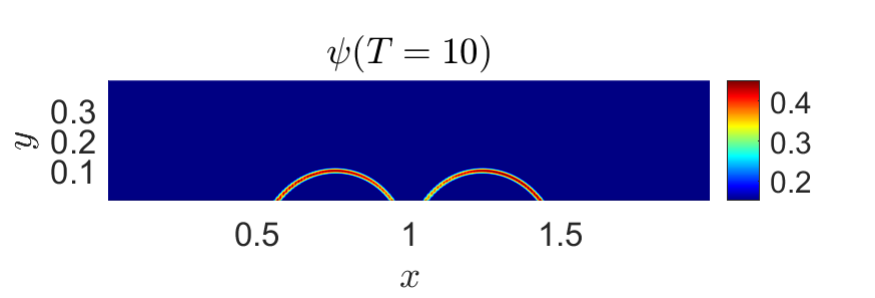}
		\end{minipage}
		\begin{minipage}{0.31\textwidth}
			\includegraphics[scale=0.35]{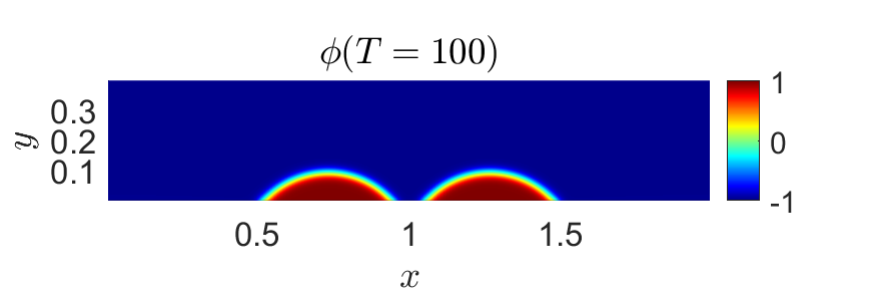}
			\includegraphics[scale=0.35]{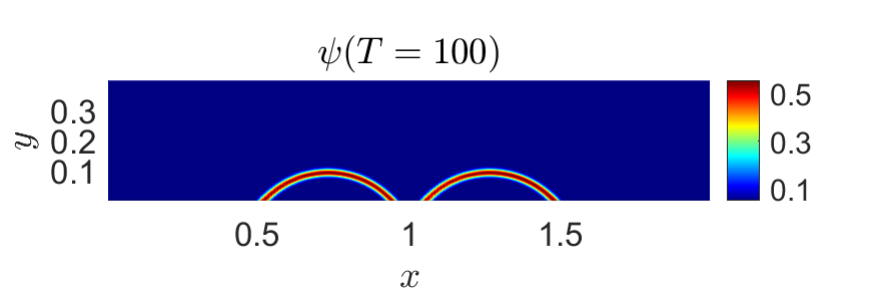}
		\end{minipage}
	}\\
        \subfigure[Energy dissipation]{
		\begin{minipage}{0.31\textwidth}
			\includegraphics[scale=0.35]{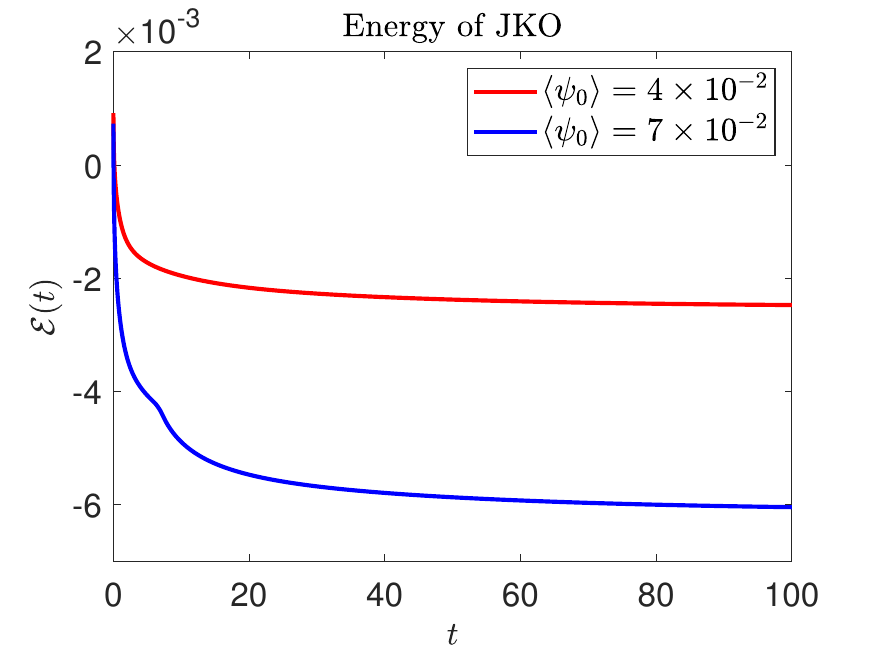}
		\end{minipage}
	}
        \subfigure[Evolution snapshots for $\langle\psi_{0}\rangle=0.07$]{
		\begin{minipage}{0.31\textwidth}
			\includegraphics[scale=0.35]{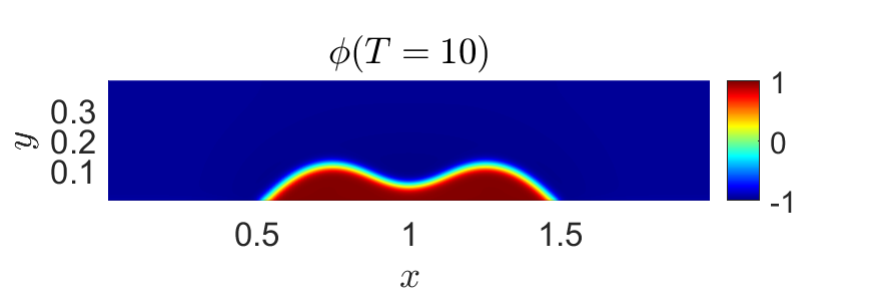}
			\includegraphics[scale=0.35]{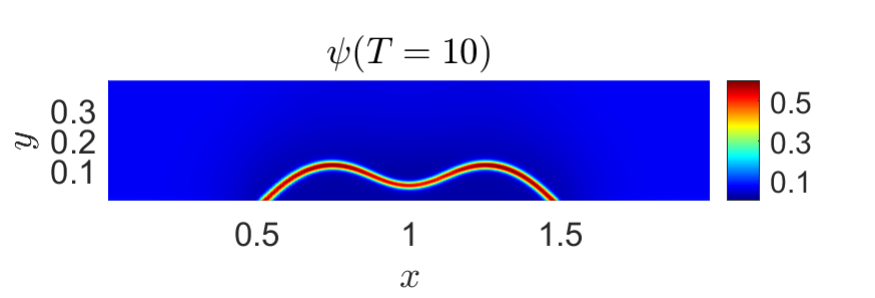}
		\end{minipage}
		\begin{minipage}{0.31\textwidth}
			\includegraphics[scale=0.35]{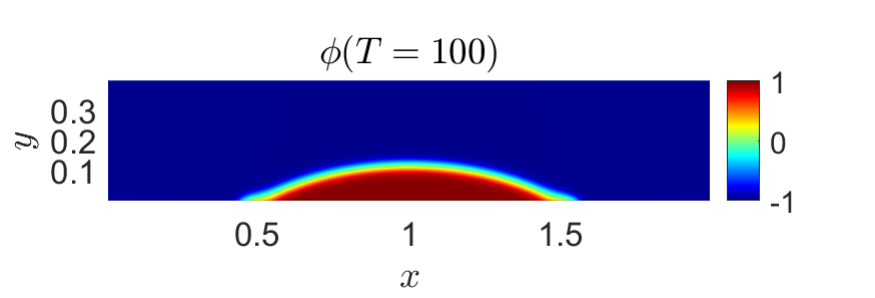}
			\includegraphics[scale=0.35]{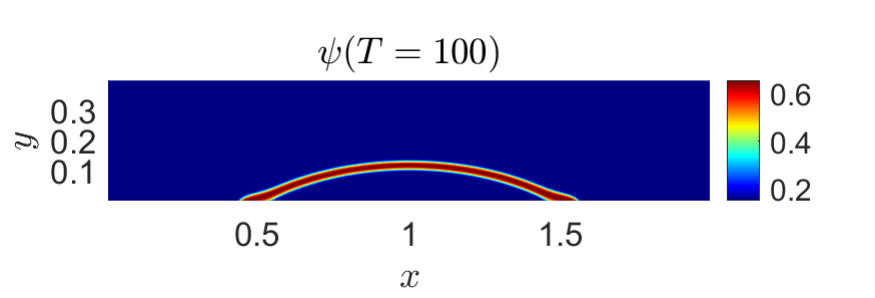}
		\end{minipage}
	}
	\caption{Evolution snapshots for the wetting dynamics of two droplets with $\langle\psi_{0}\rangle=0.04$ and $\langle\psi_{0}\rangle=0.07$ and their energy evolution.}
	\label{fig:Twodroplet004}
\end{figure}

We first consider the wetting of two droplets in the computational domain $\left[0,2\right]\times\left[0,0.4\right]$ with the following initial conditions (Fig.~\ref{fig:Twodroplet004} (a))
\begin{eqnarray}
	\begin{aligned}
		&\phi_{0}(x,y)=1+\sum_{i=1}^{2}\mathrm{tanh}\Big(10-\dfrac{\sqrt{(x-x_{i})^{2}+y^{2}}}{\sqrt{2}\mathrm{Cn}}\Big),\quad\text{where $x_{1}=0.75$ and $x_{2}=1.25$}.\\
		&\psi_{0}(x,y)=\langle\psi_{0}\rangle+0.001\xi,
	\end{aligned}
\end{eqnarray}
where $\langle\psi_{0}\rangle$ is the averaged initial concentration. We set $\mathrm{Cn}=0.01$, $\theta_{s}=60^{\circ}$, $\Delta x=\Delta y=0.005$ and $\Delta t=0.01$.

The wetting dynamics of the two droplets with different initial surfactant concentrations $\langle\psi_{0}\rangle=0.04$ and $\langle\psi_{0}\rangle=0.07$ are shown in Fig.~\ref{fig:Twodroplet004} (b) and (d), respectively. For the case with $\langle\psi_{0}\rangle=0.04$, the two droplets naturally relax to their steady states with contact angles less than $60^{\circ}$, but they do not merge and remain distant (even though very close) from each other all the time. In contrast, for the case with $\langle\psi_{0}\rangle=0.07$, the contact angle is altered more with higher surfactant concentration so that two droplets come into contact and begin to merge around $T=10$, eventually forming a larger droplet. The two different dynamics can also be revealed from the energy evolution in Fig.~\ref{fig:Twodroplet004} (c), where an additional abrupt energy decay around $t = 10$ is clearly observed for $\langle\psi_{0}\rangle=0.07$, corresponding to the coalescence of two droplets, while the energy curve for $\langle\psi_{0}\rangle=0.04$ has only one energy decay.

\begin{figure}[htbp]
	\centering
	\subfigure[Initial profile]{
		\begin{minipage}{0.28\textwidth}
			\includegraphics[scale=0.35]{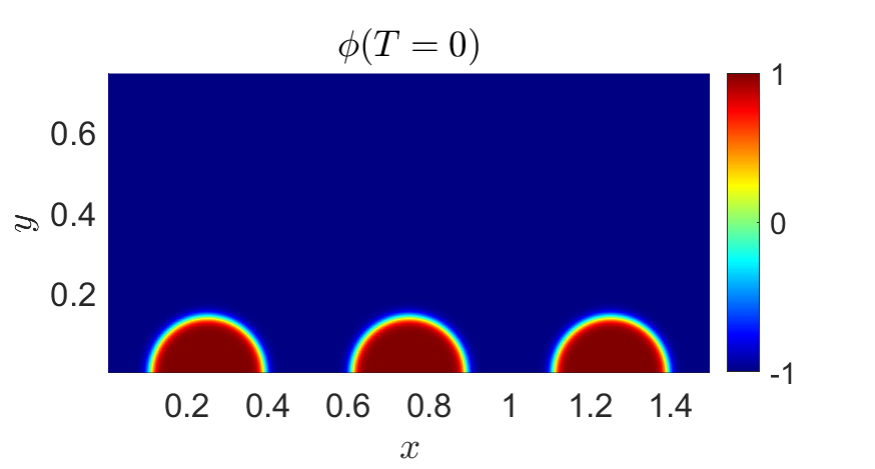}
		\end{minipage}
	}
        \subfigure[Evolution snapshots for $\langle\psi_{0}\rangle=0.02$]{
		\begin{minipage}{0.31\textwidth}
			\includegraphics[scale=0.35]{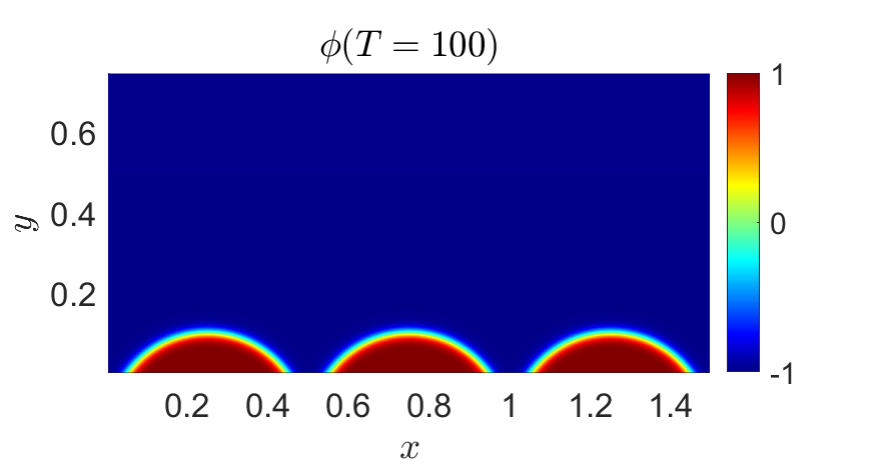}
		\end{minipage}
            \begin{minipage}{0.31\textwidth}
			\includegraphics[scale=0.35]{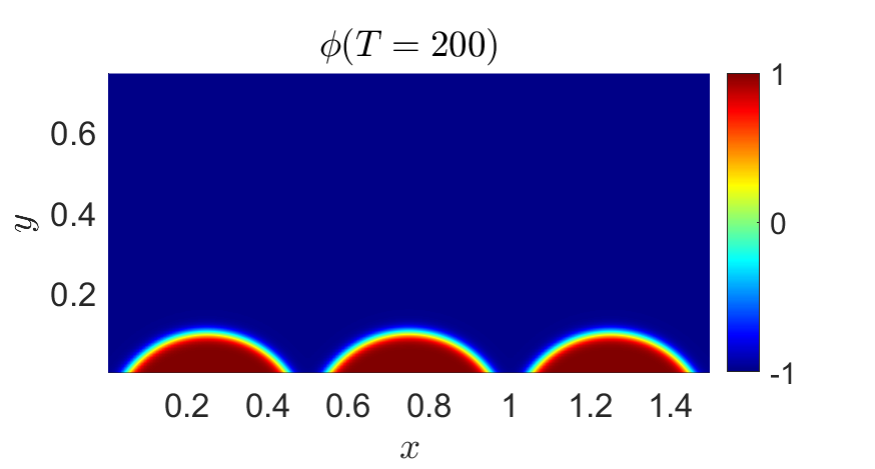}
		\end{minipage}
	}\\
        \subfigure[Energy dissipation]{
		\begin{minipage}{0.31\textwidth}
			\includegraphics[scale=0.35]{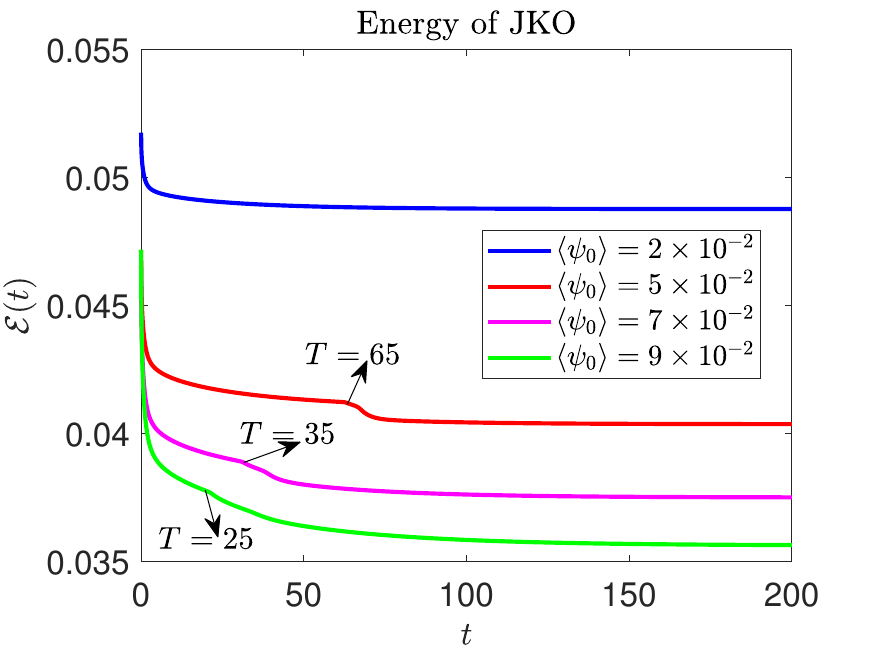}
		\end{minipage}
	}
	\subfigure[Evolution snapshots for $\langle\psi_{0}\rangle=0.05$]{
		\begin{minipage}{0.31\textwidth}
			\includegraphics[scale=0.35]{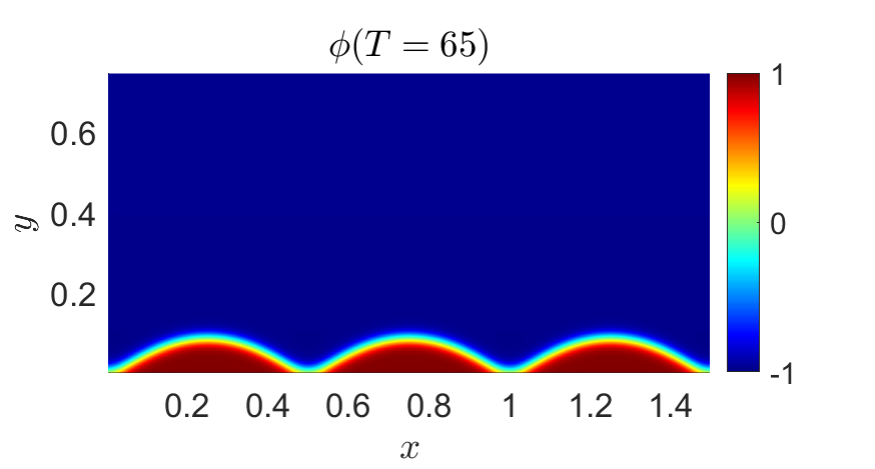}
		\end{minipage}
            \begin{minipage}{0.31\textwidth}
			\includegraphics[scale=0.35]{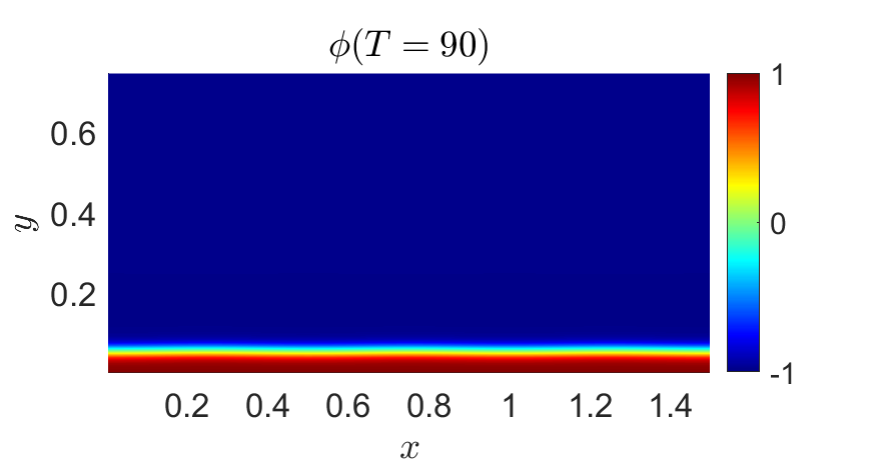}
		\end{minipage}
	}
	\caption{Evolution snapshots of $\phi$ for the wetting dynamics of three droplets with $\langle\psi_{0}\rangle=0.02$ and $\langle\psi_{0}\rangle=0.05$ and their energy dissipation.}
	\label{fig:Threedroplet}
\end{figure}

More comprehensive numerical experiments are performed for the wetting of three droplets, with initial conditions in computational domain $\left[0,1.5\right]\times\left[0,0.75\right]$ (Fig.~\ref{fig:Threedroplet} (a))
\begin{eqnarray}
	\begin{aligned}
		&\phi_{0}(x,y)=2+\sum_{i=1}^{3}\mathrm{tanh}\Big(10-\dfrac{\sqrt{(x-x_{i})^{2}+y^{2}}}{\sqrt{2}\mathrm{Cn}}\Big),\quad\text{where $x_{1}=0.25$, $x_{2}=0.75$, $x_{3}=1.25$},\\
		&\psi_{0}(x,y)=\langle\psi_{0}\rangle+0.001\xi.
	\end{aligned}
\end{eqnarray}

When surfactant concentration is relatively low as $\langle\psi_{0}\rangle=0.02$, three droplets remain separated during the wetting process, as shown in Fig.~\ref{fig:Threedroplet} (b). As we increase the surfactant concentration, the droplets begin to merge, as shown in Fig.~\ref{fig:Threedroplet} (d) for $\langle\psi_{0}\rangle=0.05$. Furthermore, we find that the amount of surfactant concentration not only determines the equilibrium droplet profile but also controls the dynamic process. More specifically, the timing of droplets coalescence and complete substrate wetting is advanced with higher surfactant concentrations, wich can be clearly seen in the energy curves in Fig.~\ref{fig:Threedroplet} (c). The time at which the second energy abrupt decay occurs is advanced from around $T=65$ for $\langle\psi_{0}\rangle=0.05$ to around $T=25$ for $\langle\psi_{0}\rangle=0.09$. 

\subsubsection{Dewetting dynamics of liquid thin film with surfactants}
Our last group of numerical experiments simulate the dewetting dynamics of a liquid film in the presence of surfactants. We consider the initial condition of $\phi$ as a plate-shaped thin film with length $l=2.5$ and thickness $h=0.03$ in the computational domain $\left[-1.5,1.5\right]\times\left[0,0.5\right]$. We choose $\mathrm{Cn}=0.006$, $\theta_s=120^{\circ}$ and $\Delta x=0.01, \Delta y=0.005$ in simulations and consider two uniform initial profiles $\psi_{0}(x,y)\equiv 0.07$ and $\psi_{0}(x,y)\equiv 0.05$ for surfactant concentrations.

\begin{figure}[htbp]
	\centering
	\subfigure{
		\begin{minipage}{0.31\textwidth}
			\includegraphics[scale=0.35]{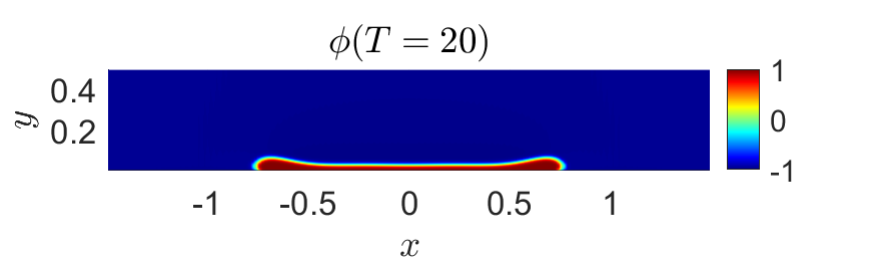}
			\includegraphics[scale=0.35]{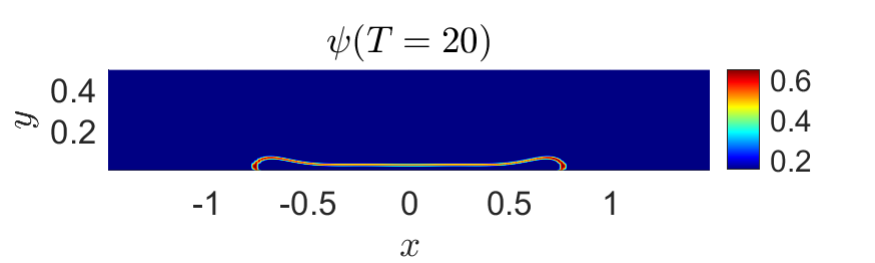}
		\end{minipage}
	}
	\subfigure{
		\begin{minipage}{0.31\textwidth}
			\includegraphics[scale=0.35]{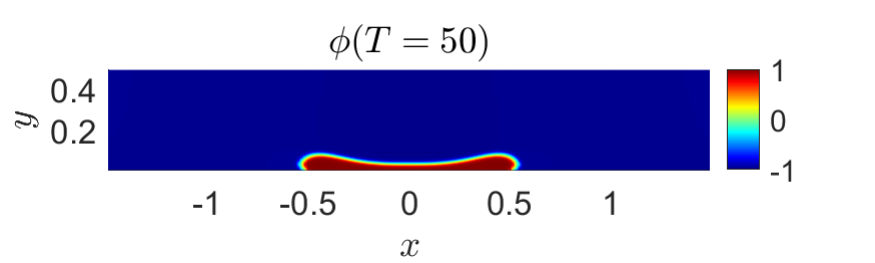}
			\includegraphics[scale=0.35]{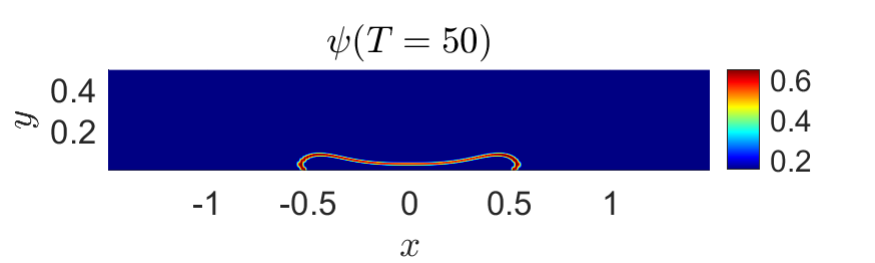}
		\end{minipage}
	}
	\subfigure{
		\begin{minipage}{0.31\textwidth}
			\includegraphics[scale=0.35]{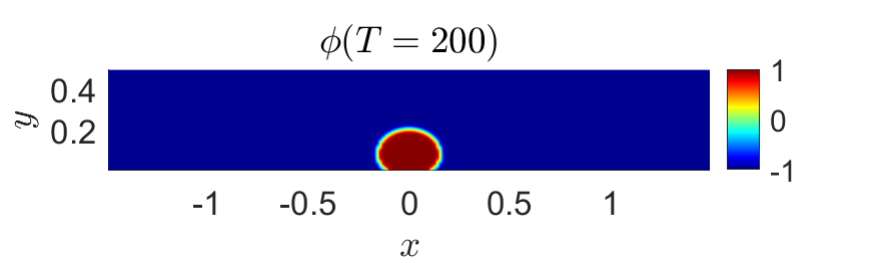}
			\includegraphics[scale=0.35]{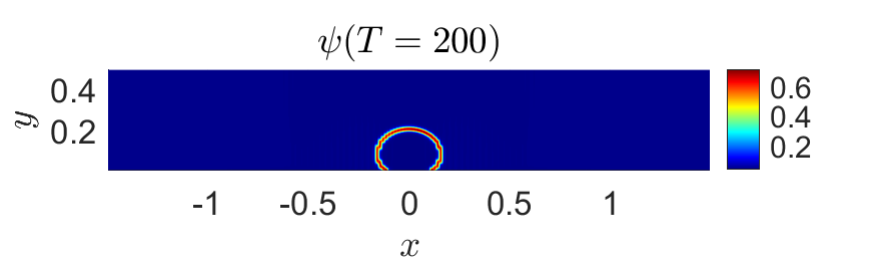}
		\end{minipage}
	}
	\caption{Evolution snapshots for the dewetting dynamics of thin film of height $h=0.03$ with $\psi_{0}=0.05$.}
	\label{fig:H003psi005}
\end{figure}

\begin{figure}[htbp]
	\centering
	\subfigure{
		\begin{minipage}{0.31\textwidth}
			\includegraphics[scale=0.35]{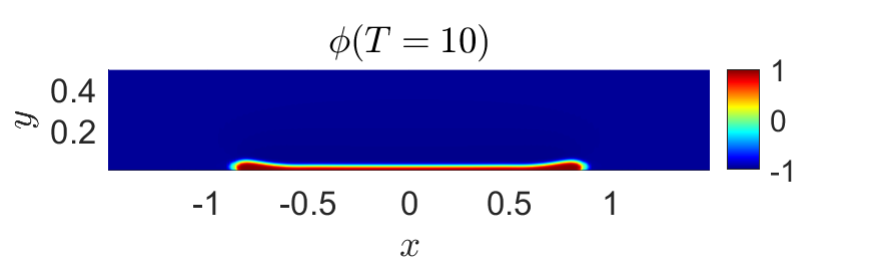}
			\includegraphics[scale=0.35]{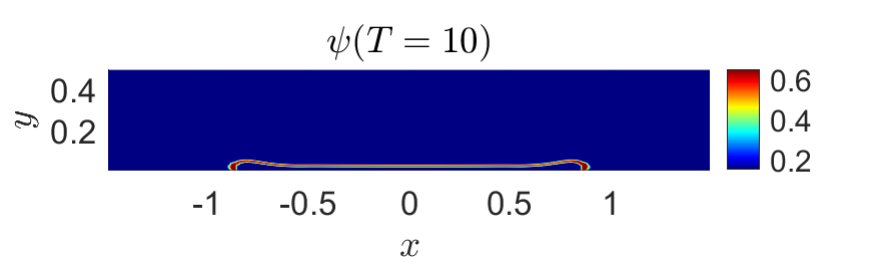}
		\end{minipage}
	}
	\subfigure{
		\begin{minipage}{0.31\textwidth}
			\includegraphics[scale=0.35]{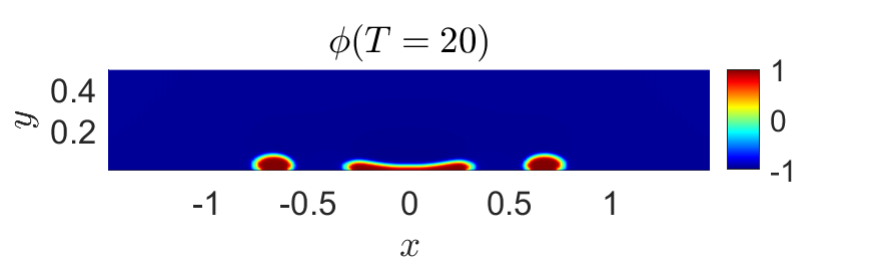}
			\includegraphics[scale=0.35]{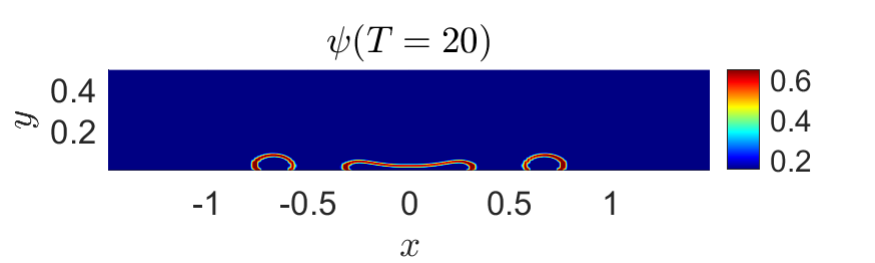}
		\end{minipage}
	}
	\subfigure{
		\begin{minipage}{0.31\textwidth}
			\includegraphics[scale=0.35]{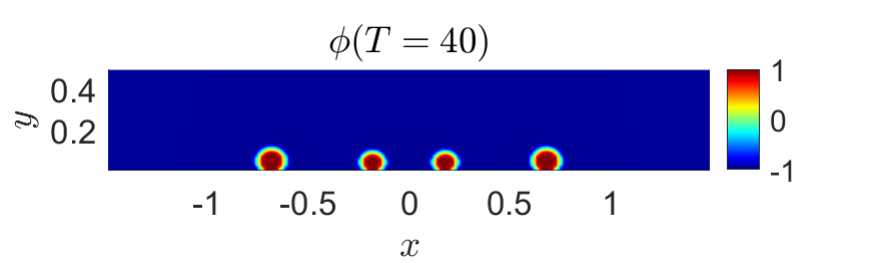}
			\includegraphics[scale=0.35]{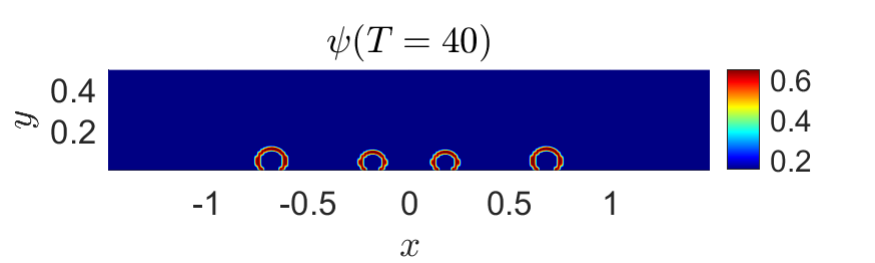}
		\end{minipage}
	}
    \includegraphics[width=0.35\linewidth]{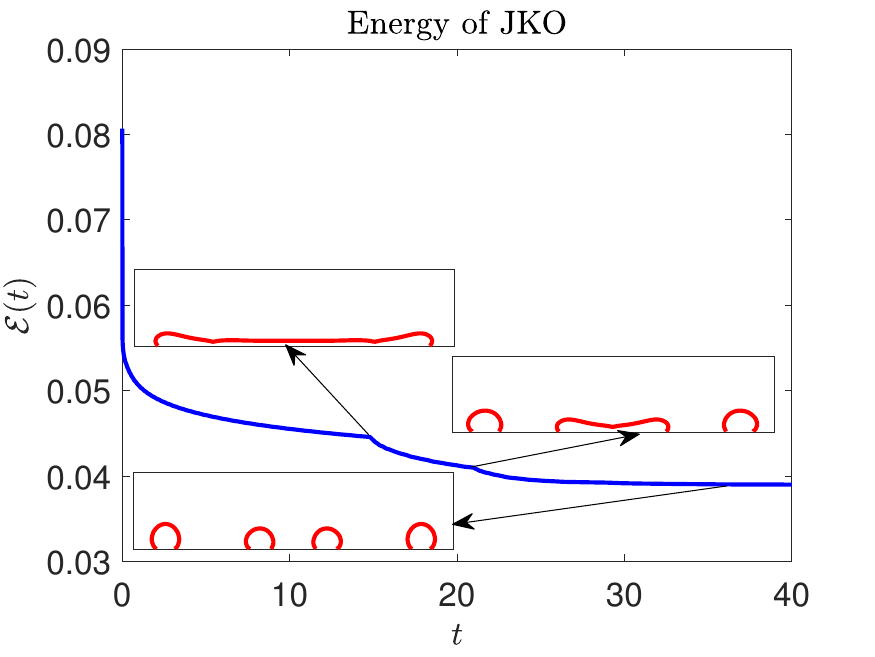}
	\caption{Evolution snapshots and energy dissipation for the dewetting dynamics of thin film of height $h=0.03$ with $\psi_{0}=0.07$.}
	\label{fig:H003psi007}
\end{figure}

Due to the unbalanced Young's force in \eqref{eq:natural_bc} at the contact line, the thin film starts to retract from its edge and the mass transport near the edge forms a ridge followed by a valley ($T=20$ in Fig.~\ref{fig:H003psi005}). As the edge retraction continues, we see two different dynamics for $\psi_{0}=0.05$ and $\psi_{0}=0.07$. When $\psi_0=0.05$, the thin film will slowly retract ($T=50$) and eventually becomes a large droplet ($T=200$), leading to the dewetting phenomenon. However, as the surfactant concentration increases, the thin film becomes more hydrophobic and undergoes a different dynamics (see Fig.~\ref{fig:H003psi007}). When $\psi_0=0.07$,  as the ridge grows and the valley deepens, the valley eventually touches the substrate and results in the pinch-off of the thin film ($T=20$). This pinch-off process will repeat once again, and the thin film eventually breaks up into four isolated small droplets ($T=40$). We also plot the energy evolution for $\psi_0=0.07$, which has abrupt decays around $T=15$ and $T=25$. The first energy abrupt decay corresponds to the first pinch-off process of the initial thin film, and the second one corresponds to the second pinch-off event of the remaining shorter film. 

\begin{figure}[htbp]
	\centering
	\subfigure{
		\begin{minipage}{0.31\textwidth}
			\includegraphics[scale=0.35]{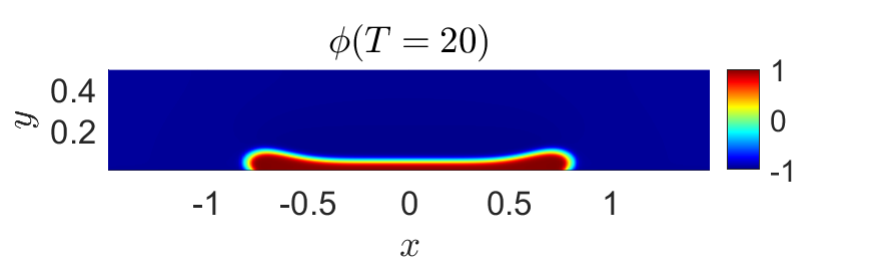}
			\includegraphics[scale=0.35]{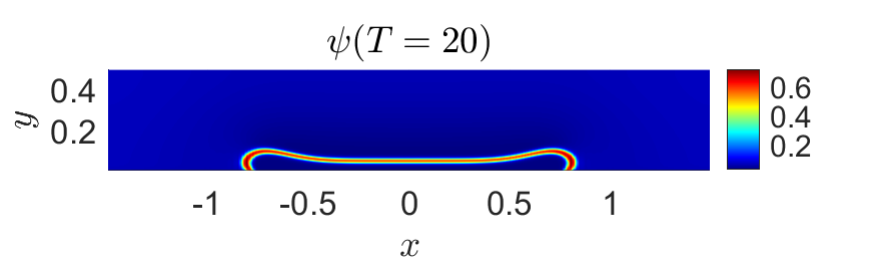}
		\end{minipage}
	}
	\subfigure{
		\begin{minipage}{0.31\textwidth}
			\includegraphics[scale=0.35]{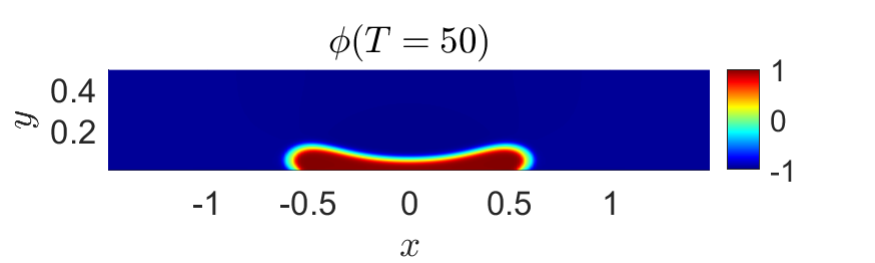}
			\includegraphics[scale=0.35]{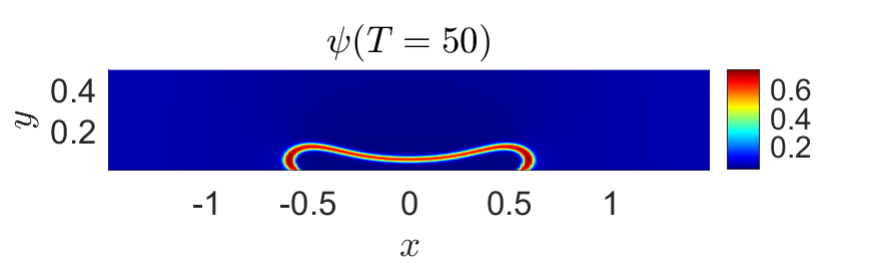}
		\end{minipage}
	}
	\subfigure{
		\begin{minipage}{0.31\textwidth}
			\includegraphics[scale=0.35]{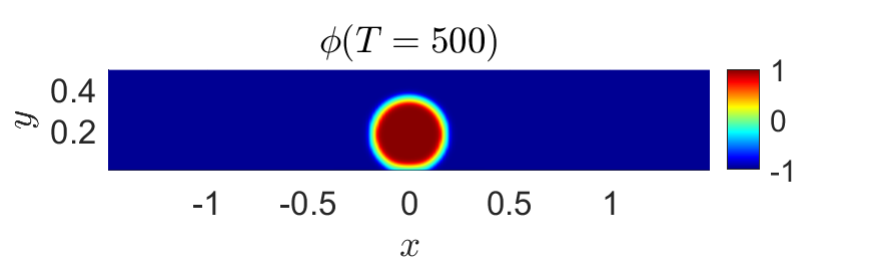}
			\includegraphics[scale=0.35]{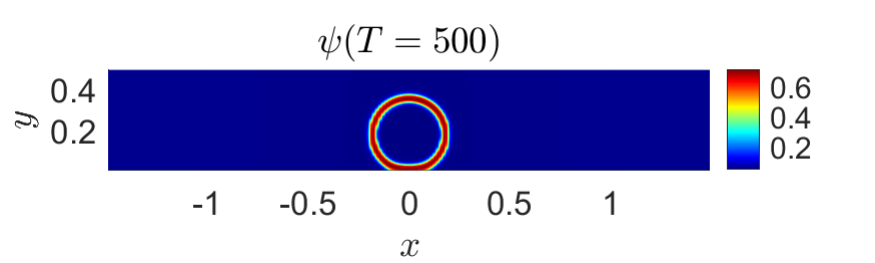}
		\end{minipage}
	}
	\caption{Evolution snapshots for the dewetting dynamics of thin film of height $h=0.05$ with $\psi_{0}=0.07$.}
	\label{fig:H005psi007}
\end{figure}

Finally, we investigate the dewetting dynamics of a thicker film of height $h=0.05$ with surfactant concentration $\psi=0.07$ (Fig.~\ref{fig:H005psi007}). We find that the pinch-off event does not occur during the whole dewetting process due to the increase of film thickness. More interestingly, when the thin film reaches the equilibrium shape of a droplet ($T=200$), we observe that the entire droplet is encapsulated by the surfactants and does not touch the substrate at all. This implies the effect of the presence of surfactants on the fluid-fluid interfacial tension \cite{Zhu2019ThermodynamicallyMCL}, which results in the inward flow of surfactants through the bottom of the droplet. This interesting phenomenon and its cause should be further confirmed and understood with the help of experimental studies. 

\section{Concluding remarks}
In this paper, we propose a novel structure-preserving scheme based on the variational formulation for the phase-field surfactant model with moving contact lines. Our proposed scheme is a generalization of the JKO scheme for Wasserstein gradient flows \cite{Carrillo2022PrimalDual,Carrillo2024StructurePD} to the coupled system of generalized Wasserstein gradient flows with diagonal matrix-valued mobility for phase-field and surfactants and $\mathcal{L}^{2}-$gradient flow for the dynamic boundary condition of moving contact lines. Unlike traditional methods, our scheme converts the problem of solving a nonlinear PDE system to a series of optimization problems with convex objective functions and linear PDE constraints, which can be efficiently solved by our primal-dual splitting method. Furthermore, the efficiency of our scheme can be significantly improved in practice by utilizing the preconditioning techniques with FFT-based fast algorithms and adaptive time stepping strategies. Owing to the well-posedness and variational structure of the JKO scheme, our method inherits the desirable features of solutions such as original energy dissipation, global bounds, and mass conservation at the discrete level. 

We validate the accuracy, performance, and structure-preserving properties of our numerical scheme through a series of numerical experiments. We also investigate the effects of surfactants on the contact angle, as well as the wetting and dewetting dynamics by performing several interesting simulations. 

We expect to conduct our future work in two directions concerning the generalization of the variational scheme and the acceleration of primal-dual algorithms. The present JKO scheme is only first-order accurate in time. Designing such second-order variational scheme for Wasserstein gradient flows that guarantees energy dissipation and bound-preserving is still an open problem \cite{Legendre2017SecondGradient,Matthes2019Variational,Carrillo2022PrimalDual}. Furthermore, the JKO scheme is expected to be extended for coupled gradient flows with symmetric positive definite mobility, e.g., the more sophisticated model for the hydrodynamics of thin film with insoluble surfactants \cite{Thiele2012Surfactant}. On the other hand, we can further improve the efficiency of our method by using adaptive mesh refinement in space \cite{Chen2016AdaptiveMesh, Wise2007Multigrid} and other acceleration techniques for the optimization algorithms, e.g., the transformed primal-dual method \cite{Chen2023Transformed} and the modified primal-dual method with a convergence rate independent of grid size \cite{Jacobs2019Gprox}. 

\section*{Declaration of competing interest}
The authors declare that they have no known competing financial interests or personal relationships that could have appeared to influence the work reported in this paper.

\section*{Acknowledgements}
CW is partially supported by the National Natural Science Foundation of China under grants 12371392 and 12431015. 
The work of ZZ is partially supported by National Key R\&D Program of China (2023YFA1011403), the NSFC grant (92470112), and the Shenzhen Sci-Tech Inno-Commission Fund (20231120102244002). 


\appendices
\section{Equilibrium boundary condition}\label{appd:eqbc}

As the mobility of contact line $\mathrm{Pe}_{s}\rightarrow 0$, the dynamic boundary equation reduces to the equilibrium boundary condition \eqref{eq:eq_bc}
\begin{equation}\label{eq:eq_bc2}
    \mathrm{Cn}\nabla\phi\cdot\boldsymbol{\nu}=-\gamma^{'}_{\omega f}(\phi),\qquad \text{on $\Gamma$}.
\end{equation}
When considering the equilibrium boundary condition, the boundary values $(\phi_{bc})_i=\phi_{i,\frac{1}{2}}$ will not be considered as independent variables of the minimization of the JKO scheme, but determined in advance by the boundary conditions, which will be used in the calculations of $\mathcal{E}^h$ and $\nabla \mathcal{E}^h$. 

Following the treatment of equilibrium boundary condition for a cubic wall energy in \cite{Carrillo2024StructurePD}, here we consider a different wall energy
\begin{eqnarray}
	\begin{aligned}
		&\gamma_{\omega f}(\phi)=-\dfrac{\sqrt{2}}{3}\mathrm{cos}(\theta_{s})\mathrm{sin}(\dfrac{\pi\phi}{2})+\dfrac{\gamma_{1}+\gamma_{2}}{2}\quad\text{with}\quad \mathrm{cos}({\theta_{s}})=\dfrac{3\sqrt{2}(\gamma_{2}-\gamma_{1})}{4}. 
	\end{aligned}
\end{eqnarray}
Discretizing the equilibrium boundary condition \eqref{eq:eq_bc2} yields
\begin{eqnarray}
	\begin{aligned}
		\mathrm{Cn}\Big(\dfrac{\phi_{i,1}-\phi_{i,0}}{\Delta y}\Big)=\gamma^{'}_{\omega f}(\phi_{i,\frac{1}{2}})\quad\text{where}\quad\phi_{i,0}=2\phi_{i,\frac{1}{2}}-\phi_{i,1},
	\end{aligned}
\end{eqnarray}	
which reduces to an equation for $X=\phi_{i,\frac{1}{2}}$
\begin{eqnarray}\label{bc1root1}
	\begin{aligned}
		f(X)=\phi_{i,1}-X-\alpha\mathrm{cos}\Big(\dfrac{\pi X}{2}\Big)=0, \quad \text{where $\alpha=-\frac{\sqrt{2}\pi\mathrm{cos}\theta_{s}\Delta y}{12\mathrm{Cn}}$.}
	\end{aligned}
\end{eqnarray}

\begin{figure}[htbp]
	\centering  
	\includegraphics[width=0.35\textwidth]{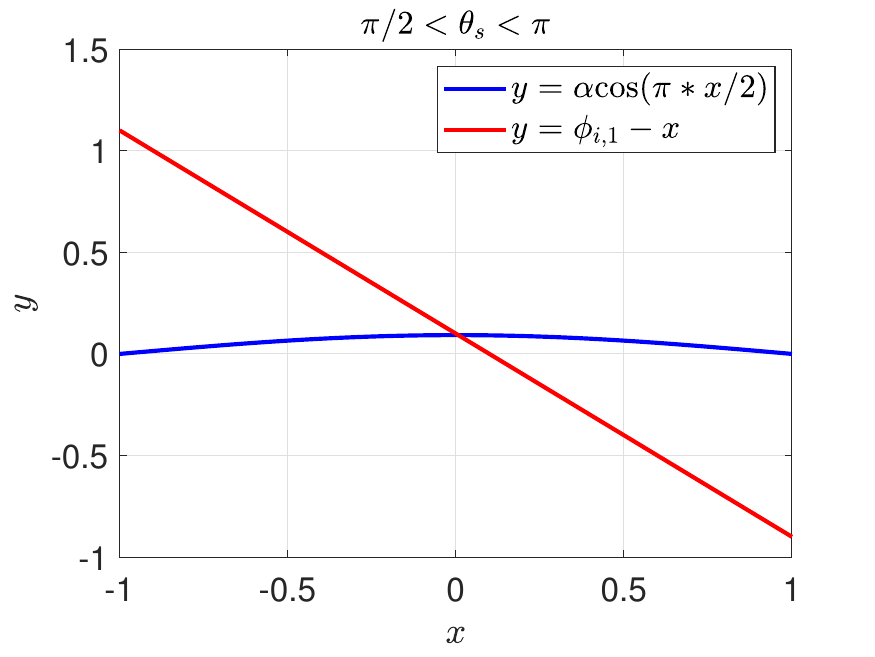}
	\includegraphics[width=0.35\textwidth]{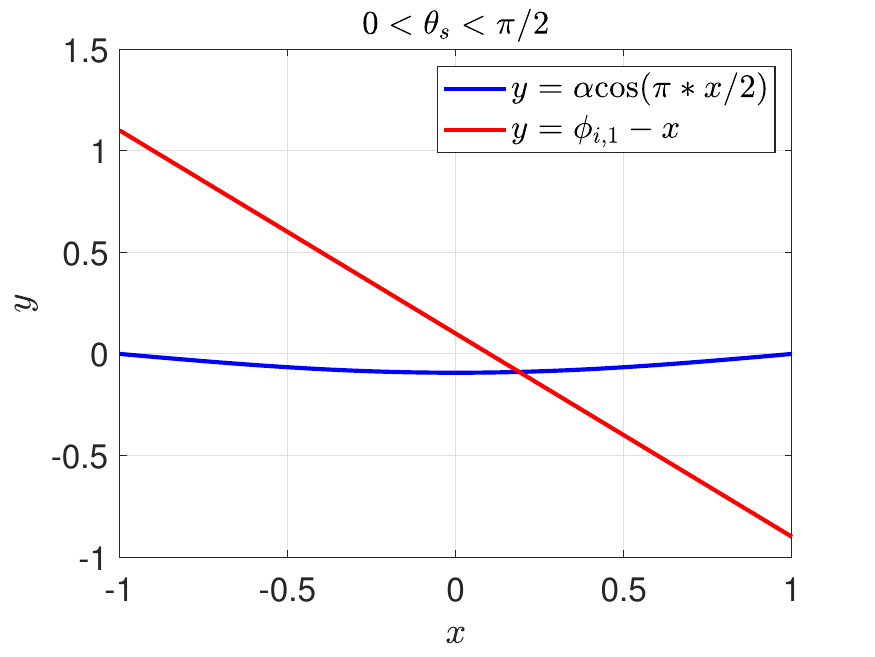}
	\caption{Intersection of curves of \eqref{bc1root1}. Left: $\alpha<0$. Right: $\alpha>0$.}
	\label{fig:RootVirtual}
\end{figure}

From Fig.~\ref{fig:RootVirtual}, it can be observed that the root of $f(X)$ uniquely exists (which can be proved later) and we can use Newton's method to find the root $X^*=\phi_{i,\frac{1}{2}}$. To choose appropriate initial guesses, we look at the derivatives of $f(X)$:
\begin{equation}
		f^{'}(X)=-1+\dfrac{\pi\alpha}{2}\mathrm{sin}\Big(\dfrac{\pi X}{2}\Big),\quad
		f^{''}(X)=\dfrac{\alpha\pi^{2}}{4}\mathrm{cos}\Big(\dfrac{\pi X}{2}\Big).
\end{equation}
Given that $\Delta y<\mathrm{Cn}$ in practice, we can show that 

(i) $f^{'}(X)<0$ for all $X$, which confirms the uniqueness of the root;

(ii) $f(-|\phi_{i,1}|-|\alpha|)\geq 0$ and $f(|\phi_{i,1}|+|\alpha|)\leq 0$, indicating $-|\phi_{i,1}|-|\alpha|\leq X^*\leq |\phi_{i,1}|+|\alpha|$.

Then we have the following strategy for choosing initial guesses to guarantee the convergence of Newton iterations:
\begin{description}
	\item[Case 1:] When $-1\leq\phi_{i,1}\leq 1$ and $0<\theta_{s}<\dfrac{\pi}{2}$, we have $f(-1)\geq 0$, $f(1)\leq 0$ and hence $X^*\in[-1,1]$. Since $f^{''}(X)<0$ on $[-1,1]$, we set the initial guess $\tilde{X}_0=1$.
	\item[Case 2:] When $-1\leq\phi_{i,1}\leq 1$ and $\dfrac{\pi}{2}\leq \theta_{s}<\pi$, we have $f^{''}(X)\geq 0$ on $[-1,1]$, then we set the initial guess $\tilde{X}_0=-1$.
	\item[Case 3:] When $\phi_{i,1}<-1$ and $0<\theta_{s}<\dfrac{\pi}{2}$, we have $X^*\in[-|\phi_{i,1}|-|\alpha|,-1)\subset (-2,-1)$ and hence $f^{''}(X)> 0$ on $[-1,1]$, then we set the initial guess $\tilde{X}_0=-|\phi_{i,1}|-|\alpha|$.
	\item[Case 4:] When $\phi_{i,1}<-1$ and $\dfrac{\pi}{2}\leq \theta_{s}<\pi$, we have $X^*\in[-|\phi_{i,1}|-|\alpha|,-1)$ and hence $f^{''}(X)\leq 0$ on $[-1,1]$, then we set the initial guess $\tilde{X}_0=-1$.
    \item[Case 5:] When $\phi_{i,1}>1$ and $0<\theta_{s}<\dfrac{\pi}{2}$, we have $X^*\in(1,|\phi_{i,1}|+|\alpha|\big]\subset (1,2)$ and hence $f^{''}(X)> 0$ on $[-1,1]$, then we set the initial guess $\tilde{X}_0=1$.
    \item[Case 6:] When $\phi_{i,1}>1$ and $\dfrac{\pi}{2}\leq \theta_{s}<\pi$, we have $X^*\in(1,|\phi_{i,1}|+|\alpha|\big]$ and hence $f^{''}(X)\leq 0$ on $[-1,1]$, then we set the initial guess $\tilde{X}_0=|\phi_{i,1}|+|\alpha|$.  
\end{description}	

\section{Matrix for the discretization of continuity equation}\label{appd:disA}

We provide the construction of matrix $A$ in \eqref{eq:minproblem} for the discrete continuity equation and boundary conditions in Problems~2. The matrix $A\in\mathbb{R}^{2N\times(6N+N_{x})}$ is in the form below
\begin{eqnarray}\label{eq:A_matix}
	\begin{aligned}
		A=\left[\begin{array}{c:c:c:c:c:c:c}
			\mathcal{I}_{N} & \mathcal{D}_{x} & \mathcal{D}_y & \mathbf{0}_{N\times N} & \mathbf{0}_{N\times N} & \mathbf{0}_{N\times N} & \mathbf{0}_{N\times N_{x}}\\
			\hdashline
			\mathbf{0}_{N\times N} & \mathbf{0}_{N\times N} & \mathbf{0}_{N\times N} & \mathcal{I}_{N} & \mathcal{D}_x & \mathcal{D}_y & \mathbf{0}_{N\times N_{x}}\\
		\end{array}\right]
	\end{aligned}
\end{eqnarray}
where $\mathcal{I}_{N}\in\mathbb{R}^{N\times N}$ represents the unit matrix and $\mathcal{D}_x$, $\mathcal{D}_y$ represent the discrete differential operator as follows
\begin{eqnarray}
	\begin{aligned}
		\mathcal{D}_x=\dfrac{1}{2\Delta x}\mathcal{I}_{N_y}\otimes\mathcal{D}_{N_x},\quad 
        \mathcal{D}_y=\dfrac{1}{2\Delta y}\mathcal{D}_{N_y}\otimes\mathcal{I}_{N_x},
	\end{aligned}
\end{eqnarray}
where $\otimes$ denote the Kronecker product and the matrix $\mathcal{D}_{n}\in\mathbb{R}^{n\times n}$ for finite differencing is defined as below
\begin{eqnarray}
	\begin{aligned}
	\mathcal{D}_{n}=
        \begin{bmatrix} 
			1 & 1 & 0 &\cdots & 0 & 0 & 0\\
			-1 & 0 & 1 &\cdots & 0 & 0 & 0\\
			0 & -1 & 0 &\cdots & 0 & 0 & 0\\
			\vdots & \vdots & \vdots & \ddots & \vdots & \vdots & \vdots\\
			0 & 0 & 0 &\cdots & 0 & 1 & 0\\
			0 & 0 & 0 &\cdots & -1 & 0 & 1\\
			0 & 0 & 0 &\cdots & 0 & -1 & -1\\
		\end{bmatrix}.
	\end{aligned}
\end{eqnarray}

\section{Proof of Structure-preserving Themorem \ref{EMB}} \label{sec:proof}

Proof: (i) Since $u^{k+1}=(\phi^{k+1},\textbf{m}_{\phi}^{k+1};\psi^{k+1},\textbf{m}_{\psi}^{k+1};\phi_{bc}^{k+1})$ solves the minimization problem \eqref{eq:fullJKO}, we have the following inequality holds for any $k\geq 0$
\begin{equation}\label{eq:eninequality} 
    \sum_{\rho=\phi,\psi}\sum_{i,j}\mathcal{D}(\rho^{k+1}_{i,j},{\textbf{m}_{\rho}}^{k+1}_{i,j})\Delta x\Delta y+\dfrac{\mathrm{Pe}_{s}}{2}\sum_{i}((\phi_{bc})_i^{k+1}-(\phi^k_{bc})_i)^{2}\Delta x 
    +\Delta t\mathcal{E}(\phi^{k+1},\psi^{k+1},\phi_{bc}^{k+1})
    \leq \Delta t\mathcal{E}(\phi^{k},\psi^{k},\phi_{bc}^{k}),
\end{equation}
    which is in fact the original energy dissipation law \eqref{eq:dE_dt} at the discrete level. Due to the non-negativeness of the distance terms, we obtain the energy dissipation
\begin{align}
    \mathcal{E}(\phi^{k+1},\psi^{k+1},\phi_{bc}^{k+1})\leq\mathcal{E}(\phi^{k},\psi^{k},\phi_{bc}^{k}).
\end{align}	

(ii) For the sake of convenience, here we only consider the mass conservation of $\phi$. Summing both ends of the constraint on $\phi$ in \eqref{eq:fullJKO} over $i$ and $j$ yields
	\begin{align}
		&\sum_{i=1}^{N_x}\sum_{j=1}^{N_y}\Big(\phi_{i,j}+\dfrac{1}{2\Delta x}((m^{x}_{\phi})_{i+1,j}-(m^{x}_{\phi})_{i-1,j})+\dfrac{1}{2\Delta y}((m^{y}_{\phi})_{i,j+1}-(m^{y}_{\phi})_{i,j-1})\Big)=\sum_{i=1}^{N_x}\sum_{j=1}^{N_y}\phi^{k}_{i,j},
	\end{align}
where the summation of the flux variables $m^{x}_{\phi}$ can be written as follows	
	\begin{eqnarray}\label{eq:discrete_m}
		\begin{aligned}
			\sum_{i=1}^{N_{x}}((m^{x}_{\phi})_{i+1,j}-(m^{x}_{\phi})_{i-1,j})
            &=-(m^{x}_{\phi})_{0,j}+(m^{x}_{\phi})_{2,j}-(m^{x}_{\phi})_{1,j}+(m^{x}_{\phi})_{3,j}\cdots-(m^{x}_{\phi})_{N_{x}-1,j}+(m^{x}_{\phi})_{N_{x}+1,j} \\
            &=-((m^{x}_{\phi})_{0,j}+(m^{x}_{\phi})_{1,j})+((m^{x}_{\phi})_{N_{x},j}+(m^{x}_{\phi})_{N_{x}+1,j})=0, \quad \text{for each $j$,} 
		\end{aligned}
	\end{eqnarray}
where we have applied the non-flux boundary conditions in \eqref{eq:fullJKO}
	\begin{align}
		(m^{x}_{\phi})_{0,j}=-(m^{x}_{\phi})_{1,j},(m^{x}_{\phi})_{N_{x}+1,j}=-(m^{x}_{\phi})_{N_{x},j}.
	\end{align}
Similarly, we can have
\begin{align}
    \sum_{j=1}^{N_{y}}((m^{y}_{\phi})_{i,j+1}-(m^{y}_{\phi})_{i,j-1})=0, \quad \text{for each $i$.}
\end{align}
Thus, we obtain the mass conservation of the field variable $\phi$
\begin{align}
    \sum_{i=1}^{N_x}\sum_{j=1}^{N_y} \phi_{i,j} =\sum_{i=1}^{N_x}\sum_{j=1}^{N_y} \phi^{k}_{i,j}.
\end{align}
The same is true for surfactant concentration $\psi$.

(iii) According to the definition of weighted Wasserstein distance \eqref{D_function}, the minimizer of the JKO scheme must lie within the range where $M_\psi(\psi)\geq 0$ such that $\mathcal{D}(\psi,m_\psi)$ is finite. Furthermore, this bound-preserving property can be guaranteed in the computation of the proximal operator of $\mathcal{D}(\psi,m_\psi)$ by Newton's method with the proposed strategy (see Remark~\autoref{rmk:primal_prox} for details).\qed

\section{The spectrum and FFT-based matrix inversion of $C_2=\lambda AA^{\mathrm{T}}$}\label{appd:inversion}

From Remark~\ref{rmk:dual_prox1} and \ref{rmk:dual_prox2}, we need to use FFT-based fast algorithm to compute the proximal operator $\mathrm{Prox}^{C_{2}}_{i^{*}_{\delta}}$ either exactly or inexactly. Given that the matrix $AA^{\mathrm{T}}$ corresponds to the discrete Laplacian operator with homogeneous Neumann boundary with staggered grids in two dimensions, the corresponding discrete Fourier transform and its eigenvalue are given in Table~\ref{tal:FFTInformation}
\begin{figure}[htbp]
	\centering
	\begin{minipage}[c][0.1\textheight][c]{1\textwidth}
		\centering
		\tabcaption{Basic information about the FFT in the two-dimensional case.}
		\label{tal:FFTInformation}
		\begin{tabular}{c|c|c|c|c}\hline
			Boundary conditions & Grid & Forward & Backward & Eigenvalue\\ 
			\hline
			Neumann & Staggered & DCT-II & DCT-III & $\lambda_{mn}=\lambda^{x}_{m}+\lambda^{y}_{n}+1$\\
			\hline
		\end{tabular}
	\end{minipage}
\end{figure}

The discrete cosine transform DCT-II and DCT-III in two dimensions are defined as follows
\begin{eqnarray}\label{eq:dct}
	\begin{aligned}	
    &\hat{f}_{mn}=\frac{2}{\sqrt{N}}\sum_{i=1}^{N_x}\sum_{j=1}^{N_y}\sqrt{\frac{1}{(1+\delta_{1m})(1+\delta_{1n})}}f_{ij}\mathrm{cos}\Big(\dfrac{\pi(i-\frac{1}{2})(m-1)}{N_x}\Big)\mathrm{cos}\Big(\dfrac{\pi(j-\frac{1}{2})(n-1)}{N_y}\Big),\\	
    &f_{ij}=\frac{2}{\sqrt{N}}\sum_{m=1}^{N_x}\sum_{n=1}^{N_y}\sqrt{\frac{1}{(1+\delta_{1m})(1+\delta_{1n})}}\hat{f}_{mn}\mathrm{cos}\Big(\dfrac{\pi(m-1)(i-\frac{1}{2})}{N_x}\Big)\mathrm{cos}\Big(\dfrac{\pi(n-1)(j-\frac{1}{2})}{N_y}\Big),
	\end{aligned}
\end{eqnarray}
where $\delta_{ij}$ is the Kronecker delta. The eigenvalues $\lambda^{x}_{m}$ and $\lambda^{y}_{n}$ are given by
\begin{eqnarray}
	\begin{aligned}
		&\lambda^{x}_{m}=\dfrac{1}{2\Delta x^{2}}\Big(1-\mathrm{cos}\Big(\dfrac{2\pi(m-1)}{N_x}\Big)\Big),\\
		&\lambda^{y}_{n}=\dfrac{1}{2\Delta y^{2}}\Big(1-\mathrm{cos}\Big(\dfrac{2\pi(n-1)}{N_y}\Big)\Big).
	\end{aligned}
\end{eqnarray}

For the matrix inversion in the form $AA^{\mathrm{T}}u=b$, we first use the forward transform to obtain $\hat{b}=\mathrm{dct}(b)$, then use the inverse transform to obtain $u=\mathrm{idct}(\hat{b}/\lambda)$. We also compare the speed of different methods for matrix inversion through Table~\ref{tal:FFTCPUTime}.

\begin{figure}[htbp]
	\centering
	\begin{minipage}[c][0.1\textheight][c]{1\textwidth}
		\centering
		\tabcaption{Time comparison of different methods for solving matrix inversion.}
		\label{tal:FFTCPUTime}
		\begin{tabular}{c|c|c|c}\hline
			Method & $u^{*}=\mathrm{inv}(AA^{\mathrm{T}})b$ & $u^{*}=AA^{\mathrm{T}}\backslash b$ & FFT \\ 
			\hline
			CPU Time & 83.995 & 0.812 & 0.0275 \\
			\hline
		\end{tabular}
	\end{minipage}
\end{figure}

\section{Energy stable schemes}\label{appd:schemes}
	\begin{itemize} 
		\item Convex splitting scheme (CSS): 
		\begin{eqnarray}
			\begin{aligned}
				&\dfrac{\phi^{n+1}-\phi^{n}}{\Delta t}=\Delta\mu^{n+1},\\
				&\mu^{n+1}=-\mathrm{Cn}\Delta \phi^{n+1}+\dfrac{1}{\mathrm{Cn}}\Big(f_{+}(\phi^{n+1})-f_{-}(\phi^{n})\Big).
			\end{aligned}
		\end{eqnarray}
		$f_{+}(\phi)=F^{'}_{+}(\phi)=\phi^{3}$ and $f_{-}(\phi)=F^{'}_{-}(\phi)=\phi$. Thus, the discrete form is as follows
		\begin{align}
			\Big(I+\mathrm{Cn}\Delta tA_{h}^{2}-\dfrac{\Delta t}{\mathrm{Cn}}A_{h}D(\Phi^{n})^{2}\Big)\Phi^{n+1}=\Big(I-\dfrac{\Delta t}{\mathrm{Cn}}A_{h}\Big)\Phi^{n}.
		\end{align}
		where, $\Phi=\left\{\phi_{i,j}\right\}^{N^{2}}_{i,j=1}$. $D(\Phi)$ is a diagonal matrix with the vector $\Phi$ placed element-wise along the main diagonal. The matrix $A_{h}$ is the ordered discrete laplacian \cite{Eyre1998UnconditionallyGradient}.
		
		\item Stabilized semi-implicit scheme (SSI):
		\begin{eqnarray}
			\begin{aligned}
				&\dfrac{\phi^{n+1}-\phi^{n}}{\Delta t}=\Delta\mu^{n+1},\\
				&\mu^{n+1}=-\mathrm{Cn}\Delta \phi^{n+1}+\dfrac{1}{\mathrm{Cn}}f(\phi^{n})+S(\phi^{n+1}-\phi^{n}).
			\end{aligned}
		\end{eqnarray}

		The constant term here is $S(\phi^{n+1}-\phi^{n})$, and it introduces the truncation error $S\Delta t\Delta\phi_{t}(\xi_{n}), \xi_{n}\in(\phi^{n},\phi^{n+1})$. Therefore, the discrete form is as follows
		\begin{align}
			\Big(I+\mathrm{Cn}\Delta tA_{h}^{2}-\Delta tSA_{h}\Big)\Phi^{n+1}=(I-\Delta tSA_{h})\Phi^{n}+\dfrac{\Delta t}{\mathrm{Cn}}A_{h}f(\Phi^{n}).
		\end{align}
	\end{itemize}

\bibliographystyle{elsarticle-num.bst}
\bibliography{Manu_refer.bib}
\end{document}